\pdfoutput=1
\input epsf

\documentclass[11pt]{article}
\newcommand{\f}{\frac}

\usepackage{float}

\usepackage[dvips]{graphicx}
\usepackage{latexsym,amsmath,amsfonts,amscd,,amssymb,subfigure}
\usepackage{graphics}
\usepackage{graphicx}
\usepackage[dvips]{color}
\usepackage{epsfig}
\usepackage{changebar}
\usepackage{bm}
\numberwithin{equation}{section}

\begin{document}

\vspace{.5in}

\begin{center}

{\Large \bf Absolutely convergent fixed-point fast sweeping WENO methods for steady state of hyperbolic conservation laws}

\end{center}

\vspace{.15in}

\centerline{
Liang Li\footnote{
College of Science, Nanjing University of Aeronautics and Astronautics, Nanjing, Jiangsu 210016, P.R.
China. E-mail:liliangnuaa@163.com.},
Jun Zhu\footnote{
College of Science, Nanjing University of Aeronautics and Astronautics, Nanjing, Jiangsu 210016, P.R.
China. E-mail: zhujun@nuaa.edu.cn. Research was supported by NSFC grant 11872210 and Science Challenge Project, No. TZ2016002.},
Yong-Tao Zhang\footnote{Department of Applied and Computational
Mathematics and Statistics, University of Notre Dame, Notre Dame,
IN 46556, USA. E-mail: yzhang10@nd.edu. Research was partially supported by NSF grant DMS-1620108.}
}

 \begin{abstract}
Fixed-point iterative sweeping methods were developed in the literature to
efficiently solve steady state solutions of Hamilton-Jacobi equations and hyperbolic conservation laws. Similar as
other fast sweeping schemes, the key
components of this class of methods are the
Gauss-Seidel iterations and alternating sweeping strategy to achieve fast convergence
rate. A family of characteristics of the corresponding hyperbolic partial differential
equations (PDEs) are covered in a certain direction simultaneously in each sweeping order.
Furthermore, good properties of fixed-point iterative sweeping methods which are different from other fast sweeping methods include that they have explicit forms and do not
involve inverse operation of nonlinear local systems, and they can be applied
to general hyperbolic equations using any monotone numerical fluxes and
high order approximations easily. In a recent article [L. Wu, Y.-T. Zhang, S. Zhang and C.-W. Shu,
Commun. Comput. Phys., 20 (2016), 835-869], a fifth order fixed-point sweeping WENO scheme was designed for
solving steady state of hyperbolic conservation laws, and it was shown that the scheme converges to steady state solution much faster than the regular total variation diminishing (TVD) Runge-Kutta time-marching approach by
stability improvement of high order schemes with a forward Euler time-marching.
An open problem is that for some benchmark numerical examples, the iteration residue of
the fixed-point sweeping WENO scheme hangs at a truncation error level instead of settling down
to machine zero. This issue makes it difficult to determine the convergence criterion for the iteration
and challenging to apply the method to complex problems. To solve this issue, in this paper we apply
the multi-resolution WENO scheme developed in [J. Zhu and C.-W. Shu, J. Comput. Phys., 375 (2018), 659-683]
to the fifth order fixed-point sweeping WENO scheme and obtain an absolutely convergent fixed-point fast sweeping method for steady state of hyperbolic conservation laws, i.e., the residue of the fast sweeping iterations
converges to machine zero / round off errors for all benchmark problems tested.
Extensive numerical experiments, including solving difficult problems such as the shock reflection and supersonic flow past plates, are performed to show the accuracy, computational efficiency, and absolute convergence of the presented fifth order sweeping scheme.
\end{abstract}


\section{Introduction}
\label{sec1}

Steady state problems for hyperbolic conservation laws and Hamilton-Jacobi equations are
common mathematical models appearing in many applications, such as compressible fluid mechanics, optimal control,
geometric optics, image processing and computer vision, etc.
An essential property of these boundary value problems is that their
solution information propagates along characteristics starting
from the boundary. For spatial discretization of these hyperbolic type PDEs, weighted essentially non-oscillatory (WENO) schemes are a popular class of high order accuracy numerical methods.
They have the
advantage of attaining uniform high order accuracy in smooth regions of the solution while maintaining sharp and essentially
non-oscillatory transitions of discontinuities. The first WENO scheme was constructed in \cite{LOC} for a third-order accurate
finite volume version. In \cite{JiangShu}, third- and fifth-order accurate
finite difference WENO schemes in multi-space dimensions were
developed, with a general framework for the design of the
smoothness indicators and nonlinear weights. To deal with problems defined on complex domain geometries, WENO schemes on unstructured meshes were designed, e.g., see
 \cite{C.Hu,D.Levy,Y.Liu,Y.-T.Zhang,Y.-T.Zhangand C.-W.Shu}.
 Recently, WENO schemes on unequal-sized sub-stencils were designed \cite{JUN, JUNZ}, based on the idea in the central WENO schemes \cite{Levy1999, Capdeville2008}.  The WENO schemes on unequal-sized sub-stencils
 exhibit many interesting properties, especially advantages on unstructured meshes and solving steady state problems, see e.g. \cite{JUNZ2, JUNZ3}.

After spatial discretization of a steady state problem by a high order WENO scheme, a large nonlinear system is obtained.
An important factor which determines computational efficiency for solving steady state problems is the iterative scheme designed for the nonlinear system. For solving a highly coupled nonlinear system resulted from a high order WENO spatial discretization, one choice is to apply the Newton iterations or a more robust method such as the homotopy method \cite{HHSSXZ}. Another efficient approach is to use the fast sweeping iterative method to solve the large nonlinear system. Fast sweeping method is a class of efficient iterative methods which were original developed to solve static Hamilton-Jacobi equations (see e.g. \cite{Z,QZZ1,QZZ2,FLZ}). Fast sweeping methods utilize alternating sweeping
strategy to cover a family of characteristics in a certain direction simultaneously in each sweeping order.
Coupled with the Gauss-Seidel iterations, these methods can
achieve a fast convergence speed for computations of steady state solutions of hyperbolic PDEs.
In \cite{ZZQ}, fast sweeping methods were first time applied in solving the large nonlinear systems resulted from a WENO spatial discretization, where high order WENO fast sweeping schemes for solving static Hamilton-Jacobi equations were developed.
In the WENO fast sweeping methods \cite{ZZQ}, an explicit strategy in the iterative schemes was designed to avoid directly solving very complicated local nonlinear equations derived from high order WENO discretizations.
The methods were combined with accurate boundary treatment techniques in \cite{TMY}.
Fast sweeping techniques were also applied in discontinuous Galerkin (DG) methods to efficiently solve static Hamilton-Jacobi equations, see e.g. \cite{LSZZ,ZCLZS,WZ}.

In order to solve steady state problems of other types of hyperbolic PDEs such as the nonlinear hyperbolic conservation laws
by high order WENO fast sweeping schemes, the fixed-point sweeping methods were adopted \cite{S.Chen, WuLiang}. Fixed-point sweeping WENO methods were first designed in \cite{Y.-T.Zhang H.-K} for solving static Hamilton-Jacobi equations. This kind of
iterative methods have good properties which are different from other fast sweeping methods include that they have explicit forms and do not involve inverse operation of nonlinear local systems, and they can be applied in solving general hyperbolic equations using any monotone numerical fluxes and high order approximations (e.g. high order WENO approximations) easily. It is interesting
to find that although they were developed in different ways, the case of using the Lax-Friedrichs flux in the third order fixed-point sweeping schemes \cite{S.Chen} is equivalent to a third order Lax-Friedrichs fast sweeping method to solve steady state problems for hyperbolic conservation laws, which was designed in \cite{W.Chen}.
In the recent work \cite{WuLiang}, a fifth order fixed-point sweeping WENO method for
efficiently solving steady state problems of hyperbolic conservation laws was developed.
Numerical experiments showed that the method converges to steady state solutions much faster than regular time-marching approach,
and the acceleration of computation is essentially achieved via using fast sweeping techniques for
stability improvement of high order WENO schemes with a forward Euler time-marching to steady state solutions.
An open problem in the high order fixed-point sweeping WENO methods such as the fifth order scheme \cite{WuLiang} is that for some benchmark numerical examples, the iteration residue of
the fixed-point sweeping WENO scheme hangs at a truncation error level instead of settling down
to machine zero / round off errors. This issue makes it difficult to determine the convergence criterion for the iteration
and challenging to apply the method to complex problems in real applications.
Studies on high order WENO schemes on unequal-sized sub-stencils reveal
that they improve the convergence of classical high order WENO
schemes to steady state solutions so that the residue of TVD Runge-Kutta iterations settles down to machine zero / round off errors \cite{JUNZ3}.
In this paper, to resolve that open issue in high order fixed-point fast sweeping methods, we apply the fifth order multi-resolution WENO scheme in \cite{JUNZ} to the fifth order fixed-point sweeping scheme and obtain an absolutely convergent fixed-point fast sweeping method for steady state of hyperbolic conservation laws, namely, the residue of the fast sweeping iterations converges to machine zero / round off errors for all benchmark problems tested.
Extensive numerical experiments, including solving difficult problems such as the shock reflection and supersonic flow past plates, are performed to show the accuracy, computational efficiency, and absolute convergence of the presented fifth order sweeping scheme.

The rest of the paper is organized as follows. The detailed numerical algorithm is described in
Section 2. In Section 3 we provide extensive numerical experiments to test and study
the proposed method, and to perform comparisons of different methods. Concluding
remarks are given in Section 4.

\section {Description of the numerical methods}
\label{secfs}

We consider steady state problems of hyperbolic conservation laws with appropriate
boundary conditions. The two dimensional (2D) case has the following general form
\begin{equation}
f(u)_x+g(u)_y=R,
\label{(1.1)}
\end{equation}
where $u$ is the vector of the unknown conservative variables, $f(u)$ and $g(u)$ are the vectors of flux functions,
and $R(u,x,y)$ is the source term.
For example, the steady Euler system of equations in compressible fluid dynamics has
that $u = (\rho, \rho u, \rho v, E)^T$, $f(u) = (\rho u, \rho u^2 + p, \rho u v, u (E + p))^T$, and $g(u) = (\rho v, \rho u v, \rho v^2 +p, (E + p))^T$.
Here $\rho$ is the density of fluid, $(u,v)^T $ is the velocity vector, $p$ is the pressure, and $E=\frac{p}{\gamma'-1}+\frac{1}{2}\rho(u^2+v^2)$ is total energy where the constant $\gamma'=1.4$ for the case of air.
A spatial discretization of (\ref{(1.1)}) by a high order WENO scheme leads to a large nonlinear system of algebraic
equations with the size determined by the number of spatial grid points.

\subsection{The fifth order multi-resolution WENO discertization}
WENO schemes with unequal-sized sub-stencils \cite{JUN, JUNZ3} have shown nice property
in convergence to steady states. In this paper, we adopt the fifth order multi-resolution WENO scheme \cite{JUNZ}
for the spatial discretization of (\ref{(1.1)}) in order to achieve an absolutely convergent fixed-point fast sweeping method for steady state of hyperbolic conservation laws.

For the flux terms
${f}({u})_x+{g}({u})_y$, the conservative finite difference scheme is used to approximate them at a grid point $(x_i,y_j)$ of a uniform (or smoothly varying) mesh
\begin{equation}
\label{eq2.9}
\left(f(u)_x+g(u)_y\right)|_{x=x_i,y=y_j}\approx \frac{1}{\Delta x}(\hat f_{i+1/2,j}-\hat f_{i-1/2,j})
+\frac{1}{\Delta {y}}(\hat{g}_{i,j+1/2}-\hat{g}_{i,j-1/2}) ,
\end{equation}
where $\hat{f}_{i+1/2,j}$ and $\hat{g}_{i,j+1/2}$ are the numerical fluxes, and $\Delta x, \Delta {y}$ are grid sizes of $x, y$ directions respectively. In the following we describe the detailed procedure for reconstructing the numerical flux $\hat{f}_{i+1/2,j}$ in the $x$ direction. Similar
procedure is applied to the $y$ direction for $\hat{g}_{i,j+1/2}$.

For upwinding and linear stability of the scheme, we perform a splitting of the flux $f(u)$, i.e., $f(u)=f^{+}(u)+f^{-}(u)$, which satisfies the condition that $\frac{df^{+}(u)}{du}\geq{0}$ and $\frac{df^{-}(u)}{du}\leq{0}$ for the scalar case,
or the corresponding eigenvalues of that Jacobian matrix are positive / negative for the system case with a local characteristic
decomposition. A simple Lax-Friedrichs flux splitting $f^{\pm}(u)=\frac{1}{2}(f(u)\pm\alpha u)$ is used here, where $\alpha=max_{u}|f^{'}(u)|$.
Then each of them is approximated separately by the WENO scheme using different stencils to find numerical fluxes $\hat
f^+_{i+1/2,j}$ and $\hat f^-_{i+1/2,j}$ respectively.
The final numerical flux is $\hat{f}_{i+1/2,j}=\hat{f}_{i+1/2,j}^{+}+\hat{f}_{i+1/2,j}^{-}$.

To reconstruct $\hat{f}_{i+1/2,j}^{+}$ along the line $y=y_j$, we first identify the numerical values $f^{+}(u_{l,j})$ to be cell averages of a
function $h(x)$ on cell $I_l=[x_{l-1/2}, x_{l+1/2}]$ of the x-direction, for all $l$. Here
$x_{l+1/2}=(x_l+x_{l+1})/2$, and $u_{l,j}$ denotes the numerical solution of $u$ at the grid point $(x_l,y_j)$.
Then,
with the information of cell averages $$f^{+}(u_{l,j})=\overline{h}_{l}=\frac{1}{\triangle x}\int_{x_{l-1/2}}^{x_{l+1/2}}h(x)dx,$$
the following reconstruction algorithm is performed.

\bigskip
{\bf Reconstruction algorithm:}

{\bf Step 1.} We choose the central spatial stencils  $ T_{k}=\left\{ I_{i+1-k},\cdot\cdot\cdot ,I_{i-1+k}\right\} $, $k=1,2,3$, and reconstruct  $2k-2$ degree polynomials $q_{k}(x)$ which satisfy
$$\frac{1}{\triangle x}\int_{x_{l-1/2}}^{x_{l+1/2}}q_k(x)dx=\overline{h}_{l}, l=i-k+1, \cdot\cdot\cdot, i-1+k; k=1,2,3.$$

{\bf Step 2.} Obtain equivalent expressions for these reconstruction polynomials of different degrees. To keep consistent notation, we denote
$ p_{1}(x)=q_{1}(x)$, with similar ideas for the central WENO schemes \cite{Levy1999, Capdeville2008, D. L} as well, and compute
$$p_{2}(x)=\frac{1}{\gamma_{2,2}}q_{2}(x)-\frac{\gamma_{1,2}}{\gamma_{2,2}}p_{1}(x),$$
$$p_{3}(x)=\frac{1}{\gamma_{3,3}}q_{3}(x)-\frac{\gamma_{1,3}}{\gamma_{3,3}}p_{1}(x)-\frac{\gamma_{2,3}}{\gamma_{3,3}}p_{2}(x),$$
with $\gamma_{1,2}+\gamma_{2,2}=1$,  $\gamma_{1,3}+\gamma_{2,3}+\gamma_{3,3}=1$, and $\gamma_{2,2}\neq 0, \gamma_{3,3}\neq 0$.
In these expressions, $\gamma_{a,b}$ for $a=1,\cdot\cdot\cdot,b$ and $b=2,3$ are the linear weights. Based on a balance between the sharp and essentially non-oscillatory shock transitions in nonsmooth regions and accuracy in smooth regions, following the practice in \cite{M.Dumbser,X. Zhong,Y.Liu,JUN,JUNZ}, we take the linear weights as  $\gamma_{1,2}=\frac{1}{11}$, $\gamma_{2,2}=\frac{10}{11}$, $\gamma_{1,3}=\frac{1}{111}$, $\gamma_{2,3}=\frac{10}{111}$, $\gamma_{3,3}=\frac{100}{111}$.

{\bf Step 3.} Compute the smoothness indicators $\beta_{k}$, which measure how smooth the functions $p_{k}(x)$ for $k=2,3$ are in the interval $[x_{i-1/2},x_{i+1/2}]$. We use the same recipe for the smoothness indicators as that in \cite{JiangShu,CWS2}:
$$\beta_{k}= \sum_{\alpha=1}^{2k-2} \int_{x_{i-1/2}}^{x_{i+1/2}}{\triangle x}^{2\alpha-1}(\frac{d^{\alpha}p_{k}(x)}{dx^{\alpha}})^2dx,\qquad k=2,3. $$
The only exception is $\beta_{1}$, which is magnified from zero to a tiny value. See \cite{JUNZ} for details.

{\bf Step 4.} Compute the nonlinear weights based on the linear weights and the smoothness indicators. We adopt the WENO-Z type nonlinear weights as that in \cite{RB,MC}. First a quantity $\tau$ which depends on the absolute differences between the smoothness indicators is calculated:
$\tau=(\frac{\sum_{l_1=1}^{2}|\beta_{3}-\beta_{l_1}|}{2})^2$.
The nonlinear weights are then computed as
$$\omega_{l_1}=\frac{\bar{\omega}_{l_1}}{\sum_{l_2=1}^{3}\bar{\omega}_{l_2}}, \quad \bar{\omega}_{l_1}=\gamma_{l_1,3}(1+\frac{\tau}{\varepsilon+\beta_{l_1}}), \qquad l_1=1,2,3.$$
$\varepsilon$ is a small value to avoid that the denominator becomes zero. In this paper, $\varepsilon$ is taken to be $10^{-6}$ for all numerical examples.

{\bf Step 5.}
The final reconstructed numerical flux $\hat{f}_{i+1/2,j}^{+}$ is given by
$$\hat{f}_{i+1/2,j}^{+}=\sum_{l_1=1}^{3}\omega_{l_1}p_{l_1}(x_{i+1/2}).$$
The reconstruction of the numerical flux $\hat{f}_{i+1/2,j}^{-}$ is mirror-symmetric with respect to $x_{i+1/2}$.

\subsection{Absolutely convergent fixed-point sweeping WENO scheme}

After the spatial discretization of the PDE (\ref{(1.1)}) by a WENO scheme with unequal-sized sub-stencils such as
the fifth order multi-resolution WENO scheme described in the last section, we obtain a nonlinear algebraic system
\begin{equation}
0=-(\hat{f}_{i+1/2,j}-\hat{f}_{i-1/2,j})/  \Delta
x-(\hat{g}_{i,j+1/2}-\hat{g}_{i,j-1/2})/\Delta
y+R(u_{ij},x_i,y_j),\hspace{.1in}i=1,\cdots,N; j=1,\cdots, M,
\label{eq2}
\end{equation}
where $N, M$ are the number of grid points in the $x$ and $y$ directions respectively.

The fixed-point fast sweeping schemes for solving the nonlinear system (\ref{eq2}) are based on iterative schemes
of time marching type for solving steady state problems.
Time marching approach for solving steady state problems is essentially a Jacobi type
fixed-point iterative scheme. The right-hand-side (RHS) of (\ref{eq2}) is a nonlinear function of numerical values
at the grid points of computational stencils. We denote it by $L$. The forward Euler (FE) time marching method
with time step size $\Delta t_n=\f{\gamma}{\alpha_x/\Delta x+\alpha_y/\Delta
y}$ is actually the following Jacobi type fixed-point iterative scheme
$$u_{ij}^{n+1}=u_{ij}^{n}+\frac{\gamma}{\alpha_{x}/\Delta x +\alpha_{y}/\Delta y}L(u_{i-r,j}^{n},\cdot\cdot\cdot,u_{i+s,j}^{n};u_{ij}^{n};u_{i,j-r}^{n},\cdot\cdot\cdot u_{i,j+s}^{n}),$$
\begin{equation}
\qquad\qquad\qquad\qquad\qquad\qquad\qquad\qquad\qquad\qquad i=1,\cdot\cdot\cdot,N;j=1,\cdot\cdot\cdot,M;      \qquad\qquad\qquad \qquad
\label{eq10}
\end{equation}
where $r,s$ are values which depend on the order of the WENO reconstruction. For the fifth order multi-resolution WENO scheme used here, we have $r=s=3$. $n$ denotes the iteration step. The values $\alpha_{x}=\max_{u}{|f'(u)|}$  and $\alpha_{y}=\max_{u}{|g'(u)|}$ for the scalar equations, or they are the maximum absolute values of eigenvalues of the Jacobian matrices $f'(u)$ and $g'(u)$ for the system cases. Actually $\alpha_{x}$ and $\alpha_{y}$ represent the maximum characteristic speeds in each spatial direction, and they are updated in every iteration. $\gamma$ is a parameter. To guarantee that the fixed-point iteration is a contractive mapping and converges, suitable values of $\gamma$ need to be chosen. In the context of time marching methods, $\gamma$ is actually the Courant-Friedrichs-Lewy (CFL) number. Because accuracy of the scheme in the time discretization does not affect the final results of solutions of steady state problems, it sounds good that the FE time marching method is a better choice than high order time schemes such as Runge-Kutta (RK) schemes
since its simplicity with only single step and single stage involved.  However, high order linear schemes which are the building blocks for high order spatial discretization schemes such as WENO schemes, has linear stability issues
when they are marched by the forward Euler time discretization to solve hyperbolic PDEs. Even with the
help of a nonlinearly stable discretization such as a WENO scheme, the Jacobi fixed-point iterative scheme (\ref{eq10}) still needs many iteration steps to converge \cite{WuLiang}.

One way to resolve this issue is to use high order TVD RK schemes \cite{C.W.S, GST}, which can maintain both linear and nonlinear stability when they are coupled with high order nonlinearly stable spatial discretizations. Similar as the FE time marching, RK time marching methods are also Jacobi type fixed-point iterations. For example, the popular third order TVD RK (TVD RK3) scheme \cite{C.W.S} can be written as the following Jacobi type fixed-point iteration scheme for solving the nonlinear system (\ref{eq2}):
$$u_{ij}^{(1)}=u_{ij}^{n}+\frac{\gamma}{\alpha_x/\Delta x+\alpha_y/\Delta
y}L(u_{i-r,j}^{n},\cdots,u_{i+s,j}^{n};u_{i,j}^{n};u_{i,j-r}^{n},\cdots, u_{i,j+s}^{n}),$$
\begin{equation}
\hspace{2in} i=1,\cdots,N; j=1,\cdots,M,
\label{eq7}
\end{equation}
$$u_{ij}^{(2)}=\frac{3}{4}u_{ij}^{n}+\frac{1}{4}u_{ij}^{(1)}+
\frac{\gamma}{4(\alpha_x/\Delta
x+\alpha_y/\Delta y)}L(u_{i-r,j}^{(1)},\cdots,u_{i+s,j}^{(1)};u_{i,j}^{(1)};u_{i,j-r}^{(1)},\cdots u_{i,j+s}^{(1)}),$$
\begin{equation}
\hspace{2in} i=1,\cdots,N; j=1,\cdots,M,
\label{eq8}
\end{equation}
$$u_{ij}^{n+1}=\frac{1}{3}u_{ij}^{n}+\frac{2}{3}u_{ij}^{(2)}+
\frac{2\gamma}{3(\alpha_x/\Delta x+\alpha_y/\Delta y)}
L(u_{i-r,j}^{(2)},\cdot\cdot\cdot,u_{i+s,j}^{(2)};u_{i,j}^{(2)};u_{i,j-r}^{(2)},\cdot\cdot\cdot u_{i,j+s}^{(2)}),$$
\begin{equation}
\hspace{2in}  i=1,\cdots,N; j=1,\cdots,M.
\label{eq9}
\end{equation}
Although the TVD RK3 scheme (\ref{eq7})-(\ref{eq9}) involves three stages in one time step, it can permit a much larger $\gamma$ than the FE scheme (\ref{eq10}) in order to ensure that the fixed-point iterations converge, hence
much less iterations are needed as that shown in \cite{WuLiang} and also the next section of this paper. It is a popular method to march numerical solutions of high order WENO schemes to steady states, see e.g. \cite{JiangShu, JUNZ3}.

In order to achieve faster convergence to steady state solutions of high order WENO schemes for solving hyperbolic PDEs than the Jacobi type fixed-point iterations, one approach is to apply the fast sweeping techniques so that
the important characteristics property of hyperbolic PDEs can be utilized in the iterations.
As shown in the numerical experiments of \cite{WuLiang} and this paper, by applying the fast sweeping techniques
to the fixed-point schemes, we can obtain much more efficient iterative schemes. Especially for the FE scheme, the
number of iteration steps to the steady state is reduced significantly and the CFL number
$\gamma$ is much larger than that in the Jacobi iterations. In fact, after applying the fast sweeping techniques, the FE scheme can have a comparable CFL number $\gamma$ as that of the TVD RK3 scheme, and it is more efficient than the TVD RK3 scheme to reach steady state of the solutions.

The fast sweeping methods have two essential components, i.e., the Gauss-Seidel philosophy and alternating direction iterations. The Gauss-Seidel philosophy requires that the newest
numerical values of $u$ are used in the fifth order multi-resolution WENO reconstruction stencils as long as they are available. With the fifth order multi-resolution WENO reconstruction described in the last section, we form an absolutely convergent fixed-point fast sweeping method. The residue of the new fixed-point fast sweeping iterations converges to machine zero / round off errors for all benchmark problems tested in the next section. The form of this FE type fixed-point fast sweeping scheme is:
$$u_{ij}^{n+1}=u_{ij}^{n}+\frac{\gamma}{\alpha_{x}/\Delta x +\alpha_{y}/\Delta  y}L(u_{i-r,j}^{\ast},\cdots,u_{i+s,j}^{\ast};u_{i,j}^{n};u_{i,j-r}^{\ast},\cdots, u_{i,j+s}^{\ast}),$$
\begin{equation}
\hspace{2in} i=i_{1},\cdots,i_{N};j=j_{1},\cdots,j_{M}.
\label{eq11}
\end{equation}
Here the operator $L$ denotes the fifth order multi-resolution WENO spatial discretization. The iteration direction
that is marked as ``$i=i_{1},\cdots,i_{N};j=j_{1},\cdots,j_{M}$'' means that the iterations in the scheme (\ref{eq11}) do not just proceed in only one direction $i = 1:N, j = 1:M$ as that in the Jacobi type
schemes, but in the following four alternating directions repeatedly,
$$\mbox{(1) } i=1:N, j=1:M;$$
$$\mbox{(2) } i=N:1, j=1:M;$$
$$\mbox{(3) } i=N:1, j=M:1;$$
$$\mbox{(4) } i=1:N, j=M:1.$$
Since the strategy of alternating direction sweepings utilizes the characteristics property of hyperbolic PDEs, combining with the Gauss-Seidel philosophy,
we are able to observe the acceleration of convergence speed, which will be shown in the following numerical experiments. By the  Gauss-Seidel philosophy, we use the newest numerical values on
the computational stencil of the WENO scheme whenever they are available. That is the
reason we use the notation $u^{\ast}$ to represent the values in the scheme (\ref{eq11}), and $u_{k,l}^{\ast}$ could be $u_{k,l}^{n}$ or $u_{k,l}^{n+1}$, depending on the current sweeping direction.

\bigskip
\noindent {\bf Remark 1.} The fast sweeping techniques can also be applied to the TVD RK3 scheme (\ref{eq7})-(\ref{eq9}) as that in \cite{WuLiang}. However, it was found in \cite{WuLiang} that the scheme resulted by applying fast sweeping method to
the TVD RK3 scheme (\ref{eq7})-(\ref{eq9}) is less efficient than the FE fast sweeping method. Hence in this paper we only focus on the FE fast sweeping method coupled with the fifth order multi-resolution WENO reconstruction, i.e., the scheme (\ref{eq11}). Also we emphasize that the scheme (\ref{eq11}), the FE fast sweeping method coupled with the fifth order multi-resolution WENO reconstruction, is a novel method which is different from the method in our previous work \cite{WuLiang}. It improves the previous scheme in \cite{WuLiang} by achieving the absolute convergence in the fixed-point fast sweeping iterations to reach steady state of hyperbolic conservation laws, i.e., the residue of the fast sweeping iterations converges to machine zero / round off errors.

\bigskip
\noindent {\bf Remark 2.} It is interesting to see that in the fast sweeping scheme (\ref{eq11}), during the iteration process the fifth order multi-resolution WENO reconstructions for numerical fluxes at the cell interfaces need to be performed individually for grid points sharing the same cell interface, due to the Gauss-Seidel philosophy. For example, the newest numerical values at the WENO reconstruction stencil for computing the numerical flux $\hat{f}_{i+1/2,j}$ are different when we update the numerical values $u_{i,j}$ and $u_{i+1,j}$. Since the Gauss-Seidel philosophy in the fast sweeping scheme requires that the newest numerical values are always used whenever they are available, the numerical flux $\hat{f}_{i+1/2,j}$ needs to be computed separately for updating
$u_{i,j}$ and $u_{i+1,j}$ during the fast sweeping iterations. This is different from the Jacobi iterations. So during the iteration process, the fast sweeping scheme does {\it not} preserve the conservation. However, this does
{\it not} affect the conservation of the final results at all, because once the fast sweeping iterations converge,
the new and old values converge to the same value, hence numerical fluxes $\hat{f}_{i+1/2,j}$ which are computed separately for $u_{i,j}$ and $u_{i+1,j}$ also converge to the same value. Then the final conservative system (\ref{eq2}) is resolved and the conservation of the final results is preserved. Also from the numerical experiments shown in the next section, the fast sweeping iterations are much more efficient than the Jacobi type time marching iterations.

\section{Numerical experiments}

In this section we perform numerical experiments to test the performance of the proposed fifth order
absolutely convergent fixed-point fast sweeping WENO scheme
for solving steady-state problems.
These numerical examples come from the literature about numerical studies on steady state solutions of
hyperbolic conservation laws, e.g., \cite{WuLiang,JUNZ4,JUNZ3}.
Computational efficiency of the fast sweeping scheme and the other two iterative schemes discussed in the last section
is compared. For the convenience of presentation, we name
the absolutely convergent fixed-point fast sweeping WENO scheme (\ref{eq11}) ``FE fast sweeping scheme'',
the forward Euler time marching scheme (\ref{eq10}) ``FE Jacobi scheme'', the popular TVD RK3 time marching scheme (\ref{eq7})-(\ref{eq9}) ``RK Jacobi
scheme''.  We perform mesh refinement study and
 calculate $L_{1}$ and $L_{\infty}$ numerical errors and accuracy orders of these iterative schemes.  Iteration numbers and CPU times of each iterative scheme to converge are recorded
for comparison of their efficiency. The convergence of the iterations is measured by the average residue defined as
\begin{equation}
ResA=\sum_{i=1}^{\cal{N}}\frac{|R_{i}|}{\cal{N}}.
\label{averres}
\end{equation}
Here $R_i$ is the local residue at the grid point $i$ and it is defined as
$$R_{i}=\frac{u_{i}^{n+1}-u_{i}^{n}}{\Delta t_{n}}.$$
$\cal{N}$ is the total number of grid points and $n$ is the iteration step.
$\Delta t_n=\f{\gamma}{\alpha_x/\Delta x+\alpha_y/\Delta
y}$.
For a system of equations, the average residue $ResA$ is calculated based on the mean of all components of the system.
For example, a 2D Euler system has four equations with four conservative variables, so the average residue is defined as $ResA=\sum_{i=1}^{\cal{N}}\frac{|R1_{i}|+|R2_{i}|+|R3_{i}|+|R4_{i}|}{4\cal{N}}$ where the $R\ast_{i}$'s are the local residuals of different conservative variables, i.e., $R1_{i}=\frac{\rho_{i}^{n+1}-\rho_{i}^{n}}{\Delta t_n}, R2_{i}=\frac{(\rho u)_{i}^{n+1}-(\rho u)_{i}^{n}}{\Delta t_n}, R3_{i}=\frac{(\rho v)_{i}^{n+1}-(\rho v)_{i}^{n}}{\Delta t_n}, R4_{i}=\frac{E_{i}^{n+1}-E_{i}^{n}}{\Delta t_n}$.  The convergence criterion
for most examples
is set to be $ResA<10^{-12}$ except that in some examples we use an even smaller threshold value to study the
levels that the residues can reach. In this paper, the number of iterations reported in the numerical simulations counts a complete update of
numerical values in all grid points once as one iteration.

\subsection{Examples with smooth solutions}

We first test these schemes using problems with smooth solutions.

\bigskip
\noindent{\bf Example 1. A 1D Burgers' equation}

\noindent We compute the steady state solution of the following 1D Burgers' equation with a source term
$$u_{t}+\left(\frac{u^{2}}{2}\right)_{x}=\sin(x)\cos(x) ,\qquad x\in\left[\frac{1}{4}\pi, \frac{3}{4}\pi\right].$$
The initial condition $u(x,0)=\beta\sin(x)$ is used as the initial guess in the iterations. At the left boundary $x=(1/4)\pi$, an inflow boundary condition is imposed with $u((1/4)\pi,t)=\sqrt{2}/2$.  The outflow boundary condition is applied at the right boundary $x=(3/4)\pi$.  If $\beta>1$, the problem has the unique steady state solution $u(x,\infty)=\sin(x)$ . We take $\beta=2$ and use these three different iterative
schemes based on the fifth order multi-resolution WENO spatial discretization to compute the steady state solution.  For the outflow boundary point $x =(3/4)\pi$ itself and the associated ghost points to the right of x =(3/4)$\pi$
in the stencil of the WENO scheme, extrapolation by a degree $4$ polynomial is used to compute numerical values at them. In this example we use $ResA<10^{-13}$ as the iteration convergence criterion.

The results of numerical accuracy, iteration numbers and CPU costs for three different iterative schemes are presented in Table \ref{tab:Margin_settings1}.
It is observed that all schemes achieve comparable numerical errors and the fifth order
accuracy when the iterations converge. Comparing the algorithm efficiency to reach the steady state, the FE Jacobi scheme (\ref{eq10}) requires a very small CFL number $\gamma=0.1$ to have the iterations converge. The reason is that a forward Euler time discretization with a
high order linear upwind spatial discretization suffers from
linear stability issue. When the nonlinearly stable WENO discretization is applied, it alleviates the linear instability problem. As a result of balance, the FE Jacobi scheme converges under a tiny CFL number. That leads to the largest iteration numbers and the most CPU time among these three iterative
schemes. By using a high order TVD RK scheme, for example the TVD RK3, the RK Jacobi scheme
(\ref{eq7})-(\ref{eq9}) with a high order WENO discretization is both linearly and nonlinearly stable. Hence here
a much larger CFL number $\gamma = 1.0$ can be used to make the iterations converge. We can see that both the iteration numbers and CPU costs
are reduced a lot by using the RK Jacobi scheme rather than the FE Jacobi scheme. From Table \ref{tab:Margin_settings1}, it can be seen that fast
sweeping techniques accelerate the convergence to the steady state significantly.
Furthermore, it is important to notice that
the FE fast sweeping scheme can have a CFL number $\gamma = 1.0$ which is the same as the TVD RK3 scheme. This shows that the fast sweeping technique improves the linear stability of
the forward Euler scheme when it is coupled with a high order spatial discretization.
Since the time accuracy does not contribute to the accuracy of a steady state solution, the forward
Euler time marching should be preferred because of its simplicity
comparing to multi-stage Runge-Kutta schemes. However, due to the linear stability issue when it is coupled with a high order spatial discretization, the forward Euler scheme is not practically useful. Here we can see that this issue is
resolved by applying the fast sweeping technique, namely, using the FE fast sweeping scheme rather
than the FE Jacobi scheme. Actually, as shown in Tables \ref{tab:Margin_settings1} about the iteration numbers and CPU times, the FE fast sweeping scheme
is the most efficient one among all three schemes.


\begin{table}
		\centering
		\begin{tabular}{|c|c|c|c|c|c|c|}\hline
			\multicolumn{7}{|c|}{FE Jacobi, $\gamma$=0.1 }\\\hline
			N&L$_{1}$ error&L$_{1}$ order&L$_{\infty} $ error&L$_{\infty} $ order&iter$\sharp$&CPU time\\\hline

    10    & 6.27E-07 &   -    & 1.54E-06     &    -   & 1153  & 2.0E-02 \\
    \hline
    20    & 1.93E-08 & 5.02  & 8.07E-08   & 4.25  & 1458  & 5.5E-02 \\
    \hline
    40    & 8.91E-10 & 4.44  & 3.21E-09     & 4.65  & 1749  & 0.14 \\
    \hline
    80    & 3.32E-11 & 4.75  & 1.11E-10    & 4.85  & 2310  & 0.30\\
    \hline
    160   & 1.13E-12 & 4.87  & 3.66E-12   & 4.92  & 3875  & 0.98\\
    \hline
    320   & 4.98E-14 & 4.51  & 1.70E-13    & 4.43  & 7196  & 3.63 \\
    \hline

		\end{tabular}
\begin{tabular}{|c|c|c|c|c|c|c|}\hline
			\multicolumn{7}{|c|}{RK Jacobi, $\gamma$=1.0 }\\\hline
			N&L$_{1}$ error&L$_{1}$ order&L$_{\infty} $ error&L$_{\infty} $ order&iter$\sharp$&CPU time\\\hline
            10    & 8.11E-07 &   -    & 3.15E-06      &    -   & 285   & 2.0E-03 \\
    \hline
    20    & 2.29E-08 & 5.15  & 1.19E-07     & 4.73  & 330   & 8.0E-03 \\
    \hline
    40    & 9.49E-10 & 4.59  & 4.00E-09     & 4.89  & 429   & 1.6E-02 \\
    \hline
    80    & 3.41E-11 & 4.80   & 1.29E-10    & 4.95  & 630   & 6.6E-02 \\
    \hline
    160   & 1.15E-12 & 4.89   & 4.35E-12  & 4.90   & 1137  & 0.19 \\
    \hline
    320   & 4.10E-14 & 4.81  & 1.55E-13   & 4.81  & 1953  & 0.58 \\
    \hline
		\end{tabular}
\begin{tabular}{|c|c|c|c|c|c|c|}\hline
			\multicolumn{7}{|c|}{FE fast sweeping, $\gamma$=1.0 }\\\hline
			N&L$_{1}$ error&L$_{1}$ order&L$_{\infty} $ error&L$_{\infty} $ order&iter$\sharp$&CPU time\\\hline

    10    & 6.27E-07 &   -    & 1.54E-06     &   -    & 130   & 2.0E-03 \\
    \hline
    20    & 1.93E-08 & 5.02  & 8.07E-08    & 4.25  & 142   & 4.0E-03 \\
    \hline
    40    & 8.91E-10 & 4.44  & 3.21E-09   & 4.65  & 155   & 1.1E-02 \\
    \hline
    80    & 3.32E-11 & 4.75  & 1.11E-10     & 4.85  & 210   & 2.9E-02 \\
    \hline
    160   & 1.12E-12 & 4.89  & 3.64E-12    & 4.93  & 328   & 9.7E-02 \\
    \hline
    320   & 3.64E-14 & 4.95  & 1.22E-13    & 4.90   & 550   & 0.28 \\
    \hline
		\end{tabular}

		\caption{Example 1, A 1D Burgers' equation. Accuracy, iteration numbers and CPU times of three different iterative schemes. CPU time unit: second}
		\label{tab:Margin_settings1}
	\end{table}

\bigskip
\noindent {\bf Example 2. A 1D shallow water equation}

\noindent In this example, we use these iterative schemes to solve a one-dimensional
shallow water equation
\[
\left(\begin{array}{cc}
 h  \\hu
\end{array}\right)_{t}+
\left(\begin{array}{cc}
 hu  \\hu^{2}+\frac{1}{2}gh^{2}
\end{array}\right)_{x}=
\left(\begin{array}{cc}
 0 \\ -ghb_{x}
\end{array}\right),
\]
where $h$ denotes the water height, $u$ is the velocity of the water, $b(x)$ represents the bottom topography, and $g$ is the gravitational constant. The bottom
topography is smooth and given by the function
$$b(x)=5e^{-\frac{2}{5}(x-5)^2},    x\in [0,10].$$
The steady state solution of the problem is $$h+b=10, \qquad hu=0.$$

\begin{table}
		\centering
		\begin{tabular}{|c|c|c|c|c|c|c|}\hline
			\multicolumn{7}{|c|}{FE Jacobi, $\gamma$=0.1 }\\\hline
			N&L$_{1}$ error&L$_{1}$ order&L$_{\infty} $ error&L$_{\infty} $ order&iter$\sharp$&CPU time\\\hline
    20    & 3.53E-03 & -  & 2.12E-02    & -  & 5676  & 0.14 \\
    \hline
    40    & 9.31E-05 & 5.25  & 1.37E-03     & 3.95  & 4512  & 0.20 \\
    \hline
    80    & 1.58E-06 & 5.88  & 3.45E-05     & 5.31  & 7314  & 0.66 \\
    \hline
    160   & 1.59E-08 & 6.63  & 4.54E-07    & 6.25  & 13023 & 2.32 \\
    \hline
    320   & - &  -  &  -  & -  & not conv. & - \\
    \hline
\end{tabular}
\begin{tabular}{|c|c|c|c|c|c|c|}\hline
			\multicolumn{7}{|c|}{RK Jacobi, $\gamma$=1.0 }\\\hline
			N&L$_{1}$ error&L$_{1}$ order&L$_{\infty} $ error&L$_{\infty} $ order&iter$\sharp$&CPU time\\\hline           20    & 3.53E-03 & -  & 2.12E-02     & -  & 321   & 2.4E-02 \\
\hline
40    & 9.31E-05 & 5.25  & 1.37E-03     & 3.95  & 459   & 6.4E-02 \\
\hline
80    & 1.58E-06 & 5.88  & 3.45E-05     & 5.31  & 741   & 0.22 \\
\hline
160   & 1.59E-08 & 6.63  & 4.54E-07    & 6.25  & 1161  & 0.62\\
\hline
320   & 2.02E-10 & 6.30   & 6.83E-09    & 6.05  & 1734  & 1.86 \\
\hline
\end{tabular}
\begin{tabular}{|c|c|c|c|c|c|c|}\hline
			\multicolumn{7}{|c|}{FE fast sweeping, $\gamma$=1.0 }\\\hline
			N&L$_{1}$ error&L$_{1}$ order&L$_{\infty} $ error&L$_{\infty} $ order&iter$\sharp$&CPU time\\\hline
20    & 3.53E-03 & -  & 2.12E-02     & -  & 221   & 1.3E-02 \\
\hline
40    & 9.31E-05 & 5.25  & 1.37E-03    & 3.95  & 121   & 1.7E-02 \\
\hline
80    & 1.58E-06 & 5.88  & 3.45E-05     & 5.31  & 144   & 4.2E-02 \\
\hline
160   & 1.59E-08 & 6.63  & 4.54E-07    & 6.25  & 228   & 0.13 \\
\hline
320   & 2.03E-10 & 6.30   & 6.83E-09    & 6.05  & 379   & 0.42 \\
\hline
\end{tabular}
\caption{Example 2, A 1D shallow water equation. Accuracy, iteration numbers and CPU times of three different iterative schemes. CPU time unit: second}
\label{tab:Margin_settings2}
\end{table}

\noindent To test these iterative schemes, we use the exact solution as the initial guess in the iterations. Because the exact steady
state solution of the PDE itself is not the numerical steady states of the schemes, the convergence behavior of iterative schemes starting from it can be observed. We use the exact solution to impose the
numerical values at boundary points.
The simulation results of these three iterative schemes are reported in Table \ref{tab:Margin_settings2}.
 It can be seen that that these schemes give similar numerical errors and higher than
fifth order accuracy when the mesh is refined.
We observe similar convergence behavior of these three different iterative schemes as that in Example 1, except that in this example,
residues of the FE Jacobi
scheme with $N=320$ stop at the level of $10^{-9}$ and fail to reach the
convergence criterion $ResA<10^{-12}$, even using a very small CFL number $\gamma=0.1$. Large number of iterations and the most CPU time
among three iterative schemes are needed for the FE Jacobi scheme to converge. On the other hand,
the FE fast sweeping scheme can use a much larger CFL number $\gamma=1.0$ and converges with much smaller iteration numbers and CPU times, including the case with the
most refined mesh $N=320$ without any difficulty. Consistent with Example 1, the FE fast sweeping scheme is more efficient than the TVD RK3 (RK Jacobi) scheme and it turns out to be the most efficient scheme among these three iterative schemes. Specifically, CPU times of the FE fast sweeping
scheme is about $20\%$ of that of the RK Jacobi scheme on refined meshes.

\bigskip
\noindent{\bf Example 3. A 2D Burgers' equation}

\noindent We use the iterative schemes to solve the steady state problem of a two-dimensional Burgers' equation with a source term
\begin{equation*}
\begin{aligned}
u_{t}+(\frac{1}{\sqrt{2}}\frac{u^{2}}{2})_{x}+(\frac{1}{\sqrt{2}}\frac{u^{2}}{2})_{y}=\sin(\frac{x+y}{\sqrt{2}})&\cos(\frac{x+y}{\sqrt{2}}),
\\&(x,y)\in [\frac{\pi}{4\sqrt{2}},\frac{3 \pi}{4\sqrt{2}}]\times[\frac{\pi}{4\sqrt{2}},\frac{3 \pi}{4\sqrt{2}}].
\end{aligned}
\end{equation*}
The initial condition
$$u(x,y,0)=\beta \sin(\frac{x+y}{\sqrt{2}})$$
is used as the initial guess in the iterations, and the value of $\beta$ is taken as $1.5$. The exact steady state solution of the problem is
$$u(x,y,\infty)=\sin(\frac{x+y}{\sqrt{2}})$$
and it is imposed on boundary points
for the purpose of testing these three iterative schemes. $ResA<10^{-13}$ is used as the iteration convergence criterion. The results for three different iterative schemes are reported in Table \ref{tab:Margin_settings3}.
Similar as the 1D examples, all three schemes achieve comparable numerical errors and the fifth order accuracy when they converge. However the FE Jacobi scheme requires a very small CFL number to achieve the convergence, hence much more   iterations and CPU times than the other two schemes are needed to reach the steady state.
The FE fast sweeping scheme is the most efficient one among all three iterative
methods as shown in Table \ref{tab:Margin_settings3}. We can see that CPU times which are needed by the FE fast sweeping
scheme is only about $35\%$ of that of the RK Jacobi (TVD RK3) scheme on the most refined mesh.

\begin{table}
		\centering
		\begin{tabular}{|c|c|c|c|c|c|c|}\hline
			\multicolumn{7}{|c|}{FE Jacobi, $\gamma$=0.1 }\\\hline
			$N \times N$ &L$_{1}$ error&L$_{1}$ order&L$_{\infty} $ error&L$_{\infty} $ order&iter$\sharp$&CPU time\\\hline

    $10\times 10$ & 1.47E-08 &   -    & 8.60E-08  &   -    & 1054  & 0.31 \\
    \hline
    $20\times 20$ & 6.14E-10 & 4.58  & 3.28E-09 & 4.71  & 1317  & 1.53 \\
    \hline
    $40\times 40$ & 2.22E-11 & 4.79  & 1.24E-10  & 4.73  & 1850  & 8.88 \\
    \hline
    \end{tabular}
\begin{tabular}{|c|c|c|c|c|c|c|}\hline
			\multicolumn{7}{|c|}{RK Jacobi, $\gamma$=1.0 }\\\hline
			$N\times N$ &L$_{1}$ error&L$_{1}$ order&L$_{\infty} $ error&L$_{\infty} $ order&iter$\sharp$&CPU time\\\hline

    $10\times 10$ & 1.81E-08 &   -    & 1.43E-07  &   -    & 279   & 9.4E-02\\
    \hline
    $20\times 20$ & 6.87E-10 & 4.72  & 5.12E-09  & 4.80   & 348   & 0.42 \\
    \hline
    $40\times 40$ & 2.35E-11 & 4.87  & 1.64E-10 & 4.96  & 519   & 2.48 \\
    \hline
    \end{tabular}
\begin{tabular}{|c|c|c|c|c|c|c|}\hline
			\multicolumn{7}{|c|}{FE fast sweeping, $\gamma$=1.0 }\\\hline
			$N\times N$ &L$_{1}$ error&L$_{1}$ order&L$_{\infty} $ error&L$_{\infty} $ order&iter$\sharp$&CPU time\\\hline
            $10\times 10$ & 1.81E-08 &   -    & 1.43E-07  &   -    & 120   & 3.3E-02\\
    \hline
    $20 \times 20$ & 6.87E-10 & 4.72  & 5.12E-09  & 4.80   & 137   & 0.18 \\
    \hline
    $40\times 40$ & 2.35E-11 & 4.87  & 1.71E-10  & 4.91  & 182   & 0.88 \\
    \hline
    \end{tabular}
	\caption{Example 3, A 2D Burgers' equation. Accuracy, iteration numbers and CPU times of three different iterative schemes. CPU time unit: second}
		\label{tab:Margin_settings3}
	\end{table}

\bigskip
\bigskip
\noindent{\bf Example 4. A 2D Euler system of equations with source terms}

\begin{figure}
\centering
\subfigure[FE Jacobi scheme]{
\begin{minipage}[t]{0.25\linewidth}
\centering
\includegraphics[width=1.5in]{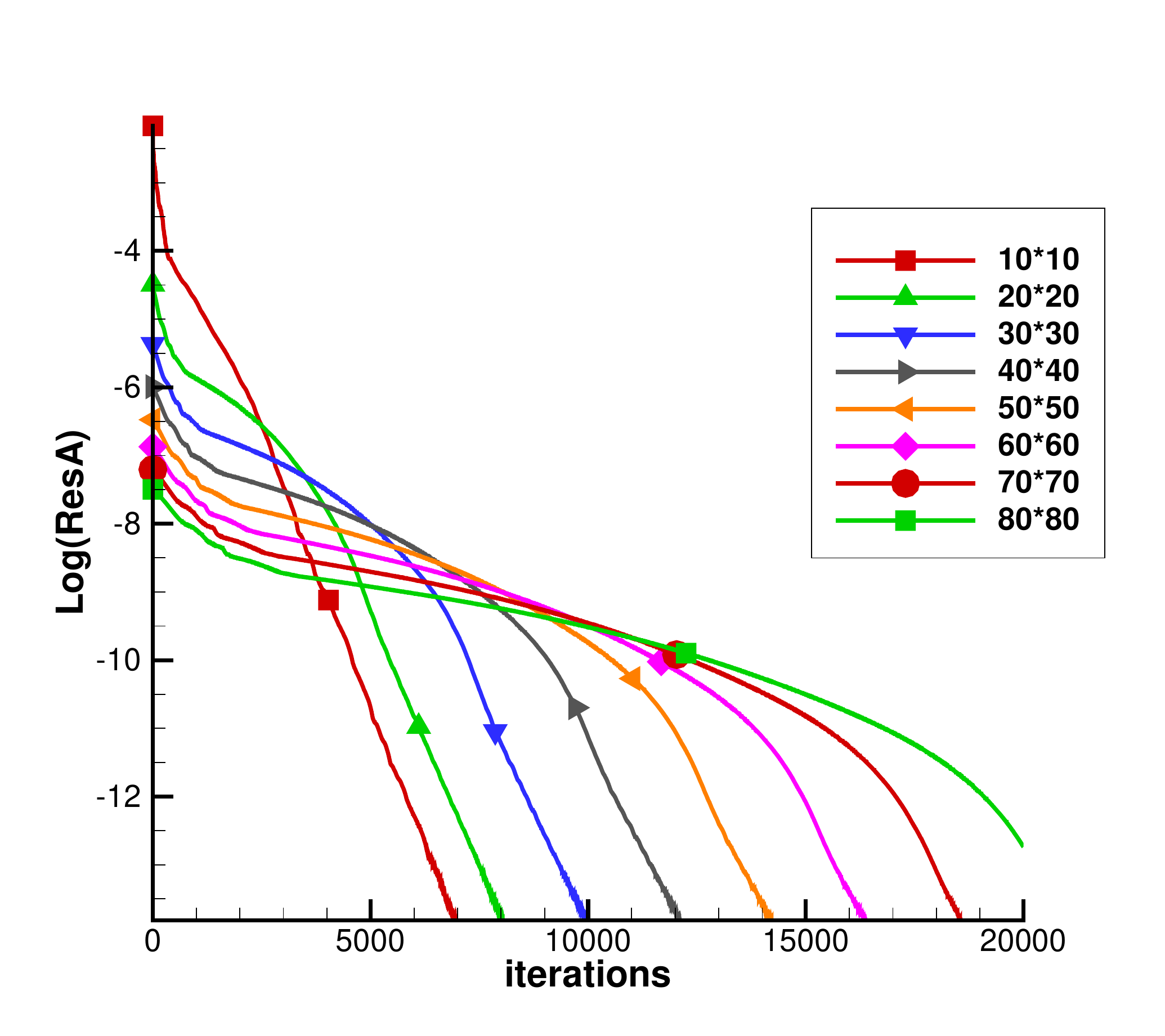}
\end{minipage}%
}%
\subfigure[RK Jacobi scheme]{
\begin{minipage}[t]{0.25\linewidth}
\centering
\includegraphics[width=1.5in]{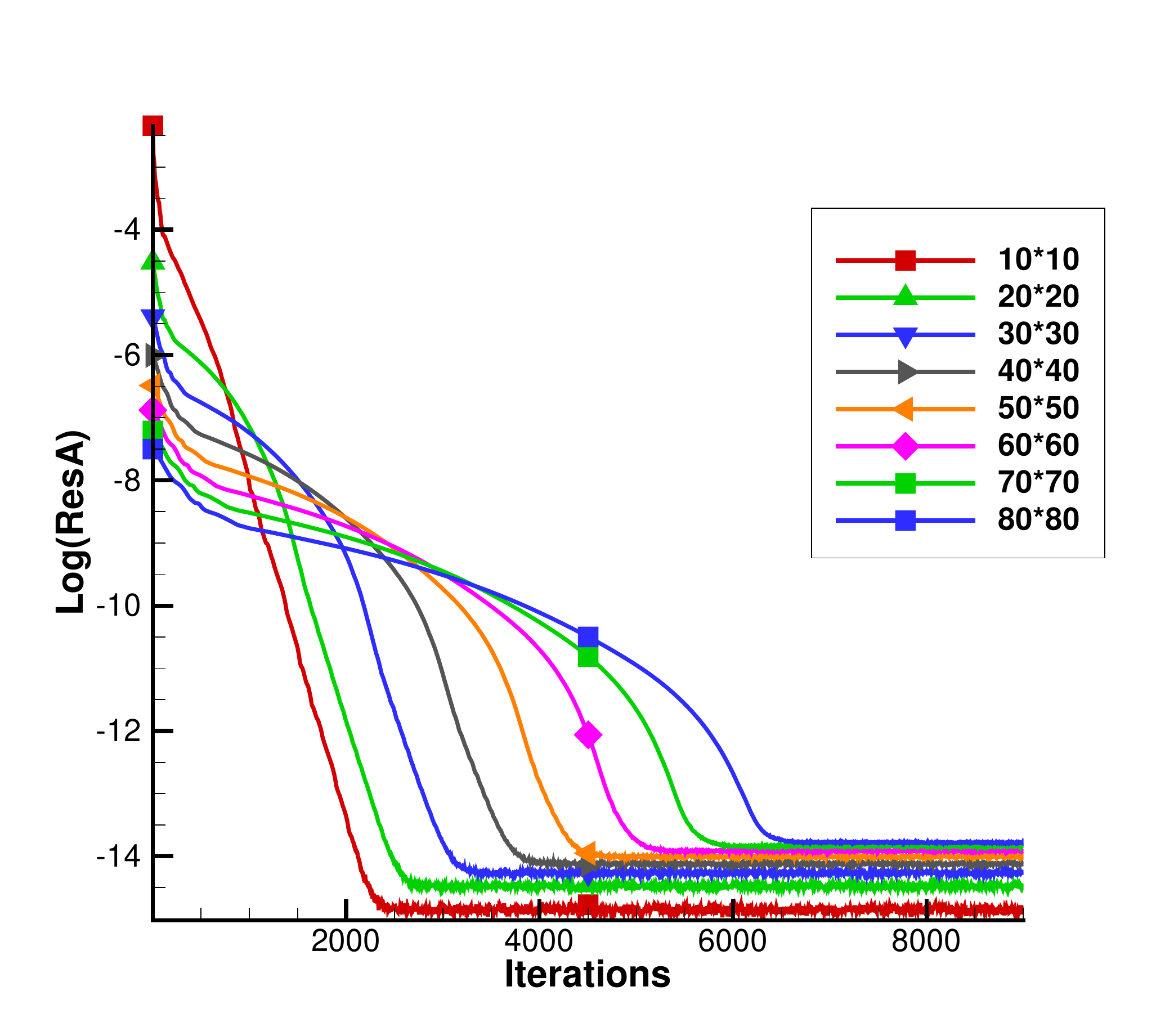}
\end{minipage}%
}%
\subfigure[FE fast sweeping scheme]{
\begin{minipage}[t]{0.25\linewidth}
\centering
\includegraphics[width=1.5in]{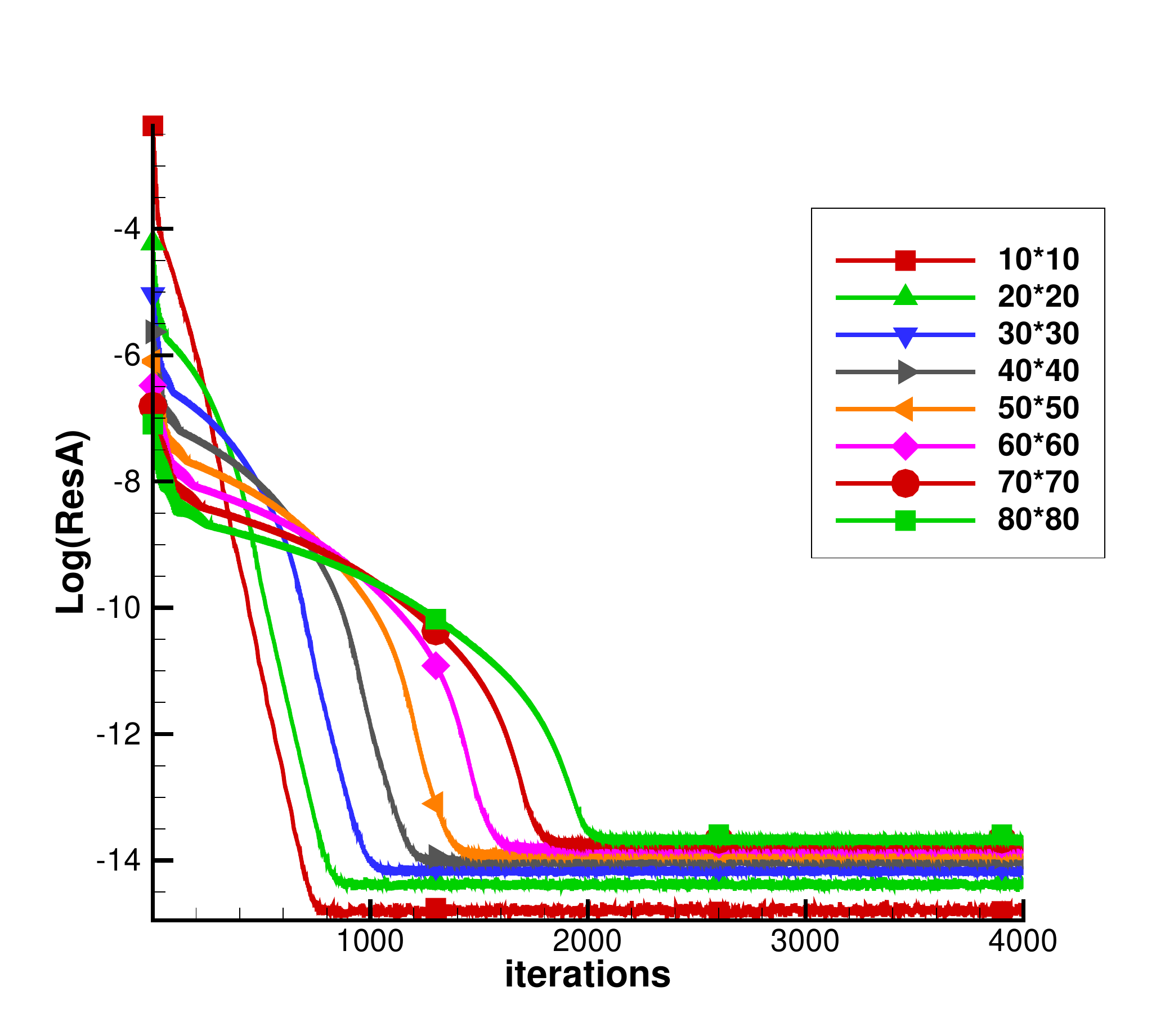}
\end{minipage}%
}%
\centering
\caption{Example 4, A 2D Euler system of equations with source terms. The convergence history of the residue
 as a function of number of iterations for three schemes on different meshes.}
\label{ex4fig}
\end{figure}

\noindent In this example, we solve the following two-dimensional Euler system of equations with source terms:
\[\setlength{\abovedisplayskip}{3pt}\setlength{\belowdisplayskip}{3pt}
\frac{\partial}{\partial{t}}\left(\begin{array}{cccc}
 \rho\\ \rho u\\ \rho v\\E
\end{array}\right)+\frac{\partial}{\partial{x}}
\left(\begin{array}{cccc}
 \rho u\\ \rho u^{2}+p\\ \rho uv\\u(E+p)
\end{array}\right)
+\frac{\partial}{\partial{y}}
\left(\begin{array}{cccc}
 \rho v\\ \rho uv\\ \rho v^{2}+p\\v(E+p)
\end{array}\right)=
\left(\begin{array}{cccc}
$0.4 cos(x+y)$\\ $0.6 cos(x+y)$\\ $0.6 cos(x+y)$\\$1.8 cos(x+y)$
\end{array}\right).
\]
The exact steady-state solutions of the system are $\rho(x,y,\infty)$ =
$1+0.2\sin(x+y), u(x,y,\infty) = 1, v(x,y,\infty) = 1,$ and $p(x,y,\infty) = 1+0.2\sin(x+y)$.
The computational
domain is $(x,y)\in [0,2\pi] \times [0,2\pi],$ and we apply the exact steady-state solutions
on the domain boundaries.
To start the iterations for the schemes, we take
the numerical initial conditions to be the same as the exact steady-state solutions, which
does not satisfy the numerical schemes and will be driven by the iterative schemes to
the numerical steady states. In Table \ref{tab:Margin_settings4}, we report the numerical accuracy for the density variable $\rho$,
iteration numbers and CPU times when these three different iterative schemes
reach the average residue threshold value $10^{-12}$ of the convergence criterion.
It can be observed that all schemes achieve similar numerical errors and the fifth order accuracy when the iterations converge. Comparing these iterative schemes' efficiency,
we obtain the same conclusion as that for the previous problems. The FE Jacobi scheme needs to use a very small CFL number to reach convergence and that leads to large number of iterations and the most CPU time cost among all three
schemes, while
the FE fast sweeping scheme is the most efficient one among them.
We can see that the FE fast sweeping scheme saves more than $50\%$ CPU time cost of that of the RK Jacobi (TVD RK3) scheme.
 In Figure \ref{ex4fig}, the convergence history of the residue (\ref{averres})
 as a function of number of iterations for the FE Jacobi scheme, the RK Jacobi scheme and the FE fast sweeping scheme on different meshes
is shown. We can see that the residue settles down to tiny numbers at the level of round off errors.

\begin{table}
		\centering
		\begin{tabular}{|c|c|c|c|c|c|c|}\hline
			\multicolumn{7}{|c|}{FE Jacobi, $\gamma$=0.1 }\\\hline
			$N\times N$ &L$_{1}$ error&L$_{1}$ order&L$_{\infty} $ error&L$_{\infty} $ order&iter$\sharp$&CPU time\\\hline
    $10\times10$ & 6.74E-04 &   -    & 2.68E-03 &   -    & 5817  & 12.92 \\
    \hline
    $20\times20$ & 1.30E-05 & 5.69  & 3.58E-05 & 6.23  & 6804  & 30.63 \\
    \hline
    $30\times30$ & 1.84E-06 & 4.83  & 4.76E-06 & 4.98  & 8583  & 68.08 \\
    \hline
    $40\times40$ & 4.49E-07 & 4.90   & 1.13E-06 & 5.00    & 10613 & 138.61\\
    \hline
    $50\times50$ & 1.50E-07 & 4.92  & 3.74E-07 & 4.96  & 12725 & 247.23 \\
    \hline
    $60\times60$ & 6.08E-08 & 4.94  & 1.51E-07 & 4.95  & 14931 & 404.63 \\
    \hline
    $70\times70$ & 2.83E-08 & 4.95  & 7.04E-08 & 4.97  & 17068 & 616.05 \\
    \hline
    $80\times80$ & 1.46E-08 & 4.96  & 3.62E-08 & 4.99  & 19093 & 886.75 \\
    \hline
		\end{tabular}
\begin{tabular}{|c|c|c|c|c|c|c|}\hline
			\multicolumn{7}{|c|}{RK Jacobi, $\gamma$=1.0 }\\\hline
			$N\times N$ & L$_{1}$ error&L$_{1}$ order&L$_{\infty} $ error&L$_{\infty} $ order&iter$\sharp$&CPU time\\\hline
    $10\times10$ & 7.41E-04 &   -    & 2.68E-03 &   -    & 1746  & 4.16\\
    \hline
    $20\times20$ & 1.31E-05 & 5.82  & 3.58E-05 & 6.23  & 2037  & 11.06\\
    \hline
    $30\times30$ & 1.85E-06 & 4.84  & 4.76E-06 & 4.98  & 2568  & 26.47 \\
    \hline
    $40\times40$ & 4.51E-07 & 4.91  & 1.13E-06 & 5.00     & 3174  & 52.63 \\
    \hline
    $50\times50$ & 1.50E-07 & 4.93  & 3.74E-07 & 4.96  & 3825  & 94.11\\
    \hline
    $60\times60$ & 6.10E-08 & 4.95  & 1.51E-07 & 4.95  & 4488  & 155.11\\
    \hline
    $70\times70$ & 2.84E-08 & 4.96  & 7.04E-08 & 4.97  & 5130  & 236.75 \\
    \hline
    $80\times80$ & 1.46E-08 & 4.96  & 3.62E-08 & 4.99  & 5739  & 340.70 \\
    \hline
		\end{tabular}
\begin{tabular}{|c|c|c|c|c|c|c|}\hline
			\multicolumn{7}{|c|}{FE fast sweeping, $\gamma$=1.0 }\\\hline
			$N\times N $&L$_{1}$ error&L$_{1}$ order&L$_{\infty} $ error&L$_{\infty} $ order&iter$\sharp$&CPU time\\\hline
    $10\times10$ & 6.62E-04 &   -    & 2.68E-03 &   -    & 560   & 1.26 \\
    \hline
    $20\times20$ & 1.30E-05 & 5.67  & 3.58E-05 & 6.23  & 653   & 4.02 \\
    \hline
    $30\times30$ & 1.84E-06 & 4.83  & 4.76E-06 & 4.98  & 821   & 10.69\\
    \hline
    $40\times40$ & 4.49E-07 & 4.90   & 1.13E-06 & 5.00     & 1010  & 22.53\\
    \hline
    $50\times50$ & 1.50E-07 & 4.92  & 3.74E-07 & 4.96  & 1213  & 42.39 \\
    \hline
    $60\times60$ & 6.08E-07 & 4.94  & 1.51E-07 & 4.95  & 1421  & 71.19 \\
    \hline
    $70\times70$ & 2.83E-08 & 4.95  & 7.04E-08 & 4.97  & 1622  & 110.55 \\
    \hline
    $80\times80$ & 1.46E-08 & 4.96  & 3.62E-08 & 4.99  & 1814  & 160.91 \\
    \hline
		\end{tabular}
		\caption{Example 4, A 2D Euler system of equations with source terms. Accuracy, iteration numbers and CPU times of three different iterative schemes. CPU time unit: second}
		\label{tab:Margin_settings4}
	\end{table}

\bigskip
\noindent{\bf Example 5. A 2D Euler system of equations without source terms}

\begin{figure}
\centering
\subfigure[FE Jacobi scheme]{
\begin{minipage}[t]{0.25\linewidth}
\centering
\includegraphics[width=1.5in]{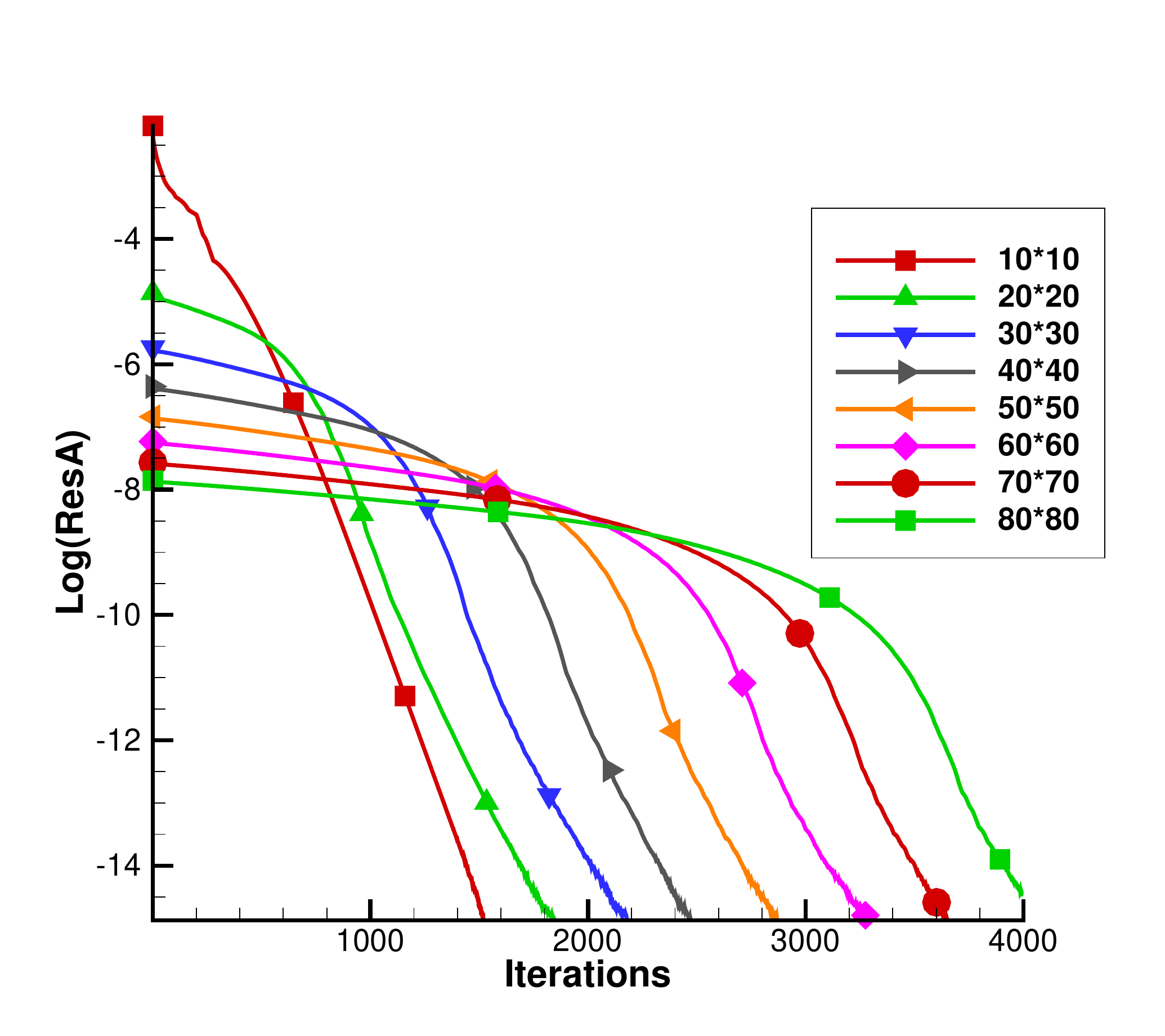}
\end{minipage}%
}%
\subfigure[RK Jacobi scheme]{
\begin{minipage}[t]{0.25\linewidth}
\centering
\includegraphics[width=1.5in]{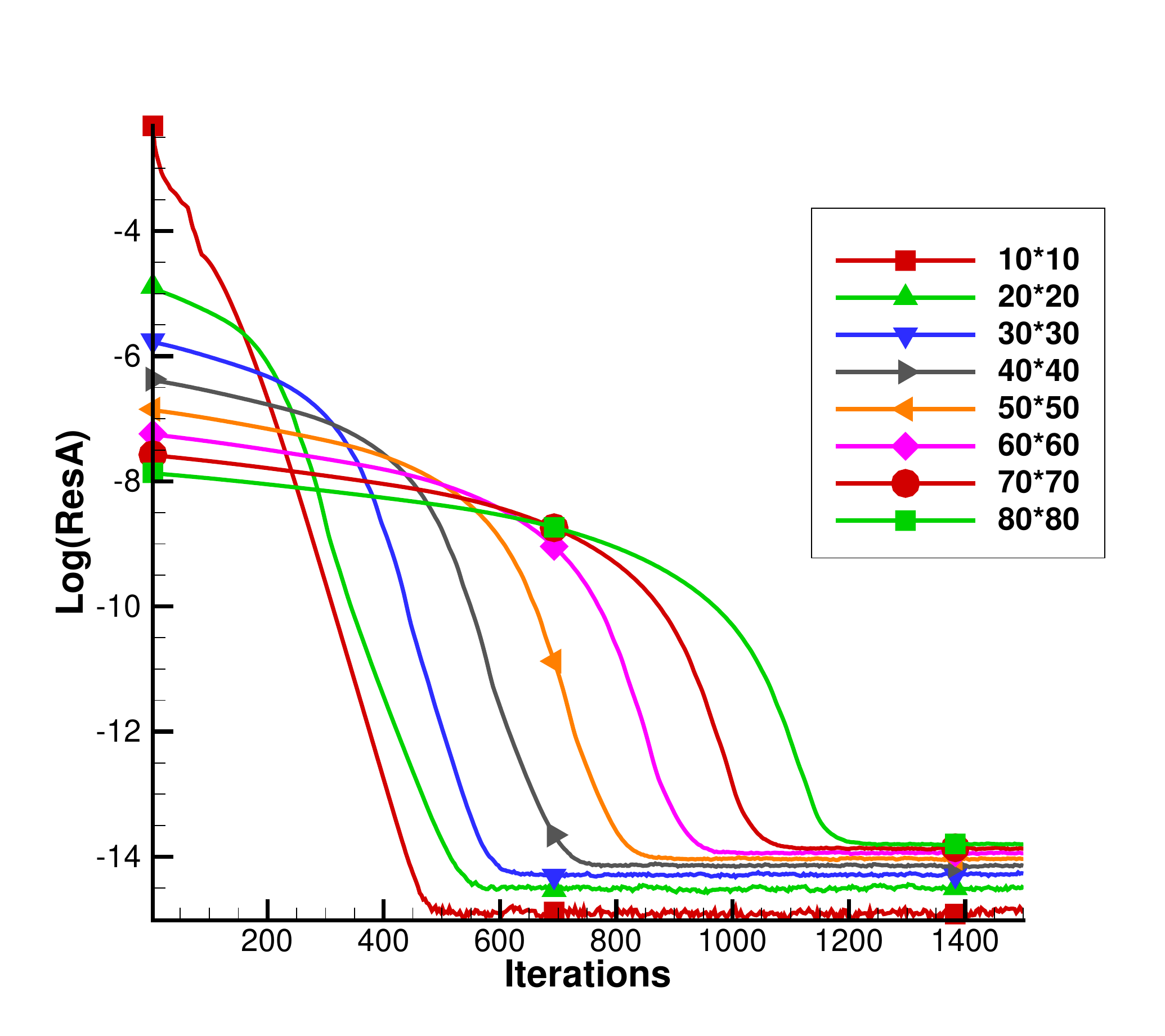}
\end{minipage}%
}%
\subfigure[FE fast sweeping scheme]{
\begin{minipage}[t]{0.25\linewidth}
\centering
\includegraphics[width=1.5in]{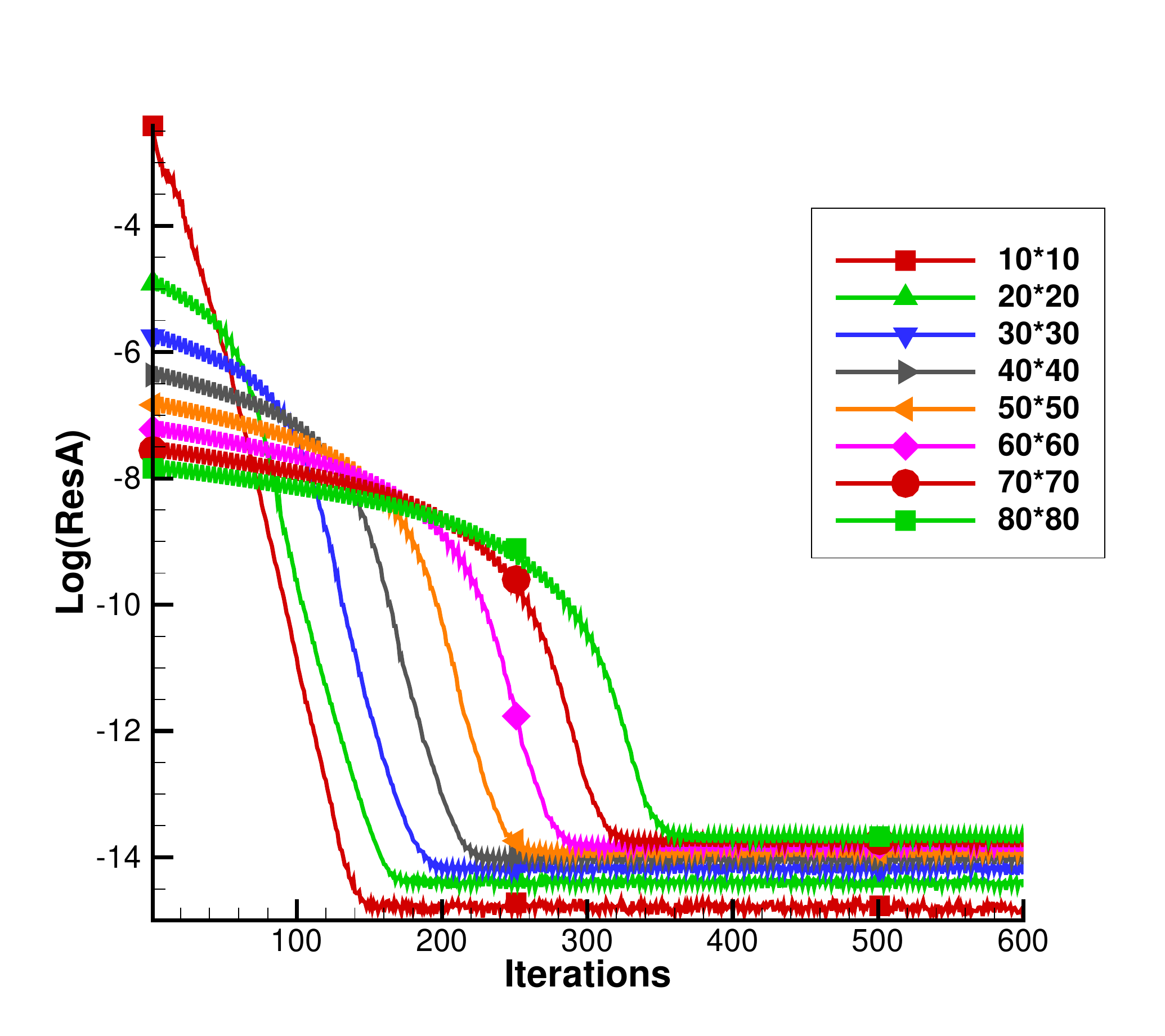}
\end{minipage}%
}%
\centering
\caption{Example 5, A 2D Euler system of equations without source terms. The convergence history of the residue
 as a function of number of iterations for three schemes on different meshes.}
 \label{ex5fig}
\end{figure}

\noindent In this example, we consider the following two-dimensional
Euler system of equations
\[\setlength{\abovedisplayskip}{3pt}\setlength{\belowdisplayskip}{3pt}
\frac{\partial}{\partial{t}}\left(\begin{array}{cccc}
 \rho\\ \rho u\\ \rho v\\E
\end{array}\right)+\frac{\partial}{\partial{x}}
\left(\begin{array}{cccc}
 \rho u\\ \rho u^{2}+p\\ \rho uv\\u(E+p)
\end{array}\right)
+\frac{\partial}{\partial{y}}
\left(\begin{array}{cccc}
 \rho v\\ \rho uv\\ \rho v^{2}+p\\v(E+p)
\end{array}\right)=
0
\]
with the exact steady-state solutions $\rho(x,y,\infty)=1+0.2\sin(x-y), u(x,y,\infty)=1,
v(x,y,\infty)=1,$ and $p(x,y,\infty)=1$.
The computational domain is $(x,y) \in [0,2\pi] \times [0,2\pi]$.
Similar as the last example, we apply the exact steady-state solutions
on the domain boundaries, and take
the numerical initial conditions to be the same as the exact steady-state solutions in the iterative schemes.
In Table \ref{tab:Margin_settings5}, we show the numerical accuracy for the density variable $\rho$,
iteration numbers and CPU times when these three different iterative schemes
reach the average residue threshold value $10^{-12}$ of the convergence criterion.
Similar conclusions can be drawn as the last example. All schemes have similar numerical errors and the fifth order accuracy when the iterations converge. However, The FE Jacobi scheme has to use a very small CFL number to reach convergence, which leads to large number of iterations and the most CPU time cost on refined meshes among all three
schemes. Again, the FE fast sweeping scheme is the most efficient one among them, and it
only takes about $40\%$ CPU time costs of that of the RK Jacobi (TVD RK3) scheme.
In Figure \ref{ex5fig}, the convergence history of the residue (\ref{averres})
 as a function of number of iterations for these three schemes on different meshes
is shown. It is verified that the residue of iterations settles down to tiny numbers at the level of round off errors for all three schemes.

\begin{table}
		\centering
		\begin{tabular}{|c|c|c|c|c|c|c|}\hline
			\multicolumn{7}{|c|}{FE Jacobi, $\gamma$=0.1 }\\\hline
			N$\times$N&L$_{1}$ error&L$_{1}$ order&L$_{\infty} $ error&L$_{\infty} $ order&iter$\sharp$&CPU time\\\hline

    10$\times$10 & 1.68E-03 &   -    & 8.01E-03 &   -    & 1233  & 1.11 \\
    \hline
    20$\times$20 & 2.38E-05 & 6.14  & 1.39E-04 & 5.85  & 1393  & 4.28 \\
    \hline
    30$\times$30 & 3.38E-06 & 4.82  & 1.93E-05 & 4.86  & 1684  & 10.97 \\
    \hline
    40$\times$40 & 8.29E-07 & 4.88  & 4.58E-06 & 5.00     & 2033  & 22.89 \\
    \hline
    50$\times$50 & 2.77E-07 & 4.91  & 1.52E-06 & 4.94  & 2410  & 42.27\\
    \hline
    60$\times$60 & 1.13E-07 & 4.93  & 6.14E-07 & 4.97  & 2803  & 69.78 \\
    \hline
    70$\times$70 & 5.27E-08 & 4.94  & 2.85E-07 & 4.98  & 3219  & 108.63 \\
    \hline
    80$\times$80 & 2.72E-08 & 4.95  & 1.46E-07 & 4.98  & 3628  & 158.42\\
    \hline
\end{tabular}
\begin{tabular}{|c|c|c|c|c|c|c|}\hline
			\multicolumn{7}{|c|}{RK Jacobi, $\gamma$=1.0 }\\\hline
			N$\times$N&L$_{1}$ error&L$_{1}$ order&L$_{\infty} $ error&L$_{\infty} $ order&iter$\sharp$&CPU time\\\hline

    10$\times$10 & 1.85E-03 &   -    & 8.01E-03 &   -    & 378   & 1.17 \\
    \hline
    20$\times$20 & 2.45E-05 & 6.24  & 1.39E-04 & 5.85  & 426   & 2.43\\
    \hline
    30$\times$30 & 3.44E-06 & 4.84  & 1.93E-05 & 4.86  & 504   & 5.12 \\
    \hline
    40$\times$40 & 8.40E-07 & 4.90   & 4.58E-06 & 5.00     & 618   & 10.06 \\
    \hline
    50$\times$50 & 2.80E-07 & 4.92  & 1.52E-06 & 4.94  & 729   & 18.28 \\
    \hline
    60$\times$60 & 1.14E-07 & 4.94  & 6.14E-07 & 4.97  & 852   & 29.06 \\
    \hline
    70$\times$70 & 5.30E-08 & 4.95  & 2.85E-07 & 4.98  & 972   & 44.78 \\
    \hline
    80$\times$80 & 2.74E-08 & 4.96  & 1.46E-07 & 4.98  & 1095  & 64.31\\
    \hline
\end{tabular}
\begin{tabular}{|c|c|c|c|c|c|c|}\hline
			\multicolumn{7}{|c|}{FE fast sweeping, $\gamma$=1.0 }\\\hline
			N$\times$N&L$_{1}$ error&L$_{1}$ order&L$_{\infty} $ error&L$_{\infty} $ order&iter$\sharp$&CPU time\\\hline

   10$\times$10 & 1.65E-03 &   -    & 8.01E-03 &   -    & 112   & 0.14\\
    \hline
    20$\times$20 & 2.37E-05 & 6.12  & 1.39E-04 & 5.85  & 130   & 0.66 \\
    \hline
    30$\times$30 & 3.37E-06 & 4.81  & 1.93E-05 & 4.86  & 155   & 1.78\\
    \hline
    40$\times$40 & 8.28E-07 & 4.88  & 4.58E-06 & 5.00     & 186   & 3.84 \\
    \hline
    50$\times$50 & 2.77E-07 & 4.91  & 1.52E-06 & 4.94  & 220   & 7.17 \\
    \hline
    60$\times$60 & 1.13E-07 & 4.93  & 6.14E-07 & 4.97  & 254   & 11.98 \\
    \hline
    70$\times$70 & 5.26E-08 & 4.94  & 2.85E-07 & 4.98  & 290   & 18.72 \\
    \hline
    80$\times$80 & 2.72E-08 & 4.95  & 1.46E-07 & 4.98  & 327   & 28.27 \\
    \hline
		\end{tabular}
		\caption{Example 5, A 2D Euler system of equations without source terms. Accuracy, iteration numbers and CPU times of three different iterative schemes. CPU time unit: second}
		\label{tab:Margin_settings5}
	\end{table}

\subsection{Examples with non-smooth solutions}

In this section, We test these schemes using problems with discontinuous solutions.

\bigskip
\noindent{\bf Example 6. One-dimensional steady shock}

\noindent In this example we
consider a 1D steady shock problem (\cite{WuLiang,SSCW,SCW})
$$U_{t}+F(U)_{x}=0,$$
where $U=(\rho,\rho u, e)^{T}, F(U)=(\rho u,\rho u^{2}+p,u(e+p))^{T}$. Here $\rho, u, e$ are the density, velocity, and total energy respectively. $p$ is the pressure and it is related to the total energy by $e=\frac{p}{\gamma'-1}+\frac{1}{2}\rho u^{2}$, $\gamma'=1.4$ which is the ratio of specific heat.
The problem is defined on the domain $x\in[-1,1]$.  Initially the flow Mach number at the left of the shock wave is $M_{\infty} =
2$ with the shock wave located at $x = 0$. Periodic boundary conditions are applied here. The initial
condition to start the iterations in the schemes is given by the Rankine-Hugoniot jump condition
(see e.g. $\cite{M.A.Sa}$) as the following:
\[U(x,0)=\begin{cases}
U_{l}, \qquad         \text { if } x<0;\\
U_{r},  \qquad        \text { if } x>0,
\end{cases}\]
where
\[
\left(\begin{array}{ccc}
 p_{l} \\\rho_{l}\\u_{l}
\end{array}\right)=
\left(\begin{array}{ccc}
 \frac{1}{\gamma' M_{\infty}^{2}}\\1\\1
\end{array}\right),
\qquad\left(\begin{array}{ccc}
 p_{r} \\\rho_{r}\\u_{r}
\end{array}\right)=
\left(\begin{array}{ccc}
 p_{l}\frac{2\gamma' M_{\infty}^{2}-(\gamma'-1)}{\gamma'+1} \\\rho_{l}\frac{\frac{\gamma'+1}{\gamma'-1}\frac{p_{r}}{p_{l}}+1}{\frac{\gamma'+1}{\gamma'-1}+\frac{p_{r}}{p_{l}}}\\\sqrt{\gamma' \frac{(2+(\gamma'-1)M_{\infty}^{2})p_{r}}{(2\gamma' M_{\infty}^{2}+(1-\gamma'))\rho_{r}}}
\end{array}\right).
\]
In this example the initial condition is the same as the exact solution of the steady state. Similar as the previous examples, because the exact steady state solution of the PDE does not satisfy the numerical schemes, we will be able to see the convergence behavior that the iterative schemes drive the initial condition to the numerical steady states.

\begin{figure}
\centering
\subfigure[RK Jacobi scheme]{
\begin{minipage}[t]{0.4\linewidth}
\centering
\includegraphics[width=2in]{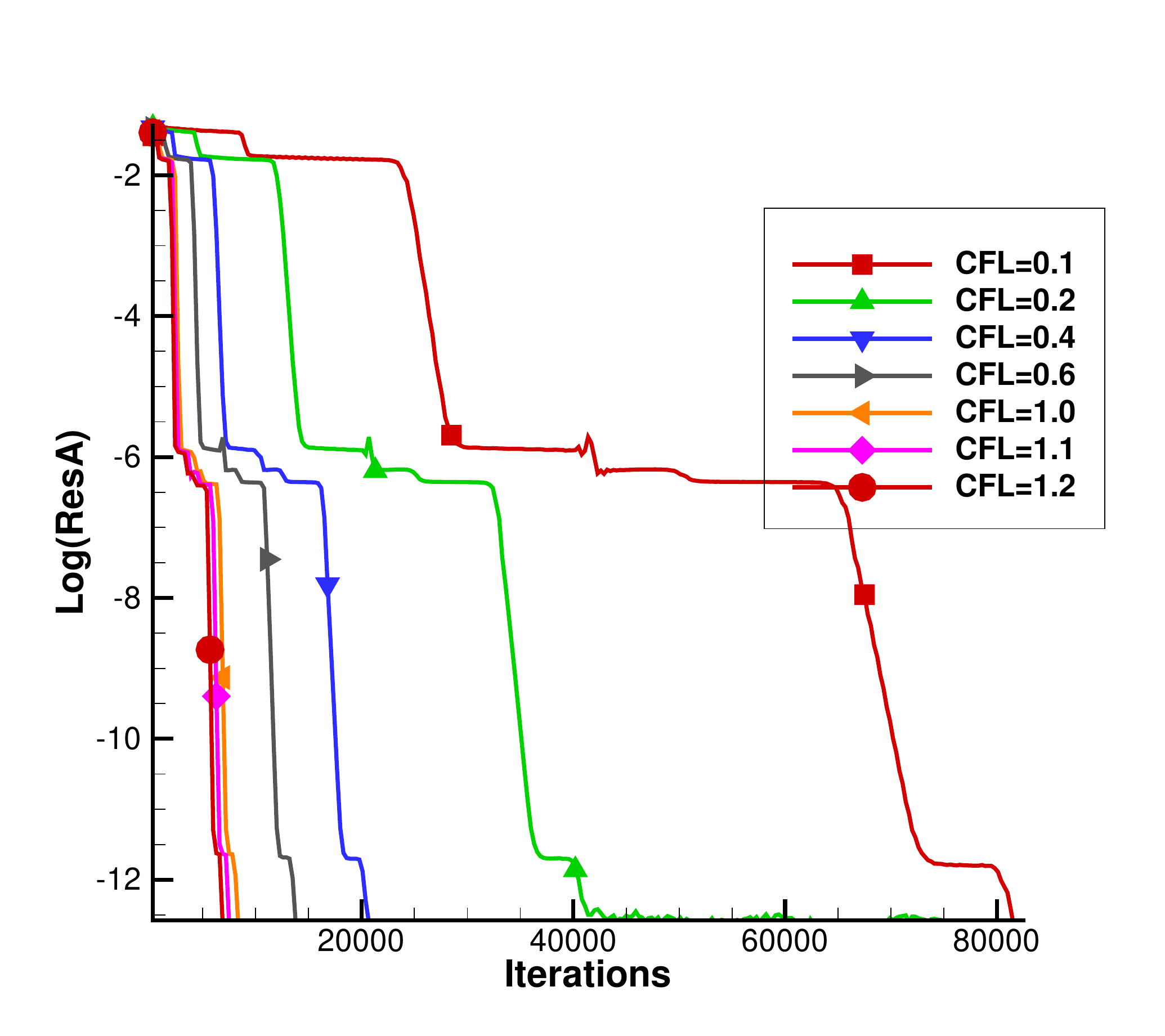}
\end{minipage}%
}%
\subfigure[FE fast sweeping scheme]{
\begin{minipage}[t]{0.4\linewidth}
\centering
\includegraphics[width=2in]{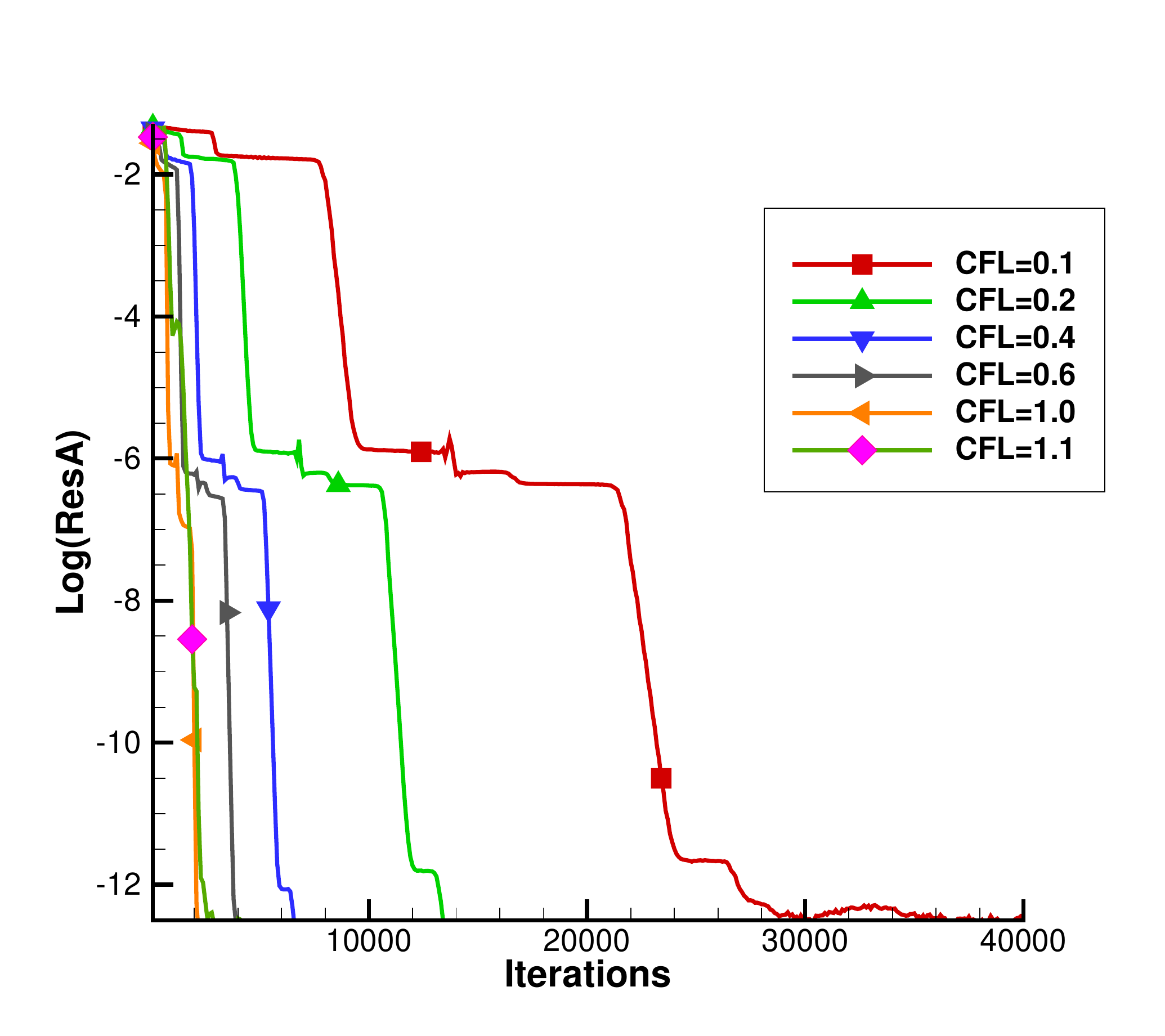}
\end{minipage}%
}%
\centering
\caption{Example 6, One-dimensional steady shock. The convergence history of the residue
 as a function of number of iterations for two schemes with different CFL numbers.}
 \label{ex6fig}
\end{figure}

In this example, $400$ uniformly spaced mesh
points are used in the spatial discretization. In Table \ref{tab:Margin_settings6},
we report number of iterations required to achieve
convergence (i.e. to reach the convergence threshold value $10^{-12}$), the final time and total CPU
time when convergence is obtained for these three iterative schemes with different
CFL numbers $\gamma$.
The FE Jacobi iterative scheme requires
a very small CFL number $\gamma$ to converge. In this example, numerical tests show that $\gamma$ needs to be less than
or equal to $0.09$. If $\gamma=0.1$, the residue hangs at ${10^{-10.9}}$ and
does not decrease till the pre-set maximum iteration number $100,000$ is reached.
The RK Jacobi scheme converges with much larger CFL numbers. Here, if $\gamma$ is less than or
equal to $1.2$, the convergence criterion $ResA < 10^{-12}$ can be satisfied. If $\gamma$ is increased to $1.3$, the residue stops decreasing at $10^{-2.41}$.
The FE fast sweeping method converges with similar CFL numbers as the RK Jacobi scheme. Here
when $\gamma$ is less than or equal to $1.1$, the FE fast sweeping scheme converges. However,
the FE fast sweeping method can converge to reach the residue threshold value $10^{-12}$ with much less
iteration numbers and CPU times than the RK Jacobi scheme. In this example, with the smallest iteration numbers
that each method can achieve to reach the convergence, the FE fast sweeping method has the CPU time cost $4.42$ seconds while the
RK Jacobi scheme needs $7.28$ seconds, hence about $40\%$ CPU time can be saved.
In Figure \ref{ex6fig}, the convergence history of the residue (\ref{averres})
 as a function of number of iterations for the RK Jacobi scheme and the FE fast sweeping scheme with different CFL numbers
is shown. Again, we see that the residue of iterations settles down to tiny numbers at the level of round off errors.

In this one-dimensional problem whose steady state has shock wave, we draw the same conclusion as the previous examples with smooth solutions. The FE fast sweeping
method is the most efficient approach for the fifth order multi-resolution
WENO computations of the steady state problems among the three methods discussed in the paper,
in terms of both iteration number and CPU time. This is further verified by the
more complicated two dimensional simulations of Euler systems in the following examples.

\begin{table}
		\centering
		\begin{tabular}{|c|c|c|c|}\hline
			\multicolumn{4}{|c|}{FE Jacobi scheme}\\\hline
            $\gamma:$ CFL number & iteration number & final time & CPU time \\\hline
            0.1  & Not convergent & - & - \\\hline
		\end{tabular}
\begin{tabular}{|c|c|c|c|}\hline
			\multicolumn{4}{|c|}{RK Jacobi scheme }\\\hline
            $\gamma:$ CFL number & iteration number & final time & CPU time \\\hline
            0.1   & 80355 & 8.93  & 83.20 \\
    \hline
    0.2   & 40515 & 9.00     & 43.17 \\
    \hline
    0.4   & 20262 & 9.00       & 22.25\\
    \hline
    1.0     & 8118  & 9.02  & 8.97 \\
    \hline
    1.1   & 7380  & 9.02  & 7.73 \\
    \hline
    1.2   & 6765  & 9.02  & 7.28 \\
    \hline
    1.3   & Not convergent &   -    & - \\
    \hline
		\end{tabular}
\begin{tabular}{|c|c|c|c|}\hline
			\multicolumn{4}{|c|}{FE fast sweeping scheme }\\\hline
            $\gamma:$ CFL number & iteration number & final time & CPU time \\\hline
           0.1   & 26904 & 8.97  & 52.59\\
    \hline
    0.2   & 13164 & 8.77  & 25.66 \\
    \hline
    0.4   & 5894  & 7.86  & 11.81\\
    \hline
    0.6   & 3776  & 7.55  & 7.81 \\
    \hline
    1.0     & 2088  & 6.96  & 4.42 \\
    \hline
    1.1   & 2426  & 8.89  & 4.50 \\
    \hline
     1.2   & Not convergent &    -   & - \\
    \hline
		\end{tabular}
		\caption{Example 6. 1D steady shock. Number of iterations, the final time and total CPU time when convergence is obtained. Convergence criterion threshold value is $10^{-12}$. CPU time unit: second}
		\label{tab:Margin_settings6}
	\end{table}

\bigskip
\noindent{\bf Example 7. A two-dimensional oblique steady shock}

\noindent In this example, we use these three iterative schemes to solve a two-dimensional
oblique steady shock problem, which is a benchmark steady state problem (see e.g. \cite{WuLiang,SSCW,SCW}). The shock has an
angle of $135^{\circ}$ with the positive x-direction. The Mach number of the flow at the left of the shock
is $M_{\infty} = 2$. The computational domain is $0\leq x\leq4$ and $0\leq y \leq2$. Initially the oblique
shock wave passes the point $(3, 0)$. The computational mesh has $200 \times 100$ equally spaced grid points
with $ \bigtriangleup x = \bigtriangleup y$. Periodic boundary conditions are applied.
The convergence criterion threshold value is set to be $10^{-12}$.

\begin{figure}
\centering
\subfigure[RK Jacobi scheme]{
\begin{minipage}[t]{0.5\linewidth}
\centering
\includegraphics[width=2.1in]{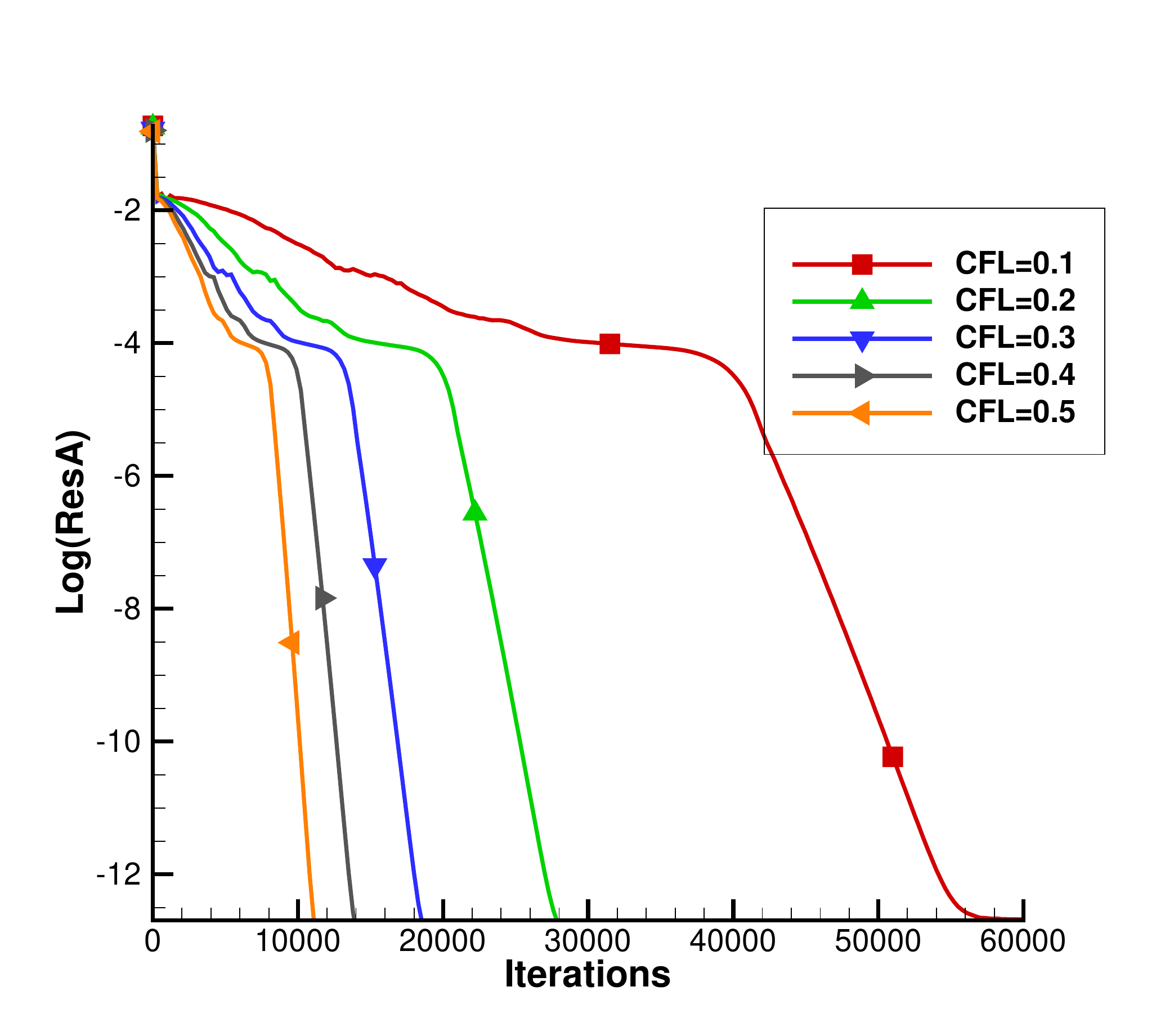}
\end{minipage}%
}%
\subfigure[FE fast sweeping scheme]{
\begin{minipage}[t]{0.5\linewidth}
\centering
\includegraphics[width=2.1in]{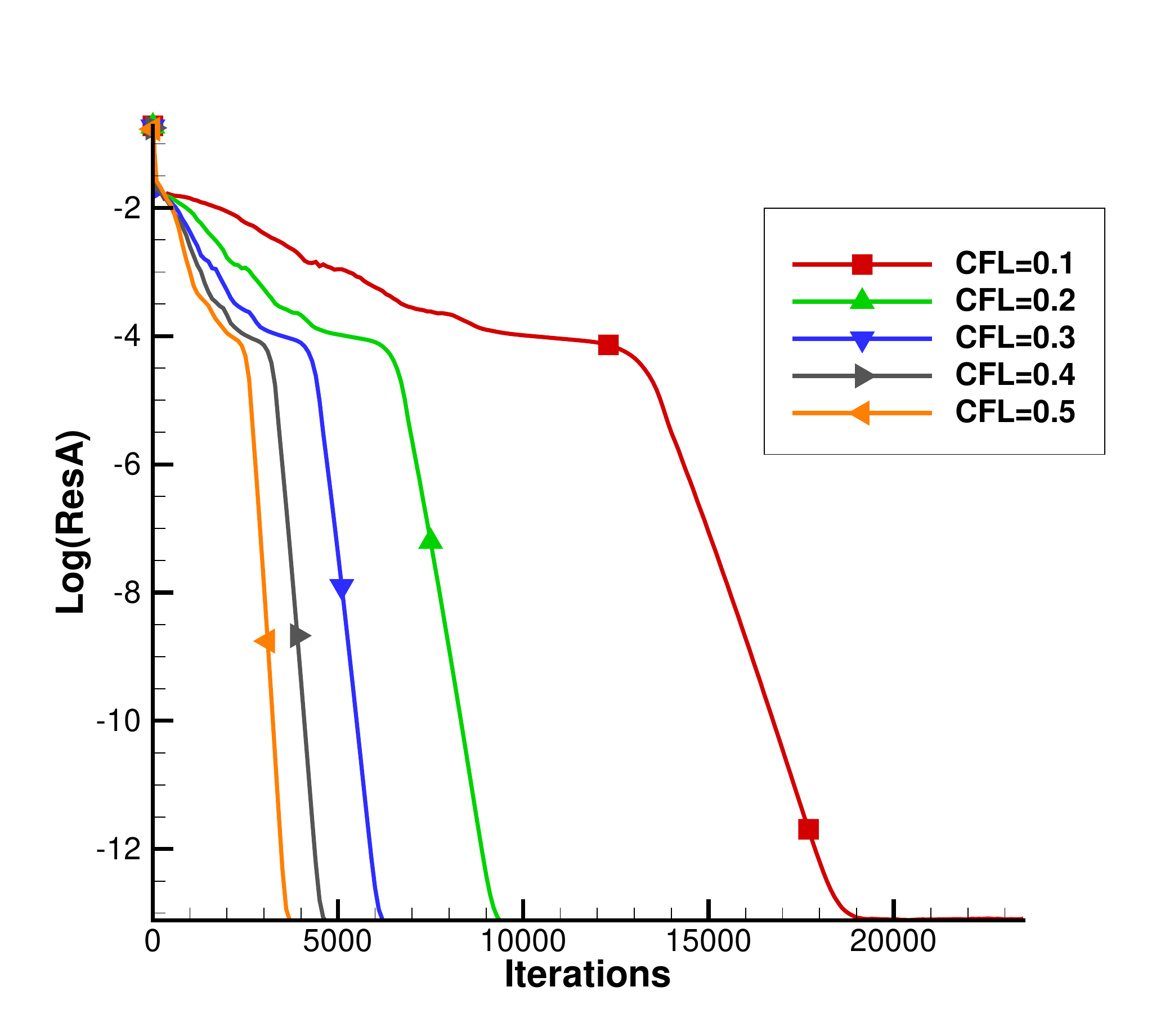}
\end{minipage}%
}%
\centering
\caption{Example 7, A two-dimensional oblique steady shock. The convergence history of the residue
 as a function of number of iterations for two schemes with different CFL numbers.}
\label{ex7fig1}
\end{figure}

\begin{figure}
\centering
\subfigure[RK Jacobi scheme]{
\begin{minipage}[t]{0.5\linewidth}
\centering
\includegraphics[width=2.1in]{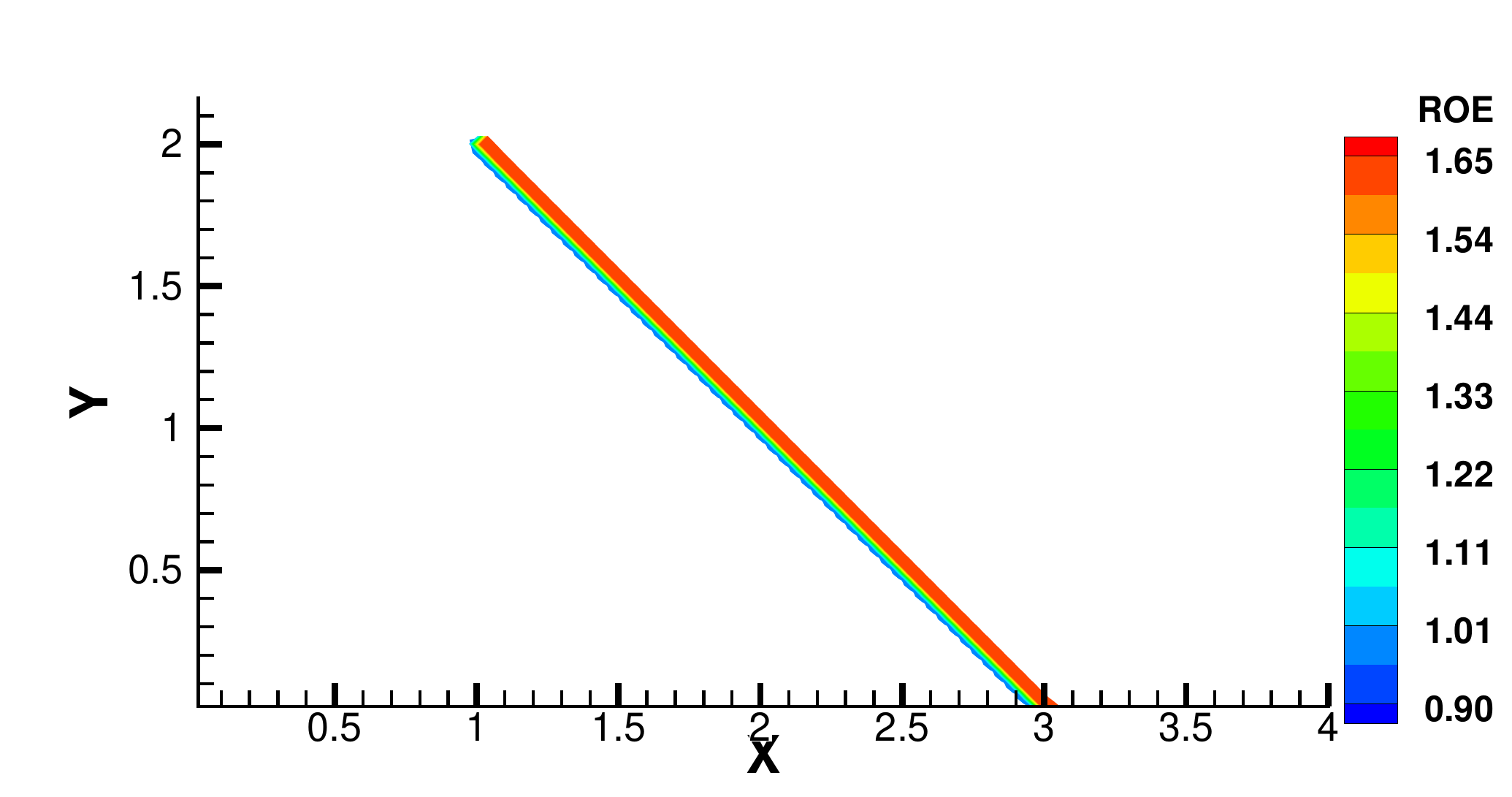}
\end{minipage}%
}%
\subfigure[FE fast sweeping scheme]{
\begin{minipage}[t]{0.5\linewidth}
\centering
\includegraphics[width=2.1in]{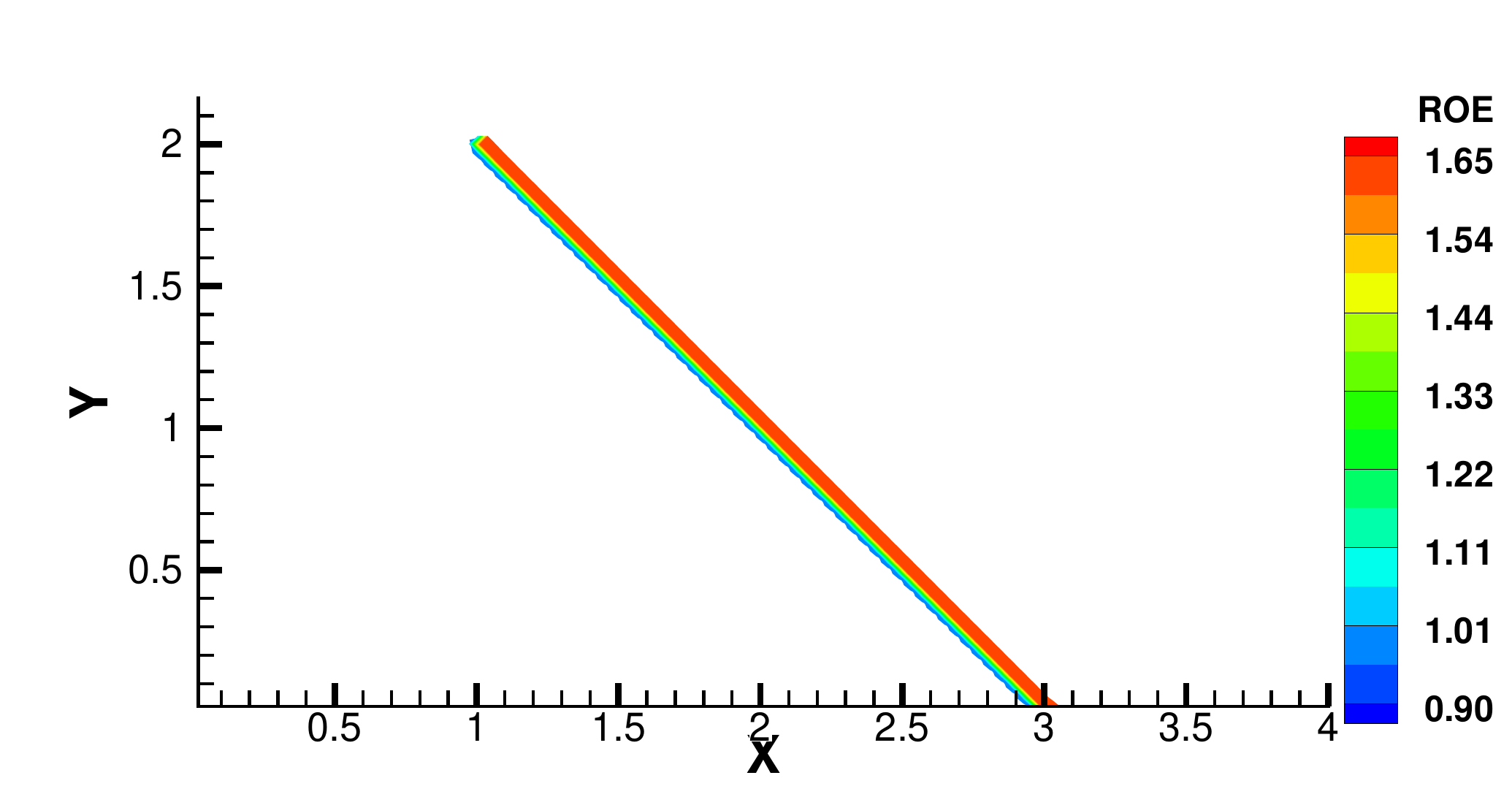}
\end{minipage}%
}%
\centering
\caption{Example 7, A two-dimensional oblique steady shock. 30 equally spaced density contours from 1.05 to 1.65
of the converged steady states of numerical solutions by two different iterative schemes.}
\label{ex7fig2}
\end{figure}

In Table \ref{tab:Margin_settings7}, number of iterations required
to reach the convergence threshold value $10^{-12}$, the final time and total CPU
time when convergence is obtained for these three iterative schemes with different
CFL numbers $\gamma$ are reported. As that in the 1D steady shock problem, the FE Jacobi scheme
requires a very small CFL number $\gamma=0.1$ to achieve the convergence, hence it needs a large
amount of iteration numbers and CPU time. The RK Jacobi scheme relaxes the CFL number constraint and converges with much larger CFL numbers
up to $\gamma=0.5$, which reduces the CPU time from $2817.84$ seconds of the FE Jacobi scheme to
$1718.41$ seconds.  The FE fast sweeping method converges with similar CFL numbers as the RK Jacobi scheme, and it is more efficient than the RK Jacobi scheme with much reduced iteration numbers and CPU times. In fact, by using the largest possible CFL number to reach the convergence, the FE fast sweeping method only needs CPU time $987.02$ seconds,
which saves about $43\%$ CPU time of the simulation by the TVD RK3 scheme (RK Jacobi). In Figure \ref{ex7fig1},
residue history in terms of iterations for the RK Jacobi and the FE fast sweeping schemes with various CFL numbers is
shown. We observe that the residue of iterations settles down to tiny numbers at the level of round off errors for both schemes. In Figure \ref{ex7fig2}, contour plots of the density variable of the converged steady state solutions of the RK Jacobi and the FE fast sweeping schemes are shown, and it is verified that the similar results are obtained.
In summary, we draw the consistent conclusion in this 2D steady shock wave example with that in the last 1D steady shock wave example, i.e., the FE fast Sweeping scheme is
the most efficient one  among these three iterative schemes.

\begin{table}
		\centering
		\begin{tabular}{|c|c|c|c|}\hline
			\multicolumn{4}{|c|}{FE Jacobi scheme }\\\hline
            $\gamma:$ CFL number & iteration number & final time & CPU time \\\hline
            0.1	&18391	&17.25&	2817.84\\\hline

            0.2  & Not convergent & - & - \\\hline
		\end{tabular}
\begin{tabular}{|c|c|c|c|}\hline
			\multicolumn{4}{|c|}{RK Jacobi scheme}\\\hline
            $\gamma:$ CFL number & iteration number & final time & CPU time \\\hline
            0.3&	19158&	17.99	&2895.34 \\
    \hline
0.4	&14244&	17.83	&2151.48 \\
    \hline
0.5	&11370	&17.79&	1718.41 \\
    \hline
    0.6   & Not convergent &    -   & - \\
    \hline
		\end{tabular}
\begin{tabular}{|c|c|c|c|}\hline
			\multicolumn{4}{|c|}{FE fast sweeping scheme }\\\hline
            $\gamma:$ CFL number & iteration number & final time & CPU time \\\hline
0.3	&6105	&17.20&	1668.17\\
    \hline
0.4	&4541	&17.06	&1244.56\\
    \hline
0.5	&3601	&16.90	&987.02\\
    \hline
     0.6   & Not convergent &   -    & - \\
    \hline
		\end{tabular}
		\caption{Example 7. A two-dimensional oblique steady shock. Number of iterations, the final time and total CPU time when convergence is obtained. Convergence criterion threshold value is $10^{-12}$. CPU time unit: second}
		\label{tab:Margin_settings7}
	\end{table}

\bigskip
\noindent{\bf Example 8. Regular shock reflection}

\noindent In this example, we solve the two dimensional regular shock reflection problem by using the discussed iterative
schemes. This problem is a typical and difficult benchmark problem of using high order schemes to simulate steady
flow. As that found in \cite{SSCW,SCW}, even with advanced techniques to improve the steady state convergence,
it is still difficult for the residue of high order WENO schemes to converge to the level of round off errors.
The fifth order fixed-point fast sweeping WENO scheme in \cite{WuLiang} is also unable to drive the iteration residue to the level of round off errors. However, with the help of the unequal-sized sub-stencil WENO method such as the multi-resolution WENO scheme \cite{JUNZ, JUNZ3}, the absolutely convergent fixed-point fast sweeping method proposed here
can make the iteration residue settle
down to the round off error level easily, and at the same time maintain the high efficiency of the original
fixed-point fast sweeping WENO scheme in \cite{WuLiang}. This will be verified in the following numerical results.

Here the computational domain is a rectangle with the length
$4$ and the height $1$. The boundary conditions in this problem consist of a reflection condition along the bottom
boundary, supersonic outflow along the right boundary and Dirichlet conditions on the other
two sides:
\begin{equation*}
(\rho,u,v,p)^{T}=
\begin{cases}
(1.0,2.9,0,5/7)^{T}\mid_{(0,y,t)^{T}},\\
(1.69997,2.61934,-0.50632,1.52819)^{T}\mid_{(x,1,t)^{T}}.\\
\end{cases}
\end{equation*}
The initial values at points in the entire domain are the same as those at the left boundary.
The convergence criterion threshold value is set to be $10^{-12}$. The computational grid is $120 \times 30$.

\begin{figure}
\centering
\subfigure[RK Jacobi scheme]{
\begin{minipage}[t]{0.5\linewidth}
\centering
\includegraphics[width=2.1in]{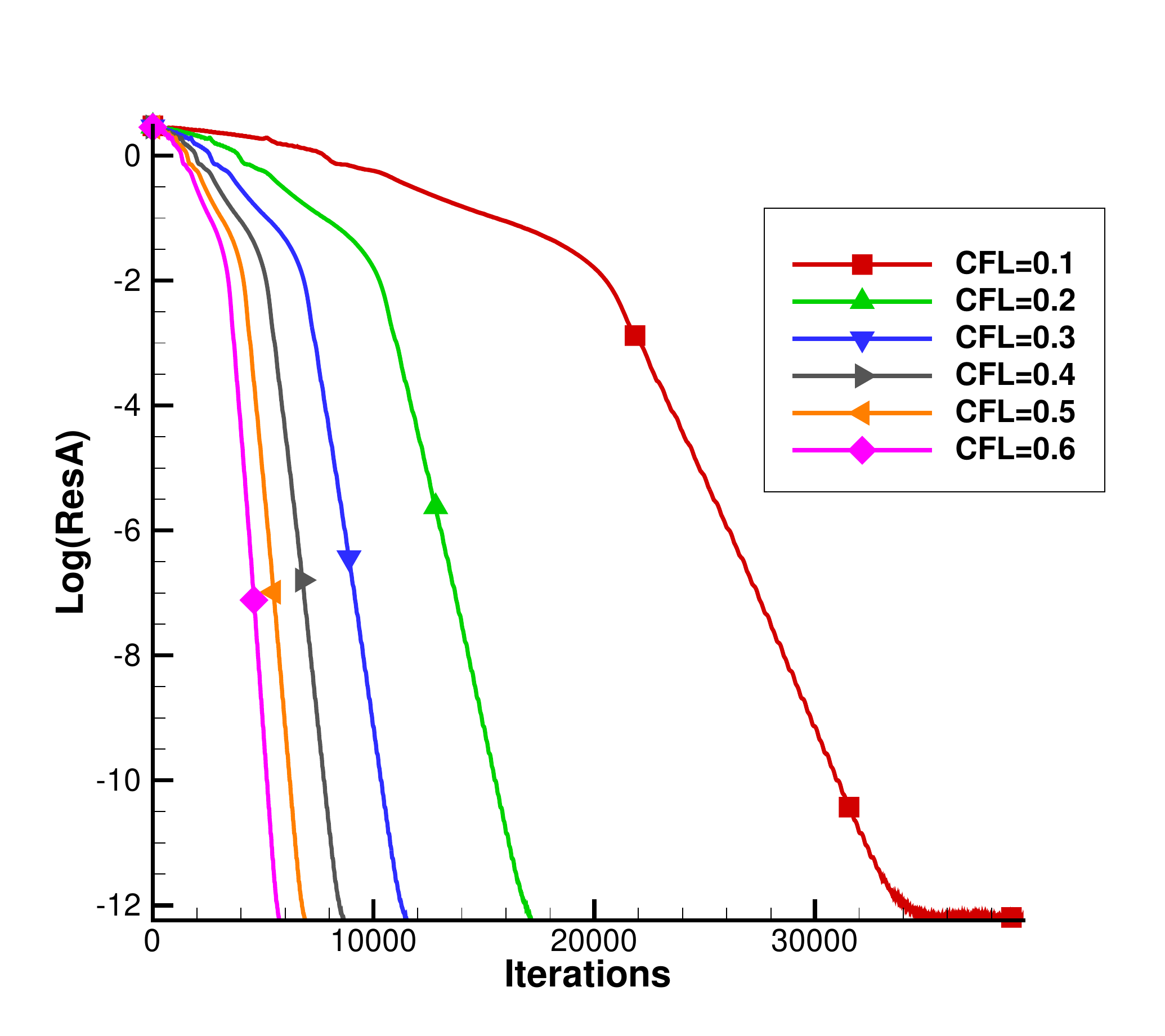}
\end{minipage}%
}%
\subfigure[FE fast sweeping scheme]{
\begin{minipage}[t]{0.5\linewidth}
\centering
\includegraphics[width=2.1in]{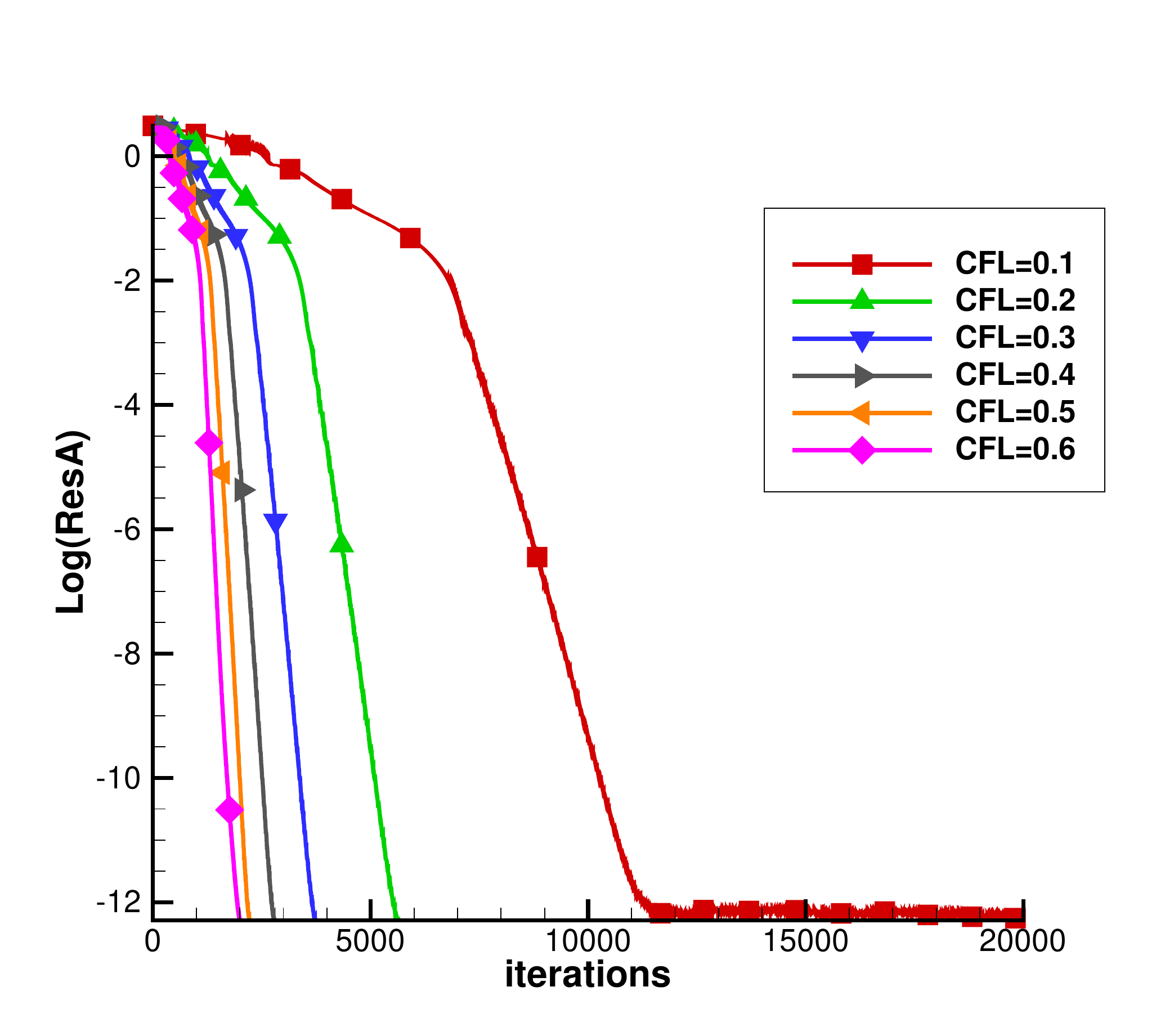}
\end{minipage}%
}%
\centering
\caption{Example 8, Regular shock reflection. The convergence history of the residue
 as a function of number of iterations for two schemes with different CFL numbers.}
\label{ex8fig1}
\end{figure}

\begin{figure}
\centering
\subfigure[RK Jacobi scheme]{
\begin{minipage}[t]{0.5\linewidth}
\centering
\includegraphics[width=2.1in]{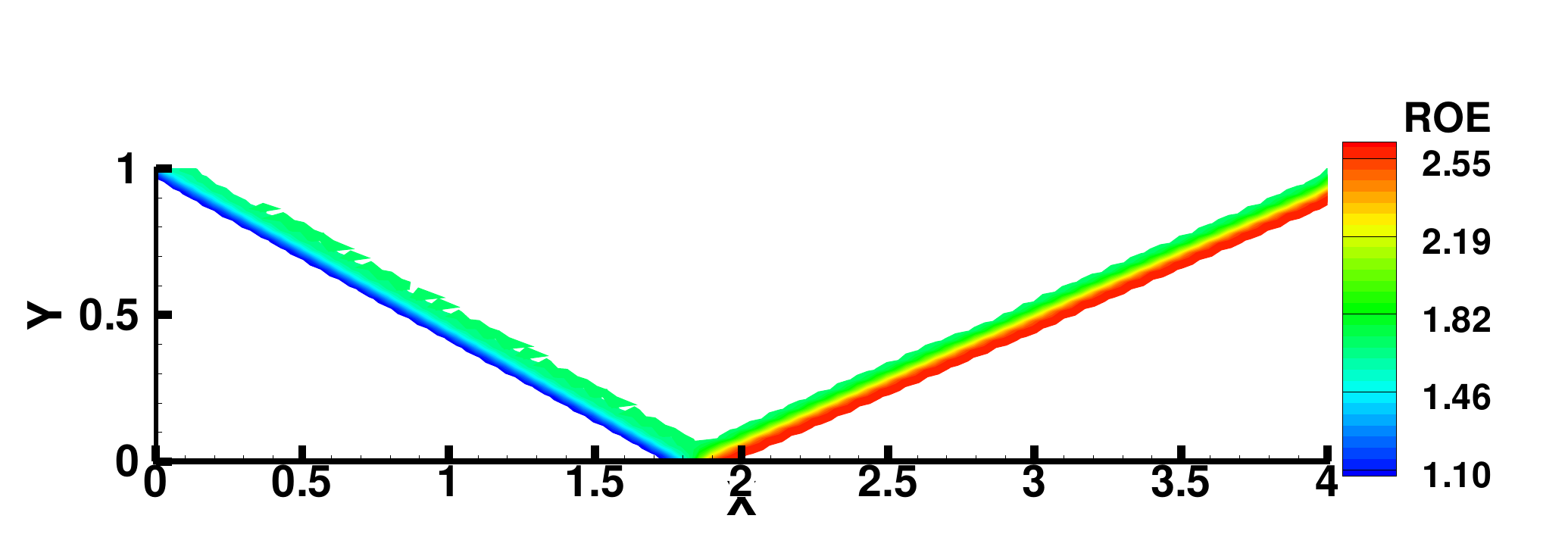}
\end{minipage}%
}%
\subfigure[FE fast sweeping scheme]{
\begin{minipage}[t]{0.5\linewidth}
\centering
\includegraphics[width=2.1in]{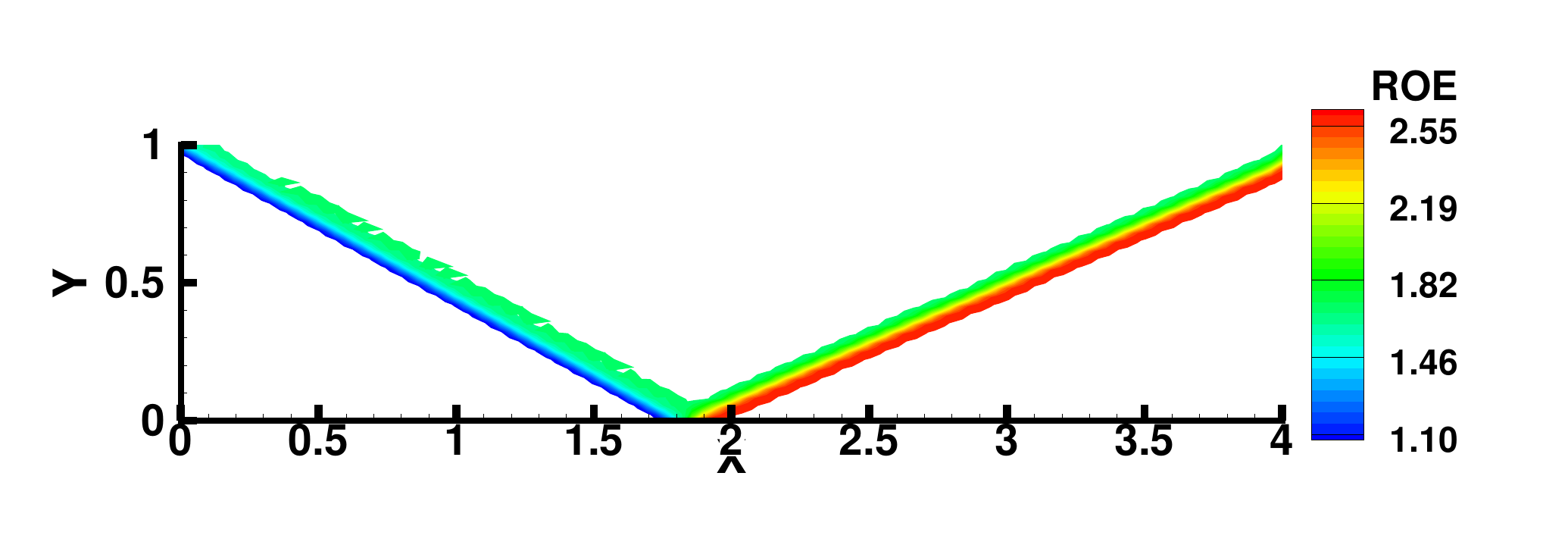}
\end{minipage}%
}%
\centering
\caption{Example 8, Regular shock reflection. 30 equally spaced density contours from 1.1 to 2.6
of the converged steady states of numerical solutions by two different iterative schemes.}
\label{ex8fig2}
\end{figure}

In Table \ref{tab:Margin_settings8}, number of iterations required
to reach the convergence threshold value $10^{-12}$, the final time and total CPU
time when convergence is obtained for these three iterative schemes with different
CFL numbers $\gamma$ are reported. Again, the FE Jacobi scheme
has to use a CFL number as small as $\gamma=0.1$ in order to converge, which leads to many
 iterations and quite large CPU time cost. We find that the residue of the FE Jacobi scheme hangs at $10^{-11.8}$ if $\gamma$ is increased to $0.2$, and oscillate between $10^{-9.7}$ and $10^{-9.8}$
if $\gamma=0.3$. By using the RK Jacobi (TVD RK3) scheme, the CFL number is enlarged to $\gamma=0.6$ and the iterations
converge much more efficiently than the FE Jacobi scheme. The FE fast sweeping scheme still turns out to be the most efficient one, which converges under the similar CFL numbers as the RK Jacobi scheme and only takes $98.92$ seconds of CPU time to reach the steady state by using the largest possible CFL number. In fact, comparing to $289.89$ seconds, the least CPU time cost case of the RK Jacobi scheme, the FE fast sweeping scheme saves about $66\%$ CPU time cost to reach the steady state solution in this example.
In Figure \ref{ex8fig1},
residue history in terms of iterations for the RK Jacobi and the FE fast sweeping schemes with different CFL numbers is
shown. We see that the residue of iterations settles down to values at the level of round off errors for all cases.
This result justifies the performance improvement of the absolutely convergent fixed-point fast sweeping method proposed here over our previous work in \cite{WuLiang}. In Figure \ref{ex8fig2}, contour plots of the density variable of the converged steady state solutions of the RK Jacobi and the FE fast sweeping schemes are shown, and the similar results are obtained as expected.

\begin{table}
\centering
\begin{tabular}{|c|c|c|c|}\hline
			\multicolumn{4}{|c|}{FE Jacobi scheme }\\\hline
            $\gamma:$ CFL number & iteration number & final time & CPU time \\\hline
0.1	&12046	&5.09	&618.36\\\hline
0.2  & Not convergent & - & - \\\hline
\end{tabular}
\begin{tabular}{|c|c|c|c|}\hline
			\multicolumn{4}{|c|}{RK Jacobi scheme }\\\hline
            $\gamma:$ CFL number & iteration number & final time & CPU time \\\hline
0.3	&11268	&4.76	&579.22\\\hline
0.4	&8454	&4.76	&436.33\\\hline
0.5	&6762	&4.76	&348.09\\\hline
0.6	&5634	&4.76	&289.89\\\hline
0.7   & Not convergent &   -    & - \\
\hline
\end{tabular}
\begin{tabular}{|c|c|c|c|}\hline
			\multicolumn{4}{|c|}{FE fast sweeping scheme }\\\hline
            $\gamma:$ CFL number & iteration number & final time & CPU time \\\hline
0.3	&3651	&4.62	&188.14\\\hline
0.4	&2722	&4.59	&140.03\\\hline
0.5	&2170	&4.57	&111.91\\\hline
0.6	&1934	&4.89	&98.92\\\hline
0.7   & Not convergent &    -   & - \\
    \hline
\end{tabular}
\caption{Example 8, Regular shock reflection. Number of iterations, the final time and total CPU time when convergence is obtained. Convergence criterion threshold value is $10^{-12}$. CPU time unit: second}
\label{tab:Margin_settings8}
\end{table}

\newpage
\noindent{\bf Example 9. Supersonic flow past a  plate with an attack angle}

\noindent In the following we test the schemes by solving a series of benchmark problems about
 supersonic flow past two-dimensional plates from \cite{JUNZ4, JUNZ3}.
The first example is a supersonic flow past a two-dimensional plate
with an attack angle of $\alpha=10 ^{\circ}$ . The free stream Mach number is $M_{\infty}$=3. The ideal gas
goes from the left toward the plate. The initial conditions are $p=\frac{1}{\gamma' M_{\infty}^{2}}, \rho=1, u=\cos(\alpha)$ and $v=\sin(\alpha)$. The computational domain is $[0,10]\times[-5,5]$. The plate is located at $x\in[1,2]$ and $y=0$.  The slip boundary condition is imposed on the plate. The physical values of the inflow and outflow boundary conditions are applied in different directions. The computational grid is $200 \times 200$.
The convergence criterion threshold value is set to be $10^{-12}$.

As that pointed out in \cite{JUNZ4, JUNZ3}, for this kind of problems, although the boundary is very far away from the plate, the waves including the shocks and the rarefaction waves propagate to the far field boundaries. This makes
it very challenging for the residue of high order schemes to converge to small values at the round off error level.
Again, with the help of the multi-resolution WENO scheme \cite{JUNZ}, the iteration residue of the absolutely convergent fixed-point fast sweeping method proposed here
settles down to the round off error level easily, and at the same time the high efficiency of the original
fixed-point fast sweeping WENO scheme in \cite{WuLiang} is preserved, as shown in the following numerical results.

\begin{figure}
\centering
\subfigure[RK Jacobi scheme]{
\begin{minipage}[t]{0.5\linewidth}
\centering
\includegraphics[width=2.1in]{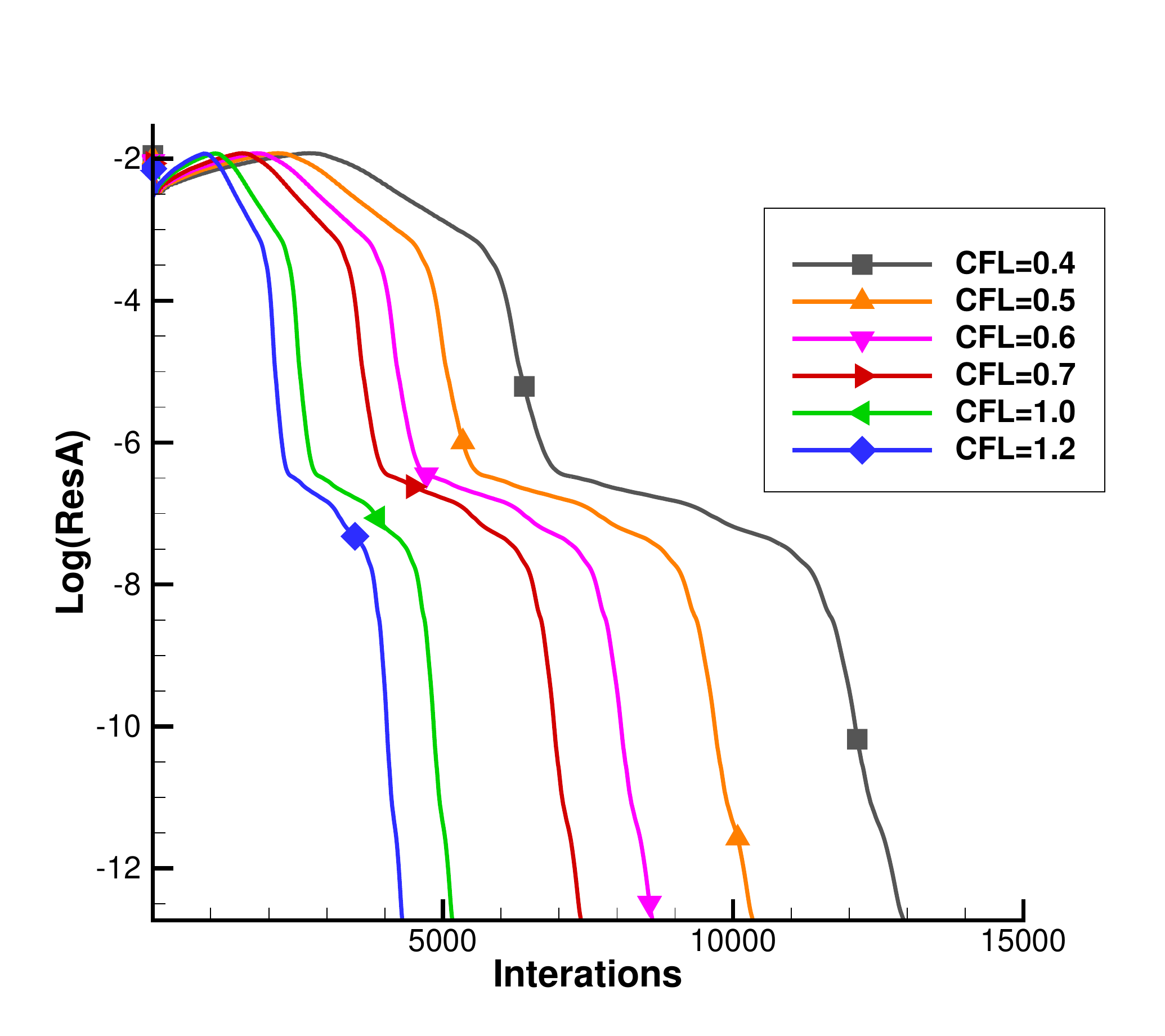}
\end{minipage}%
}%
\subfigure[FE fast sweeping scheme]{
\begin{minipage}[t]{0.5\linewidth}
\centering
\includegraphics[width=2.1in]{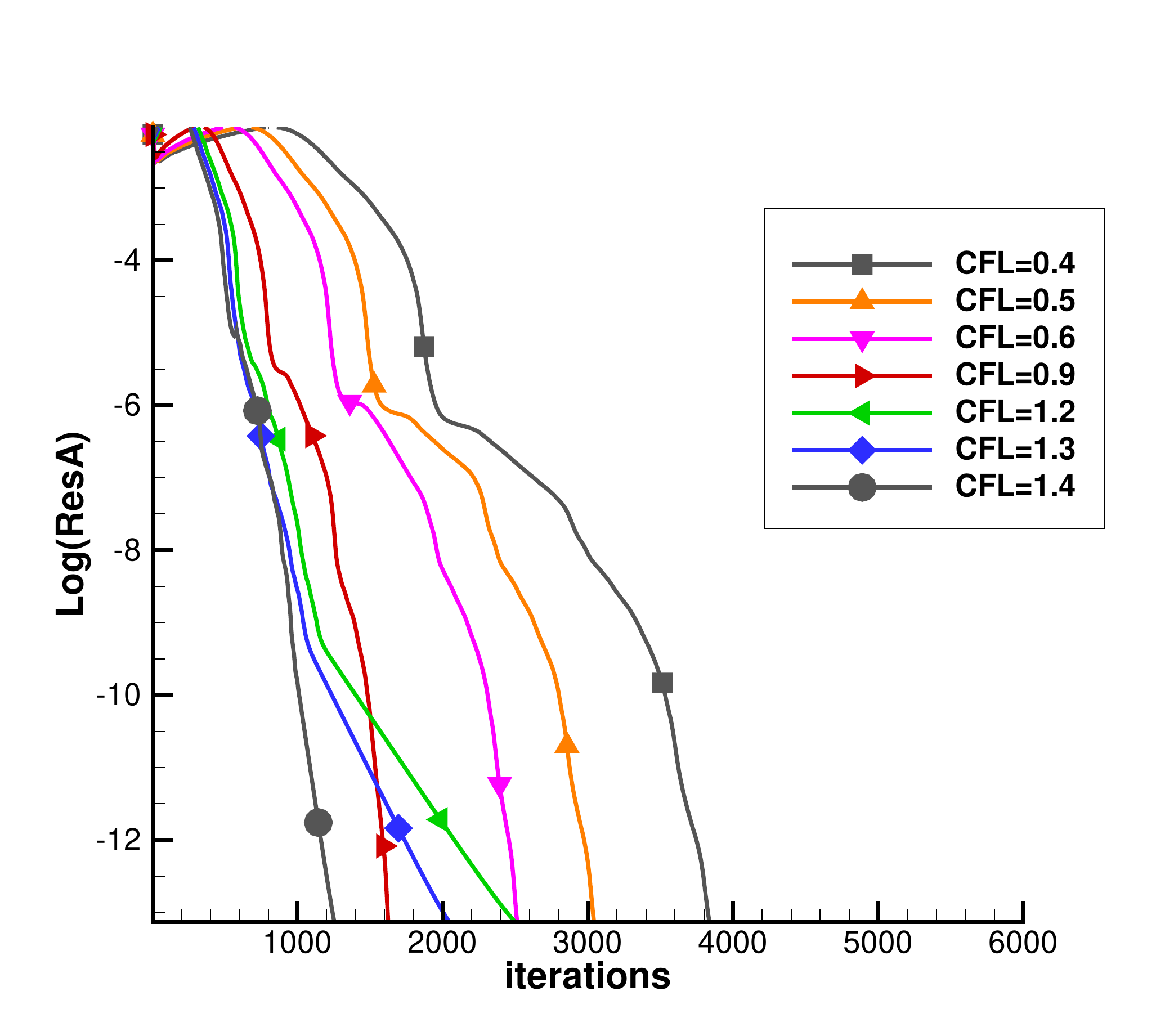}
\end{minipage}%
}%
\centering
\caption{Example 9, supersonic flow past a  plate with an attack angle. The convergence history of the residue
 as a function of number of iterations for two schemes with different CFL numbers.}
 \label{ex9fig1}
\end{figure}

\begin{figure}
\centering
\subfigure[RK Jacobi scheme]{
\begin{minipage}[t]{0.5\linewidth}
\centering
\includegraphics[width=2.1in]{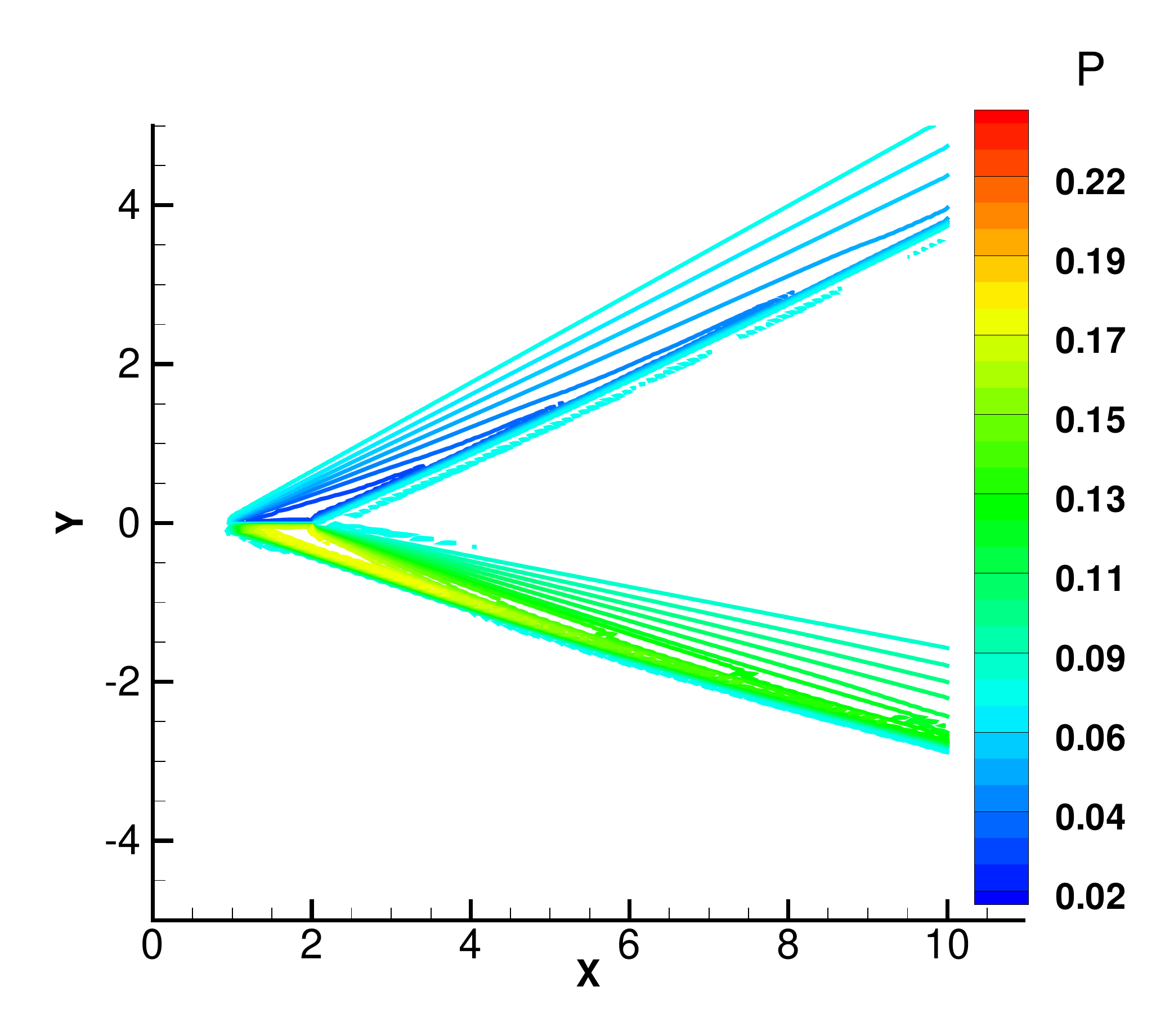}
\end{minipage}%
}%
\subfigure[FE fast sweeping scheme]{
\begin{minipage}[t]{0.5\linewidth}
\centering
\includegraphics[width=2.1in]{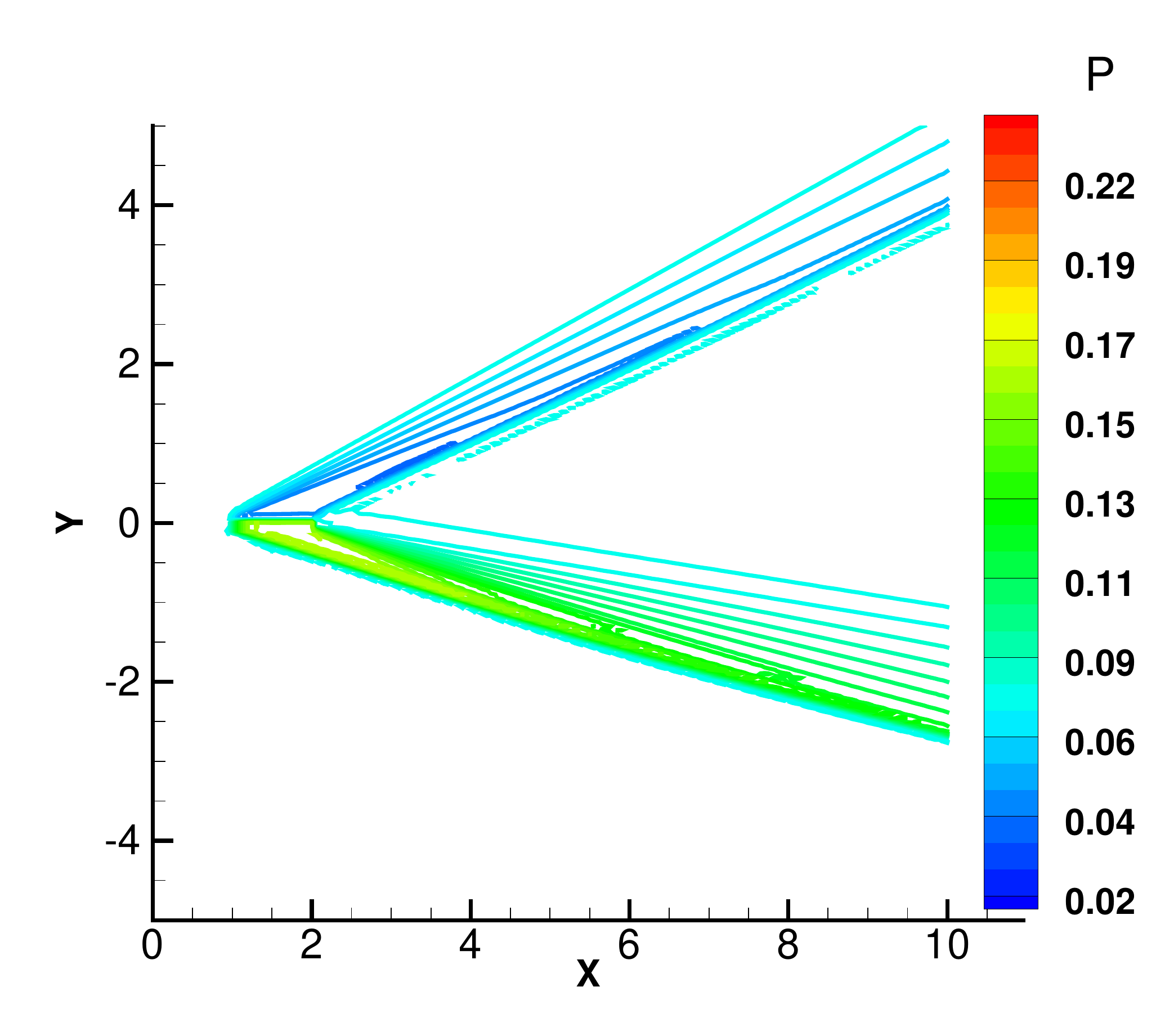}
\end{minipage}%
}%
\centering
\caption{Example 9, supersonic flow past a  plate with an attack angle. 30 equally spaced pressure contour from 0.02 to 0.23 of the converged steady states of numerical solutions by two different iterative schemes.}
\label{ex9fig2}
\end{figure}

In Table \ref{tab:Margin_settings9}, number of iterations required
to reach the convergence threshold value $10^{-12}$, the final time and total CPU
time when convergence is obtained for these three iterative schemes with different
CFL numbers $\gamma$ are reported. As previous examples, the FE Jacobi scheme
has to use a CFL number as small as $\gamma=0.1$ in order to converge, hence it takes a lot of CPU time cost.
The RK Jacobi (TVD RK3) scheme increases the CFL number to $\gamma=1.2$ and the iterations
converge much more efficiently than the FE Jacobi scheme. In this example, the FE fast sweeping scheme
can still converge with a larger CFL number ($\gamma=1.4$) than the TVD RK3 scheme, and only takes $631.50$ seconds of CPU time to reach the steady state by using the largest possible CFL number while the TVD RK3 scheme needs $1169.30$ seconds of CPU time. $46\%$ CPU time has been saved here. In Figure \ref{ex9fig1},
residue history in terms of iterations for the RK Jacobi and the FE fast sweeping schemes with different CFL numbers is
shown. It can be seen that the residue of iterations settles down to values at the level of round off errors for all cases.
Again, this result verifies the performance improvement of the absolutely convergent fixed-point fast sweeping method proposed here over our previous work in \cite{WuLiang}. In Figure \ref{ex9fig2}, contour plots of the pressure variable of the converged steady state solutions of the RK Jacobi and the FE fast sweeping schemes are shown, and the comparable results are observed.

\begin{table}
		\centering
		\begin{tabular}{|c|c|c|c|}\hline
			\multicolumn{4}{|c|}{FE Jacobi scheme }\\\hline
            $\gamma:$ CFL number & iteration number & final time & CPU time \\\hline
0.1   & 17337 & 30.53 & 4626.95 \\
\hline
0.2  & Not convergent & - & - \\\hline
		\end{tabular}
\begin{tabular}{|c|c|c|c|}\hline
			\multicolumn{4}{|c|}{RK Jacobi scheme }\\\hline
            $\gamma:$ CFL number & iteration number & final time & CPU time \\\hline
0.4   & 13179 & 30.94 & 3527.08 \\
\hline
0.5   & 10488 & 30.78 & 2811.88 \\
\hline
0.7   & 7470  & 30.69 & 1994.08 \\
\hline
1.0     & 5220  & 30.64 & 1406.23 \\
\hline
1.2   & 4347  & 30.62 & 1169.30 \\
\hline
1.3   & Not convergent &   -    & - \\
    \hline
		\end{tabular}
\begin{tabular}{|c|c|c|c|}\hline
			\multicolumn{4}{|c|}{FE fast sweeping scheme }\\\hline
            $\gamma:$ CFL number & iteration number & final time & CPU time \\\hline
          0.4   & 3748  & 27.79 & 2043.47 \\
\hline
0.5   & 2976  & 27.58 & 1622.19 \\
\hline
0.7   & 2384  & 30.93 & 1294.97 \\
\hline
0.9   & 1588  & 26.48 & 861.66 \\
\hline
1.4   & 1164  & 29.87 & 631.50 \\
\hline
1.5   & Not convergent &    -   & - \\    \hline
		\end{tabular}
		\caption{Example 9, supersonic flow past a  plate with an attack angle. Number of iterations, the final time and total CPU time when convergence is obtained. Convergence criterion threshold value is $10^{-12}$. CPU time unit: second}
		\label{tab:Margin_settings9}
	\end{table}

\bigskip
\noindent{\bf Example 10. Supersonic flow past two plates with an attack angle}

\noindent In this example we replace one plate in Example 9 by two plates, i.e., a supersonic flow past two plates with an attack angle of
$\alpha=10 ^{\circ}$. The free stream Mach number is $M_{\infty}=3$. The ideal gas goes from the left
toward these two plates. The initial conditions are $p=\frac{1}{\gamma' M_{\infty}^{2}}$, $\rho=1$, $u=\cos(\alpha)$ and $v=\sin(\alpha)$. The computational domain is $[0,10]\times [-5,5]$. Two plates are placed at $x \in[2,3]$ with $y=-2$ and at
$x \in [2,3]$ with $y=2$. The slip boundary conditions are imposed on these plates. The physical
values of the inflow and outflow boundary conditions are applied on all
 boundaries of the domain. The computational grid is $200 \times 200$.
The convergence criterion threshold value is set to be $10^{-12}$.

The solution of this example contains strong shocks, rarefaction waves and their interactions. It is more complex than   Example 9. Again, from the following numerical results, we can see that the iteration residue of the absolutely convergent fixed-point fast sweeping method proposed in this paper
settles down to the round off error level easily, and it is the most efficient one among these three schemes discussed here.

\begin{figure}
\centering
\subfigure[RK Jacobi scheme]{
\begin{minipage}[t]{0.5\linewidth}
\centering
\includegraphics[width=2.1in]{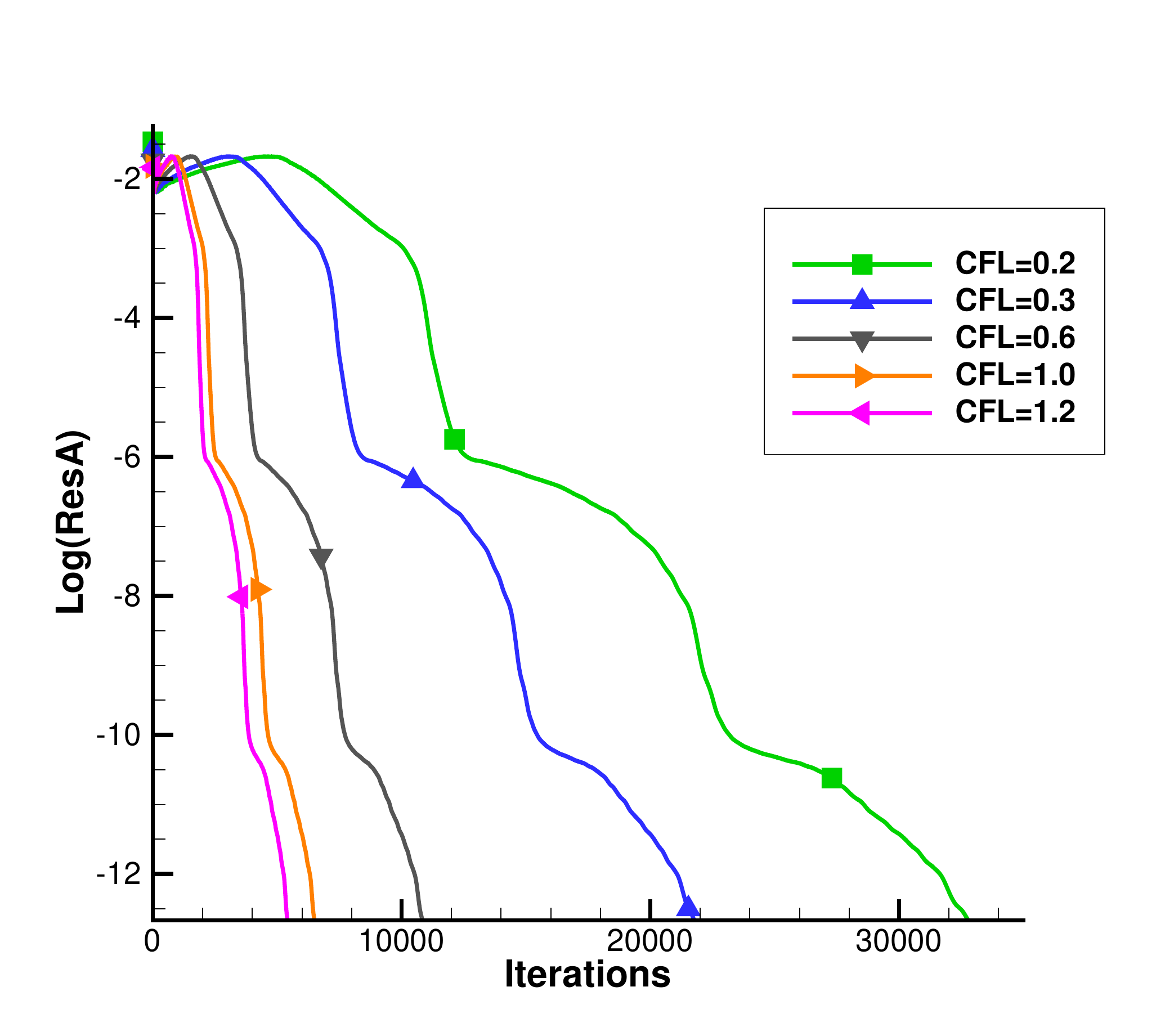}
\end{minipage}%
}%
\subfigure[FE fast sweeping scheme]{
\begin{minipage}[t]{0.5\linewidth}
\centering
\includegraphics[width=2.1in]{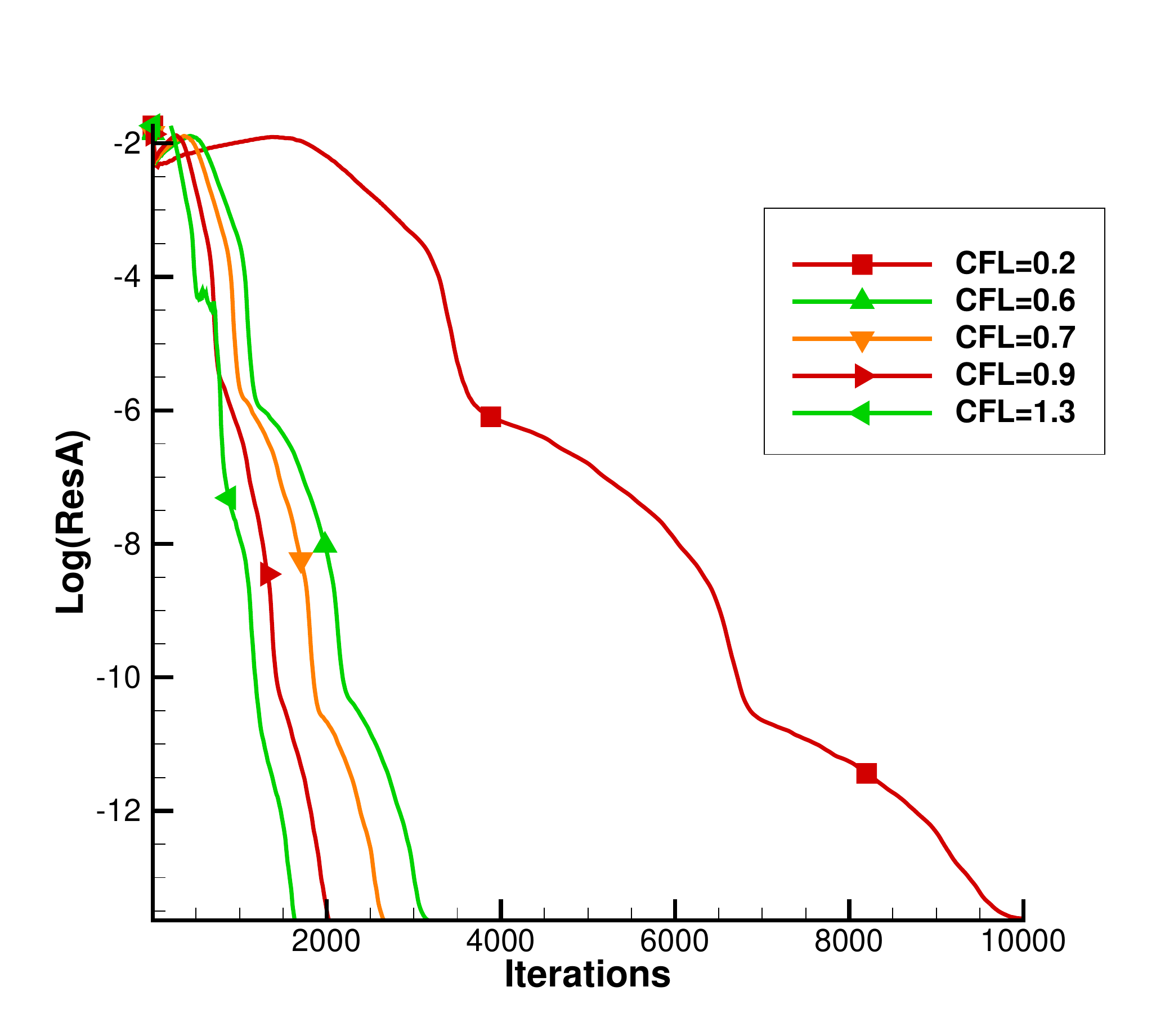}
\end{minipage}%
}%
\centering
\caption{Example 10, supersonic flow past two plates with an attack angle. The convergence history of the residue
 as a function of number of iterations for two schemes with different CFL numbers.}
 \label{ex10fig1}
\end{figure}

\begin{figure}
\centering
\subfigure[RK Jacobi scheme]{
\begin{minipage}[t]{0.5\linewidth}
\centering
\includegraphics[width=2.1in]{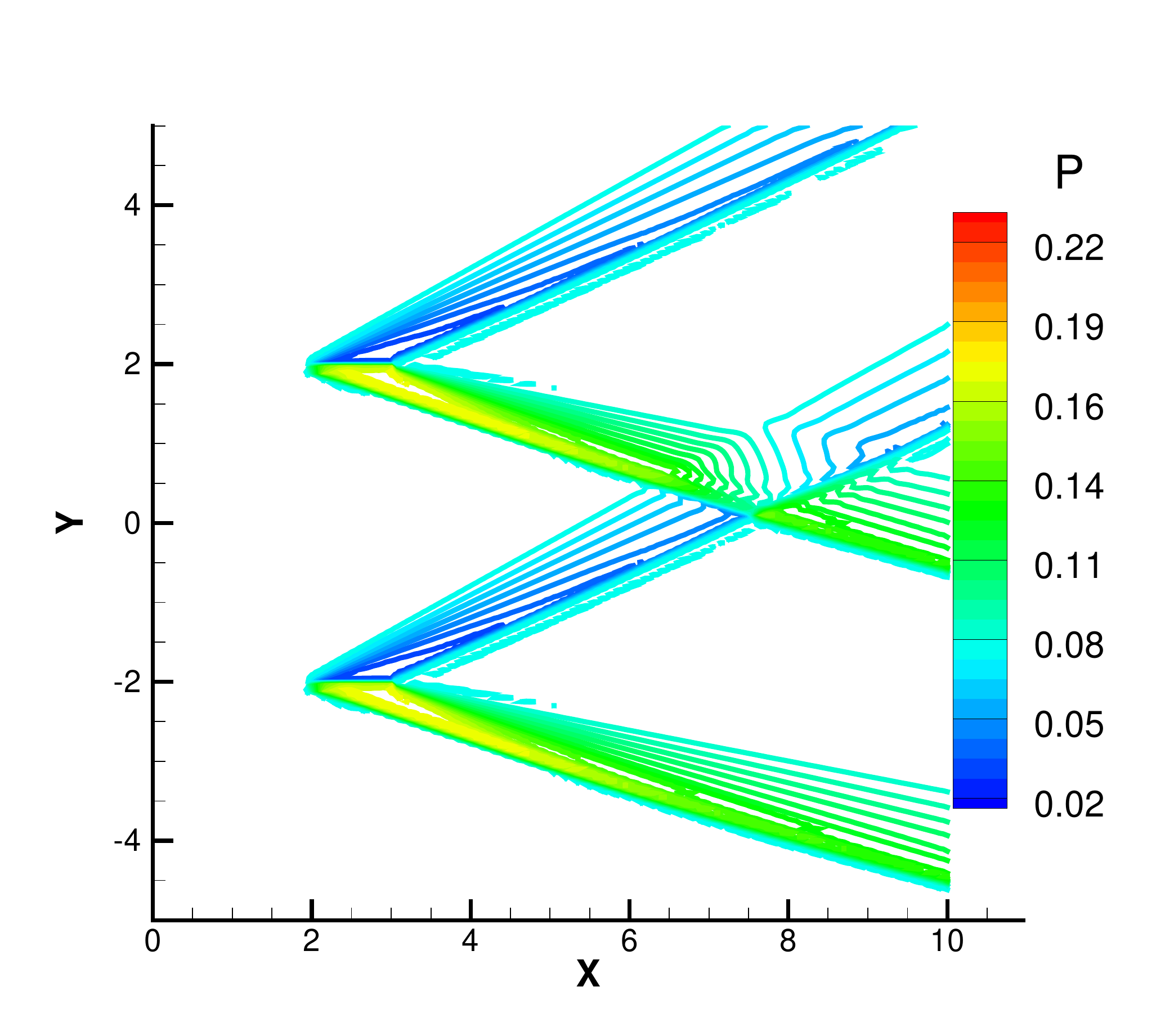}
\end{minipage}%
}%
\subfigure[FE fast sweeping scheme]{
\begin{minipage}[t]{0.5\linewidth}
\centering
\includegraphics[width=2.1in]{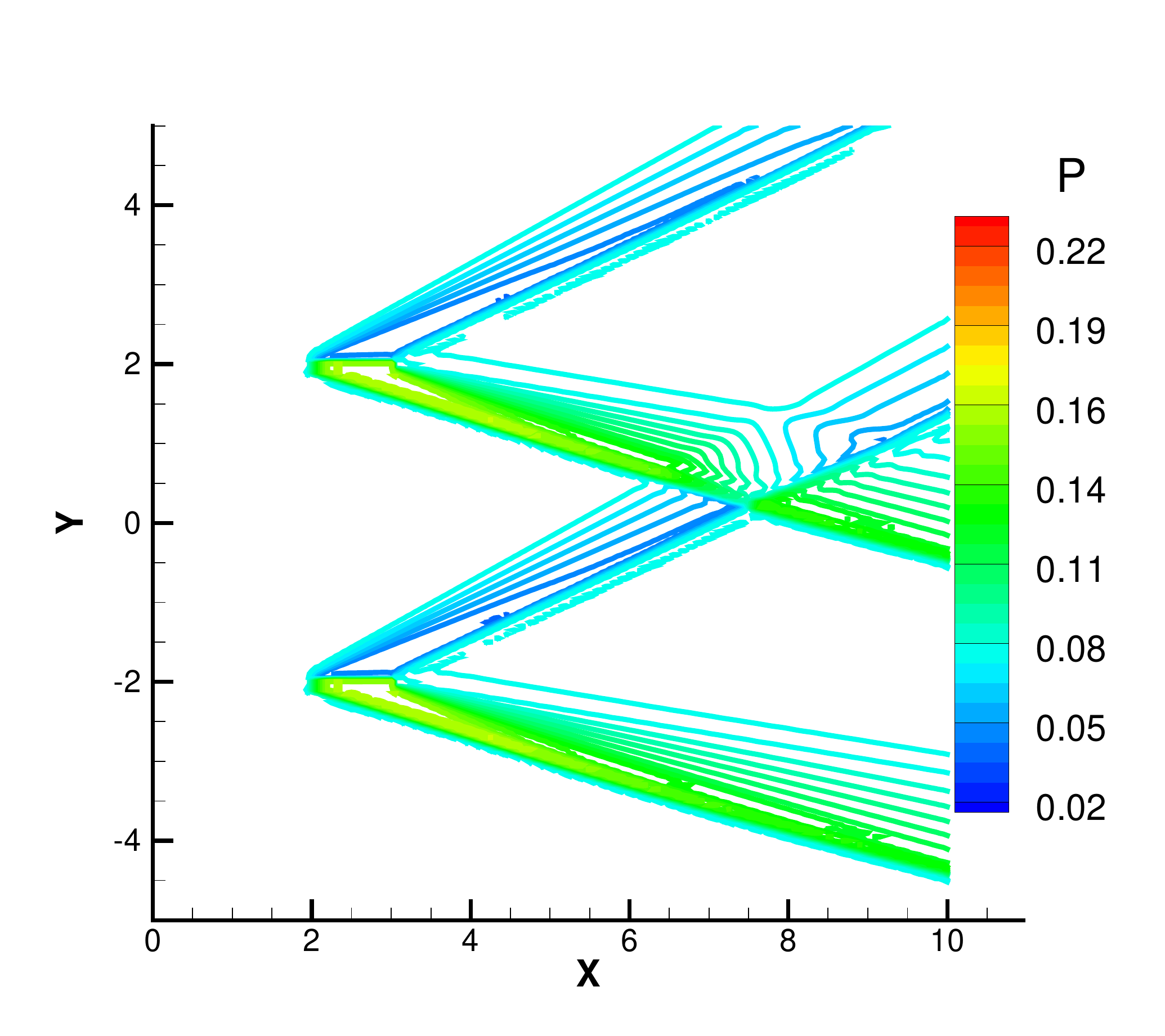}
\end{minipage}%
}%
\centering
\caption{Example 10, supersonic flow past two plates with an attack angle. 30 equally spaced  pressure contour from 0.02 to 0.23 of the converged steady states of numerical solutions by two different iterative schemes.}
\label{ex10fig2}
\end{figure}

In Table \ref{tab:Margin_settings10}, number of iterations required
to reach the convergence threshold value $10^{-12}$, the final time and total CPU
time when convergence is obtained for these three iterative schemes with different
CFL numbers $\gamma$ are reported. It is observed that the largest CFL number the FE Jacobi scheme
can achieve for convergence is $\gamma=0.1$ which makes it the most expensive one among these three schemes.
The RK Jacobi (TVD RK3) scheme can use much larger CFL numbers till $\gamma=1.2$, by which it only takes
$1431.80$ seconds of CPU time to reach the convergence threshold value. The FE fast sweeping scheme can use the
largest CFL number $\gamma=1.3$ to converge to the steady state solution in this example, which makes it the most efficient one. Only $810.30$ seconds of CPU time are needed by the FE fast sweeping scheme to satisfy the convergence criterion, which saves more than $43\%$ CPU time of that by the TVD RK3 scheme.
In Figure \ref{ex10fig1},
residue history in terms of iterations for the RK Jacobi and the FE fast sweeping schemes with different CFL numbers is
shown. We see that the residue of iterations settles down to values at the level of round off errors for these cases.
In Figure \ref{ex10fig2}, contour plots of the pressure variable of the converged steady state solutions of the RK Jacobi and the FE fast sweeping schemes are presented, which show the comparable results.

\begin{table}
		\centering
		\begin{tabular}{|c|c|c|c|}\hline
			\multicolumn{4}{|c|}{FE Jacobi scheme }\\\hline
            $\gamma:$ CFL number & iteration number & final time & CPU time \\\hline
0.1   & 21316 & 37.53 & 5764.19 \\
\hline
0.2  & Not convergent & - & - \\\hline
		\end{tabular}
\begin{tabular}{|c|c|c|c|}\hline
			\multicolumn{4}{|c|}{RK Jacobi scheme }\\\hline
            $\gamma:$ CFL number & iteration number & final time & CPU time \\\hline
0.3   & 21072 & 37.10  & 5701.58 \\
\hline
0.6   & 10524 & 37.06 & 2863.86\\
\hline
1.0     & 6315  & 37.06 & 1716.80 \\
\hline
1.2   & 5262  & 37.06 & 1431.80 \\
\hline
1.3   & Not convergent &   -    & - \\
    \hline
\end{tabular}
\begin{tabular}{|c|c|c|c|}\hline
			\multicolumn{4}{|c|}{FE fast sweeping scheme}\\\hline
            $\gamma:$ CFL number & iteration number & final time & CPU time \\\hline
           0.6   & 2836  & 31.43 & 1535.53 \\
\hline
0.7   & 2388  & 30.88 & 1303.16 \\
\hline
0.8   & 2056  & 30.38 & 1125.40 \\
\hline
0.9   & 1820  & 30.25 & 984.02 \\
\hline
1.3   & 1476  & 35.23 & 810.30 \\
\hline
1.4   & Not convergent &   -    & - \\    \hline
\end{tabular}
\caption{Example 10, supersonic flow past two plates with an attack angle. Number of iterations, the final time and total CPU time when convergence is obtained. Convergence criterion threshold value  is $10^{-12}$. CPU time unit: second}
		\label{tab:Margin_settings10}
	\end{table}

\bigskip
\noindent{\bf Example 11. Supersonic flow past three plates with an attack angle}

 \noindent In this example we place three plates in the field such that there are more interactions of different waves
 and the flow is more complex than previous examples.
  The supersonic flow passes three plates with an attack angle of $\alpha=10 ^{\circ}$. The free stream Mach number is $M_{\infty}=3$. The ideal gas goes from the left side of the domain toward the plates. The initial conditions are $p=\frac{1}{\gamma' M_{\infty}^{2}}$, $\rho=1$, $u=\cos(\alpha)$ and $v=\sin(\alpha)$. The computational domain is  $[0,10]\times [-5,5]$. Three plates are placed at $x \in [2,3]$
with $y=-2$, $x \in [1,2]$ with $y=0$, and $x \in [2,3]$ with $y=2$. We impose the slip boundary condition
 on these plates, and apply the physical values of the inflow and outflow boundary conditions
on all
 boundaries of the computational domain. The computational grid is $200 \times 200$.
The convergence criterion threshold value is set to be $10^{-12}$.

\begin{figure}
\centering
\subfigure[RK Jacobi scheme]{
\begin{minipage}[t]{0.5\linewidth}
\centering
\includegraphics[width=2.1in]{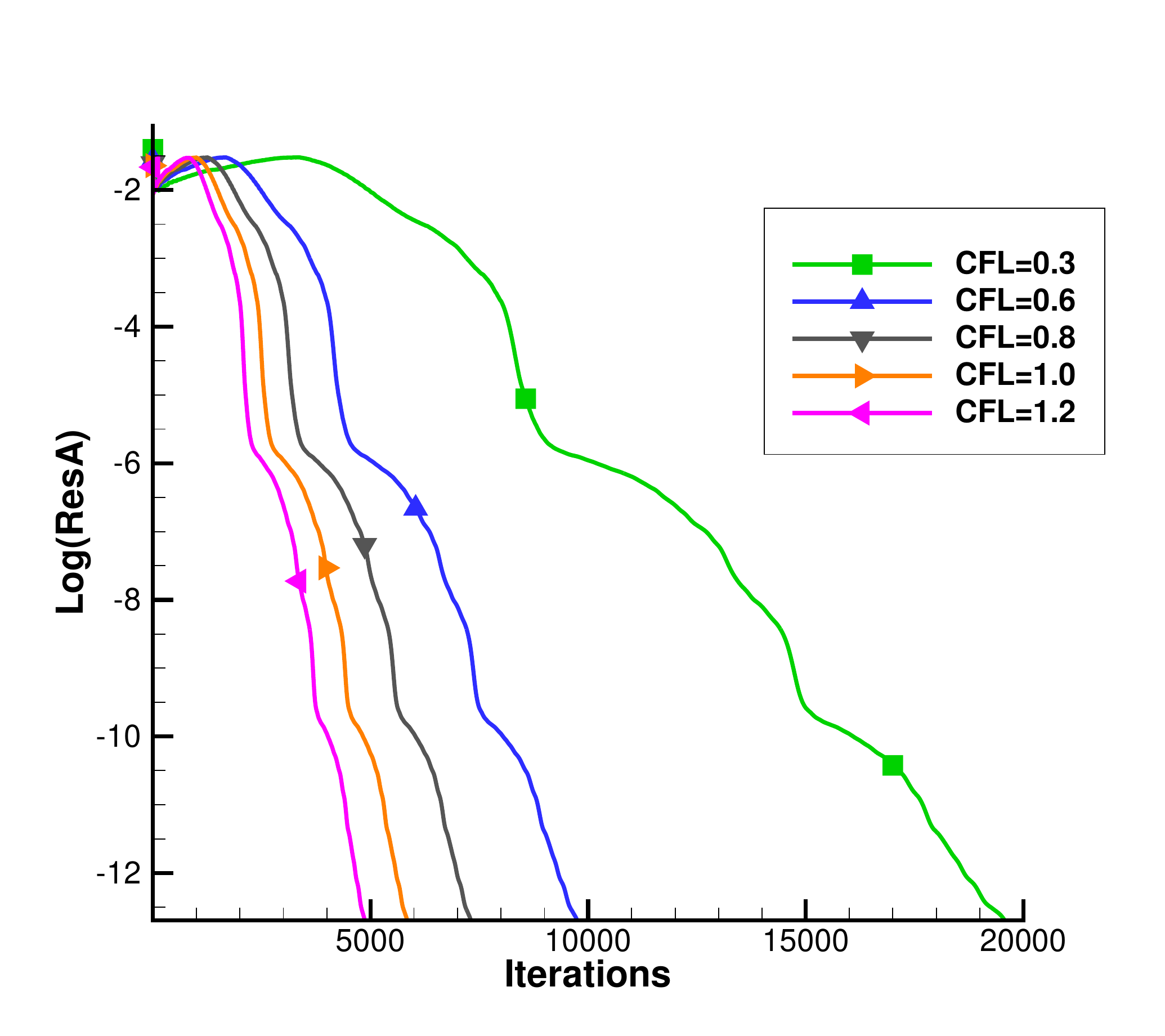}
\end{minipage}%
}%
\subfigure[FE fast sweeping scheme]{
\begin{minipage}[t]{0.5\linewidth}
\centering
\includegraphics[width=2.1in]{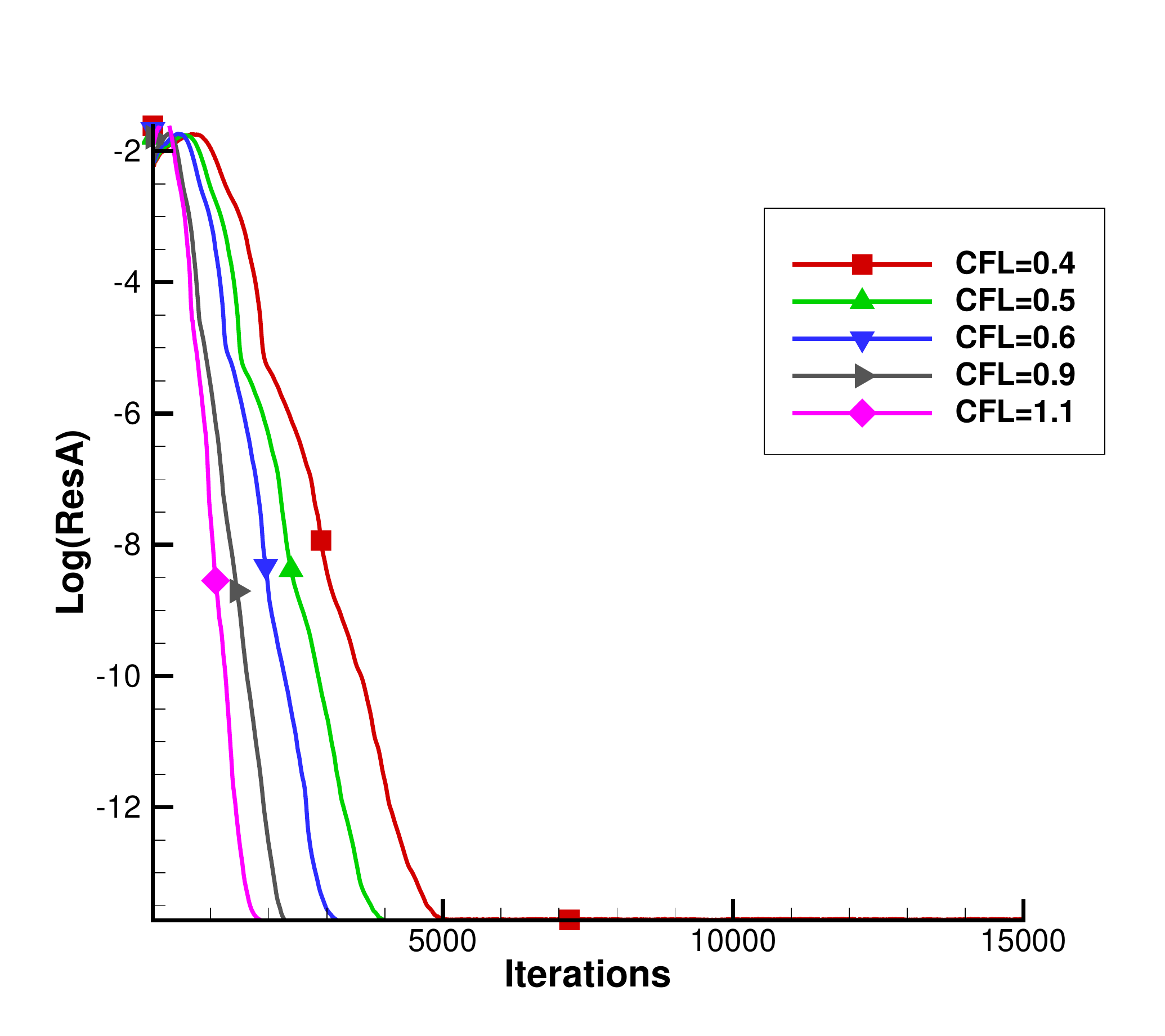}
\end{minipage}%
}%
\centering
\caption{Example 11, supersonic flow past three plates with an attack angle. The convergence history of the residue
 as a function of number of iterations for two schemes with different CFL numbers.}
 \label{ex11fig1}
\end{figure}

\begin{figure}
\centering
\subfigure[RK Jacobi scheme]{
\begin{minipage}[t]{0.5\linewidth}
\centering
\includegraphics[width=2.1in]{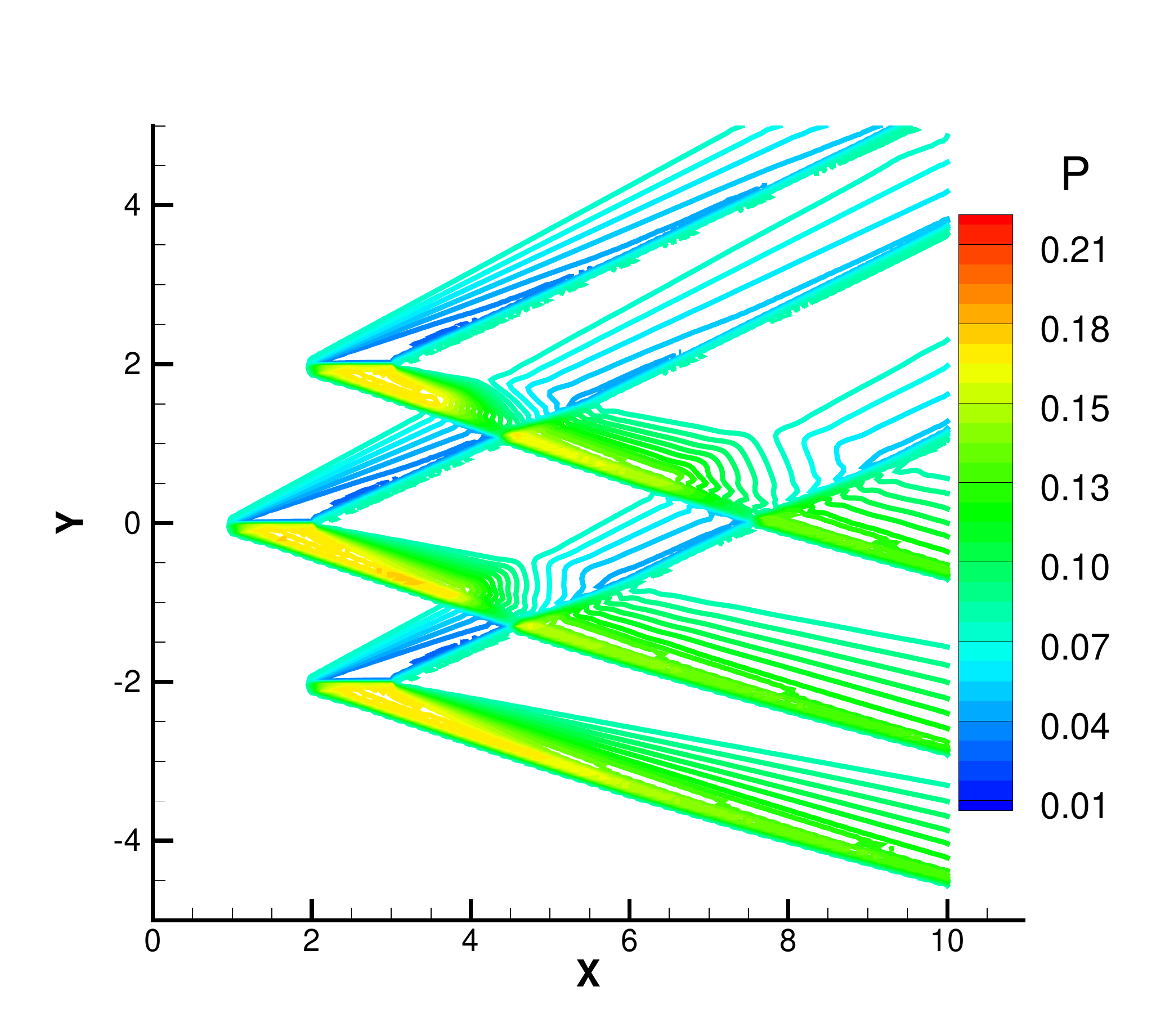}
\end{minipage}%
}%
\subfigure[FE fast sweeping scheme]{
\begin{minipage}[t]{0.5\linewidth}
\centering
\includegraphics[width=2.1in]{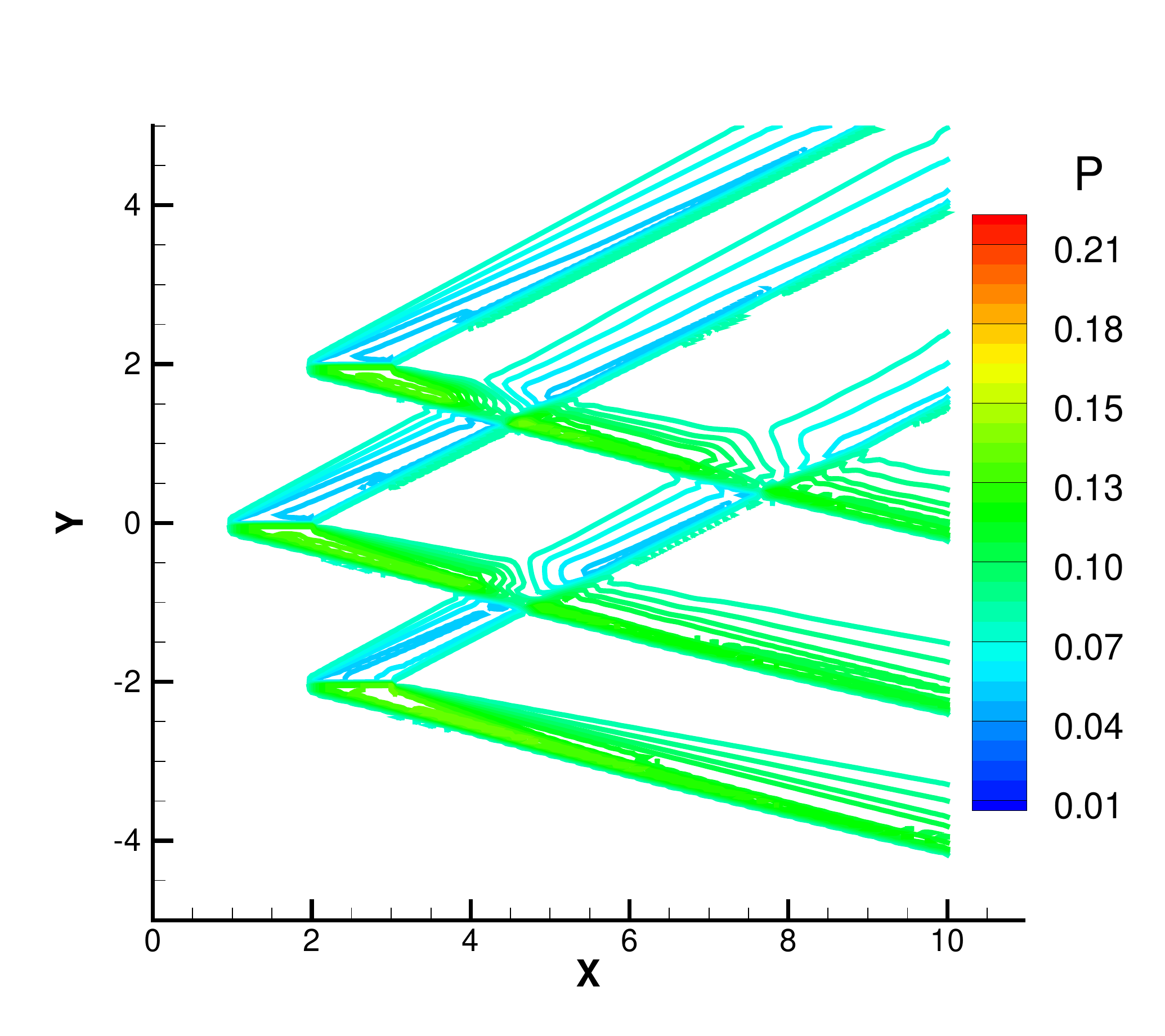}
\end{minipage}%
}%
\centering
\caption{Example 11, supersonic flow past three plates with an attack angle. 30 equally spaced  pressure contour from 0.01 to 0.22 of the converged steady states of numerical solutions by two different iterative schemes.}
\label{ex11fig2}
\end{figure}

In Table \ref{tab:Margin_settings11}, number of iterations required
to reach the convergence threshold value $10^{-12}$, the final time and total CPU
time when convergence is obtained for these three iterative schemes with different
CFL numbers $\gamma$ are reported. Residue history in terms of iterations for the RK Jacobi and the FE fast sweeping schemes with different CFL numbers is
shown in Figure \ref{ex11fig1}, and contour plots of the pressure variable of the converged steady state solutions of the RK Jacobi and the FE fast sweeping schemes are presented in Figure \ref{ex11fig2}. Based on these results for this more complex example, we draw similar conclusions as in previous examples. The proposed absolutely convergent fixed-point fast sweeping method is the most efficient scheme among all three schemes studied here. By using the largest CFL number permitted by each method to converge to steady states, the absolutely convergent fixed-point fast sweeping method saves about $40\%$ CPU time of that by the TVD RK3 (RK Jacobi) scheme, while numerical steady states
obtained by different schemes are comparable.

\begin{table}
		\centering
		\begin{tabular}{|c|c|c|c|}\hline
			\multicolumn{4}{|c|}{FE Jacobi scheme}\\\hline
            $\gamma:$ CFL number & iteration number & final time & CPU time \\\hline

0.1   & 19235 & 33.88 & 5230.52 \\
\hline
0.2  & Not convergent & - & - \\\hline
		\end{tabular}
\begin{tabular}{|c|c|c|c|}\hline
			\multicolumn{4}{|c|}{RK Jacobi scheme }\\\hline
            $\gamma:$ CFL number & iteration number & final time & CPU time \\\hline
0.3   & 18636 & 32.83 & 5056.52 \\
\hline
0.6   & 9315  & 32.82 & 2536.50 \\
\hline
1.0     & 5589  & 32.82 & 1520.80 \\
\hline
1.2   & 4656  & 32.81 & 1267.42 \\
\hline
1.3   & Not convergent &   -    & - \\
    \hline
\end{tabular}
\begin{tabular}{|c|c|c|c|}\hline
\multicolumn{4}{|c|}{FE fast sweeping scheme }\\\hline
            $\gamma:$ CFL number & iteration number & final time & CPU time \\\hline
           0.5   & 3292  & 30.42 & 1802.92 \\
\hline
0.6   & 2652  & 29.41 & 1448.73 \\
\hline
0.8   & 1980  & 29.27 & 1079.00 \\
\hline
0.9   & 1916  & 31.87 & 1050.42 \\
\hline
1.1   & 1432  & 28.90  & 781.05 \\
\hline
1.2   & Not convergent &    -   & - \\    \hline
		\end{tabular}
\caption{Example 11, supersonic flow past three plates with an attack angle. Number of iterations, the final time and total CPU time when convergence is obtained. Convergence criterion threshold value is $10^{-12}$. CPU time unit: second}
		\label{tab:Margin_settings11}
	\end{table}

\bigskip
\noindent{\bf Example 12. Supersonic flow past a long plate with an attack angle}

\noindent In this example, we test the schemes by solving the case of a long plate in the flow field, i.e., a supersonic flow past a long plate with an attack angle of $\alpha=10 ^{\circ}$.
The
free stream has the Mach number $M_{\infty}=3$. The ideal gas goes from the left toward the long
plate. The initial condition to start the iterations is $p=\frac{1}{\gamma' M_{\infty}^{2}}$, $\rho=1$, $u=\cos(\alpha)$ and $v=\sin(\alpha)$.
 The computational domain is $[0,7]\times[-5,5]$, and the long plate is set in the region $x \in [2,7]$ with $y=0$.
The slip boundary condition is imposed on the long plate. The physical values of the inflow
and outflow boundary conditions are applied at the boundaries of the computational domain.
The difference from previous examples is that the plate extends to the right boundary, so the shock waves
and the rarefaction waves pass through the right boundary on both sides of the plate.
The computational grid is $140 \times 200$.
The convergence criterion threshold value is set to be $10^{-13}$.

\begin{figure}
\centering
\subfigure[RK Jacobi scheme]{
\begin{minipage}[t]{0.5\linewidth}
\centering
\includegraphics[width=2.1in]{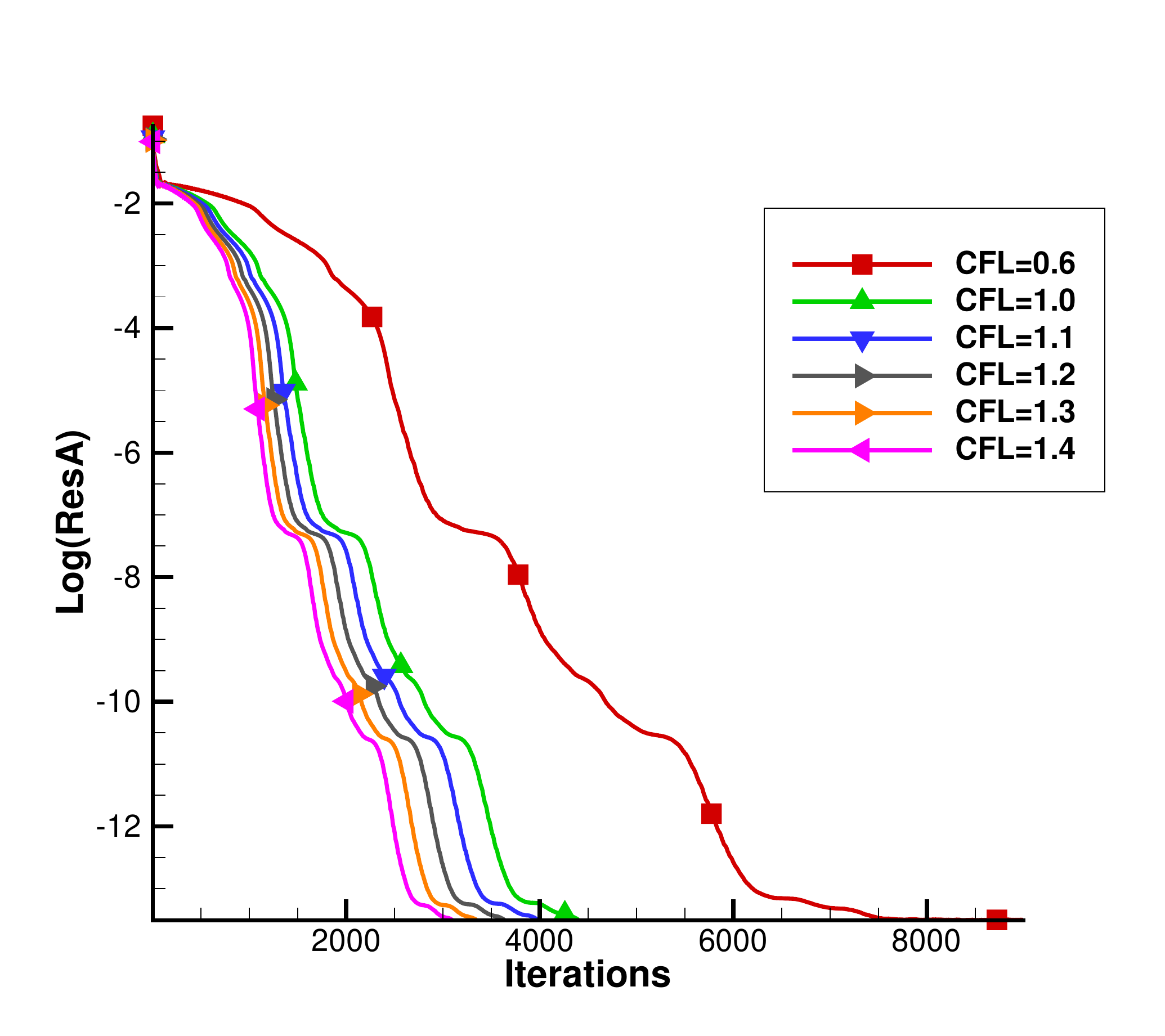}
\end{minipage}%
}%
\subfigure[FE fast sweeping scheme]{
\begin{minipage}[t]{0.5\linewidth}
\centering
\includegraphics[width=2.1in]{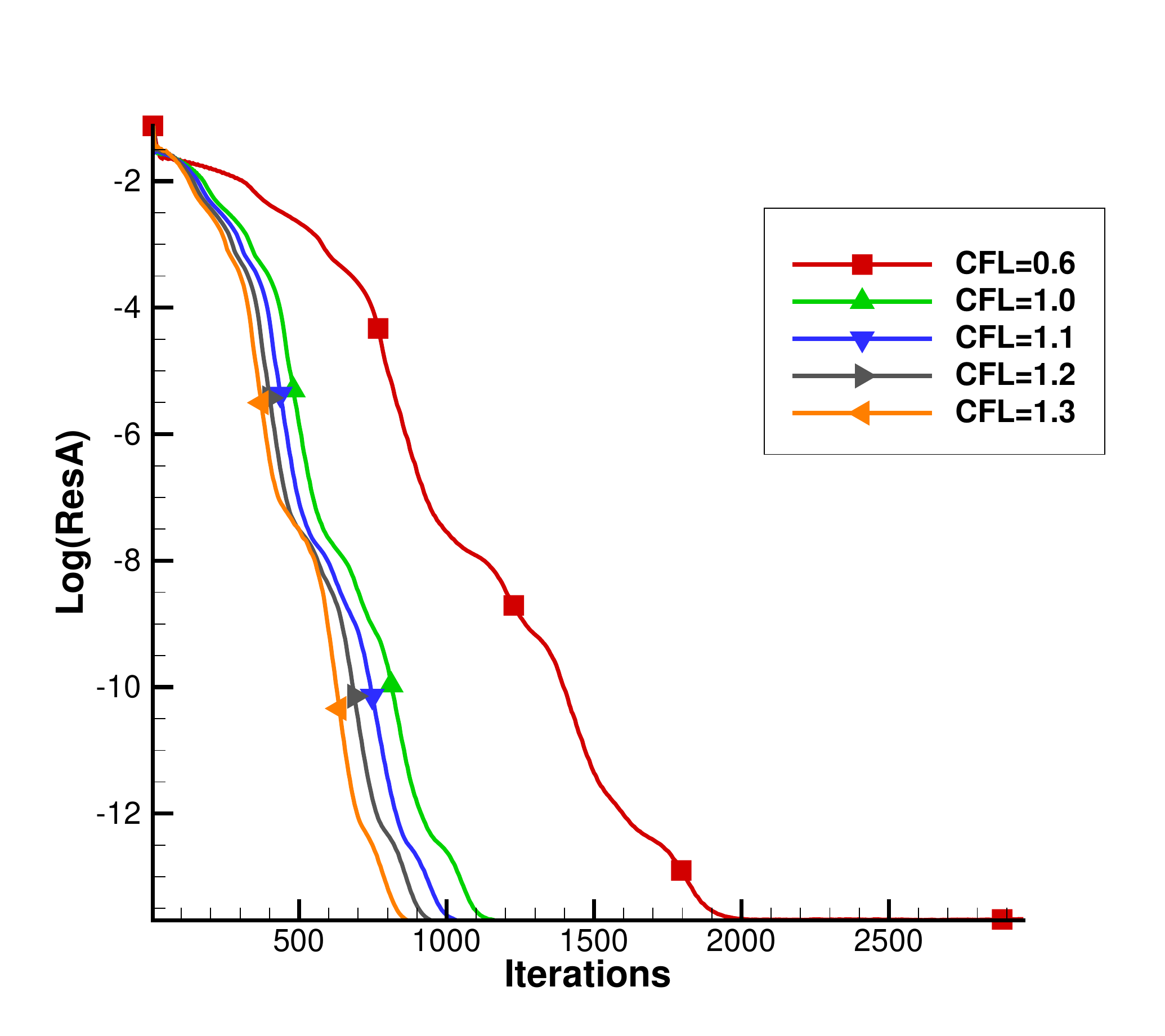}
\end{minipage}%
}%
\centering
\caption{Example 12, supersonic flow past a long plate with an attack angle. The convergence history of the residue
 as a function of number of iterations for two schemes with different CFL numbers.}
 \label{ex12fig1}
\end{figure}

\begin{figure}
\centering
\subfigure[RK Jacobi scheme]{
\begin{minipage}[t]{0.5\linewidth}
\centering
\includegraphics[width=2.1in]{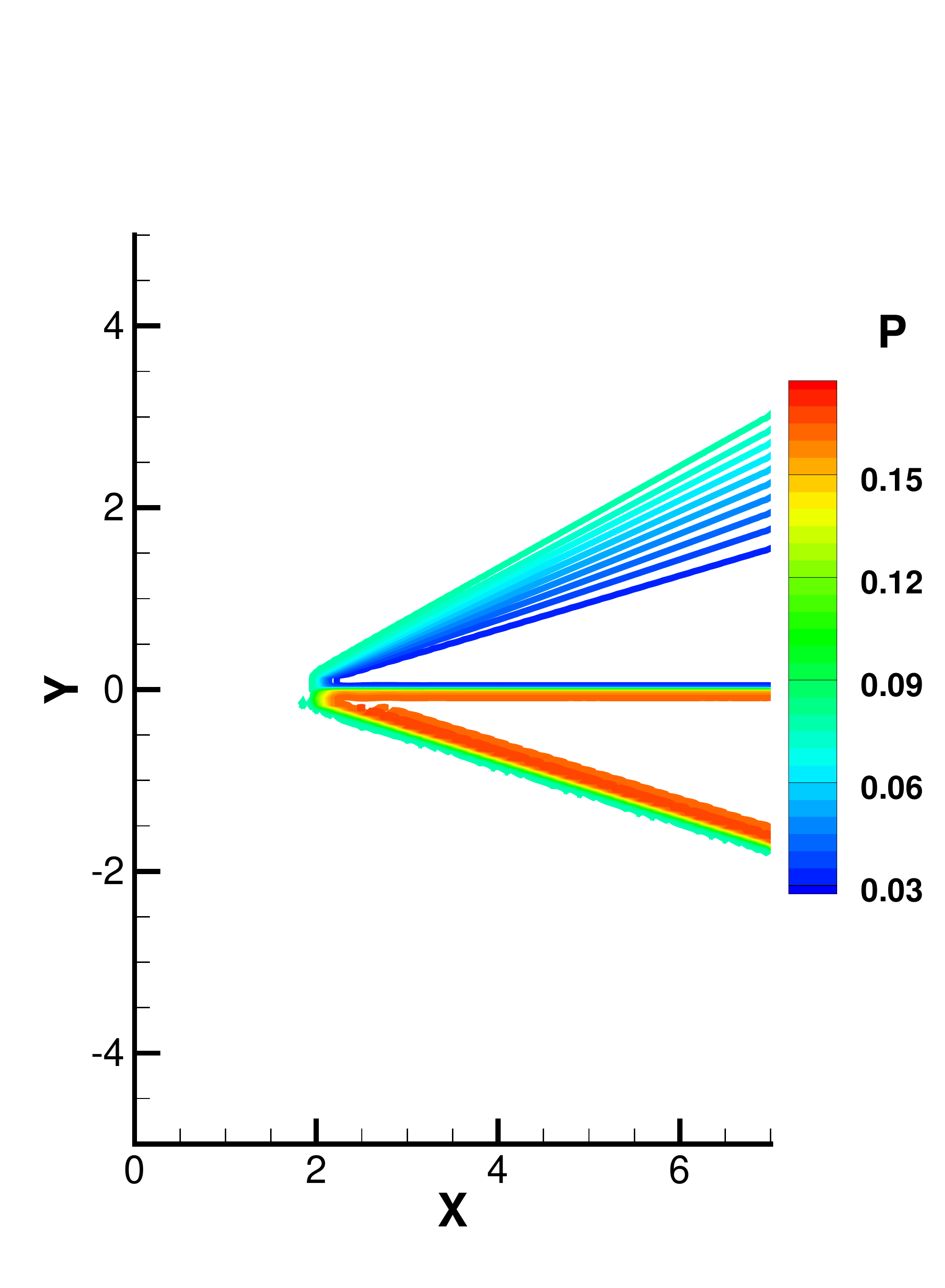}
\end{minipage}%
}%
\subfigure[FE fast sweeping scheme]{
\begin{minipage}[t]{0.5\linewidth}
\centering
\includegraphics[width=2.1in]{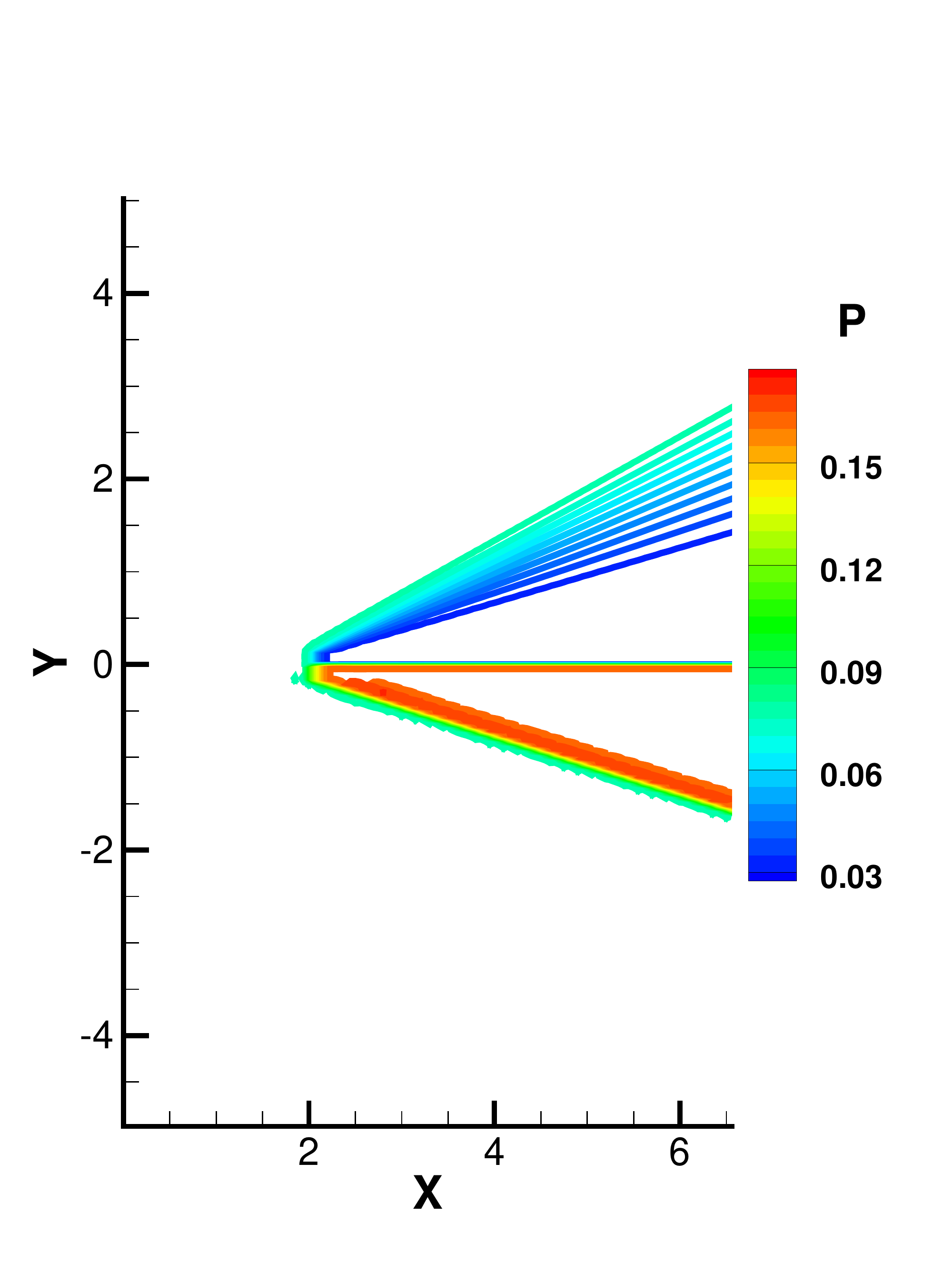}
\end{minipage}%
}%
\centering
\caption{Example 12, supersonic flow past a long plate with an attack angle. 30 equally pressure contour from 0.03 to 0.17 of the converged steady states of numerical solutions by two different iterative schemes.}
\label{ex12fig2}
\end{figure}

In Table \ref{tab:Margin_settings12}, number of iterations required
to reach the convergence threshold value $10^{-13}$, the final time and total CPU
time when convergence is obtained for these three iterative schemes with different
CFL numbers $\gamma$ are reported. Residue history in terms of iterations for the RK Jacobi and the FE fast sweeping schemes with different CFL numbers is
shown in Figure \ref{ex12fig1}, and contour plots of the pressure variable of the converged steady state solutions of the RK Jacobi and the FE fast sweeping schemes are presented in Figure \ref{ex12fig2}. Again, we observe that
the absolutely convergent fixed-point fast sweeping method (FE fast sweeping scheme) is the most efficient scheme among all three schemes studied here. By using the largest CFL number permitted by each method to converge to steady states, the absolutely convergent fixed-point fast sweeping method saves about $43\%$ CPU time of that by the TVD RK3 (RK Jacobi) scheme, while numerical steady states
obtained by these schemes are comparable.


\begin{table}
		\centering
		\begin{tabular}{|c|c|c|c|}\hline
			\multicolumn{4}{|c|}{FE Jacobi scheme }\\\hline
            $\gamma:$ CFL number & iteration number & final time & CPU time \\\hline
0.1   & 36330 & 21.09 & 2630.02 \\
\hline
 0.2  & Not convergent &- & - \\\hline
		\end{tabular}
\begin{tabular}{|c|c|c|c|}\hline
			\multicolumn{4}{|c|}{RK Jacobi scheme }\\\hline
            $\gamma:$ CFL number & iteration number & final time & CPU time \\\hline
 1.0     & 3702  & 22.76 & 731.28 \\
\hline
1.2   & 3078  & 22.70  & 603.77 \\
\hline
1.3   & 2838  & 22.67 & 556.39 \\
\hline
1.4   & 2637  & 22.68 & 518.45 \\
\hline
1.5   & Not convergent &   -    & - \\
    \hline
\end{tabular}
\begin{tabular}{|c|c|c|c|}\hline
			\multicolumn{4}{|c|}{FE fast sweeping scheme }\\\hline
            $\gamma:$ CFL number & iteration number & final time & CPU time \\\hline
         1.0     & 1048  & 19.27 & 396.86 \\
\hline
1.2   & 846   & 19.04 & 328.58 \\
\hline
1.3   & 788   & 18.78 & 297.14 \\
\hline
1.4   & Not convergent &    -   & - \\   \hline
		\end{tabular}
\caption{Example 12, supersonic flow past a long plate with an attack angle. Number of iterations, the final time and total CPU time when convergence is obtained. Convergence criterion threshold value  is $10^{-13}$. CPU time unit: second}
		\label{tab:Margin_settings12}
	\end{table}

\bigskip
\noindent{\bf Example 13. Supersonic flow past three long plates}

\noindent In the last example, we test the schemes by solving the case of a supersonic flow past three long plates with an attack angle of $\alpha=10 ^{\circ}$  .
The free stream Mach number is still $M_{\infty}=3$. The ideal gas goes from the left toward these three
long plates. The initial condition to start the iterations is $p=\frac{1}{\gamma' M_{\infty}^{2}}$, $\rho=1$, $u=\cos(\alpha)$ and $v=\sin(\alpha)$.
The computational field is $[0,5]\times[-5,5]$. The long plates are placed at $x \in [2,5]$ with $y=-2$, $x \in [2,5]$ with $y=0$, and $x \in [2,5]$ with $y=2$. The slip boundary condition is imposed on theses three long plates, and the physical values of the inflow and
outflow boundary
conditions are applied at the left, right, bottom, and top boundaries of the domain. In this example,
the plates also extend to the right boundary as Example 12, however more plates than Example 12 lead to more complicated interactions of the shocks and the rarefaction waves which pass through the right boundary on both sides of the plates.
The computational grid is $100 \times 200$.
The convergence criterion threshold value is set to be $10^{-13}$.

\begin{figure}
\centering
\subfigure[RK Jacobi scheme]{
\begin{minipage}[t]{0.5\linewidth}
\centering
\includegraphics[width=2.1in]{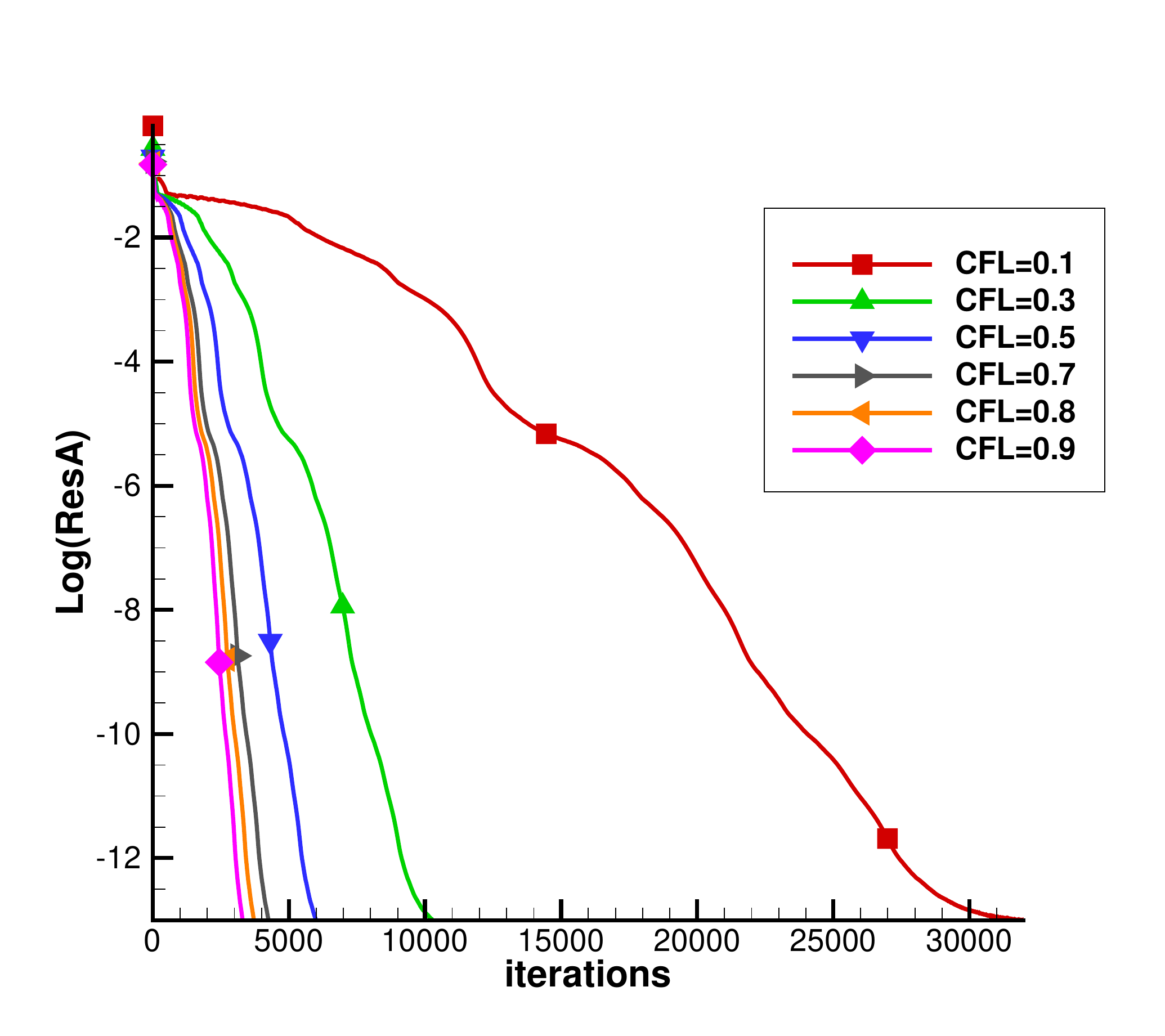}
\end{minipage}%
}%
\subfigure[FE fast sweeping scheme]{
\begin{minipage}[t]{0.5\linewidth}
\centering
\includegraphics[width=2.1in]{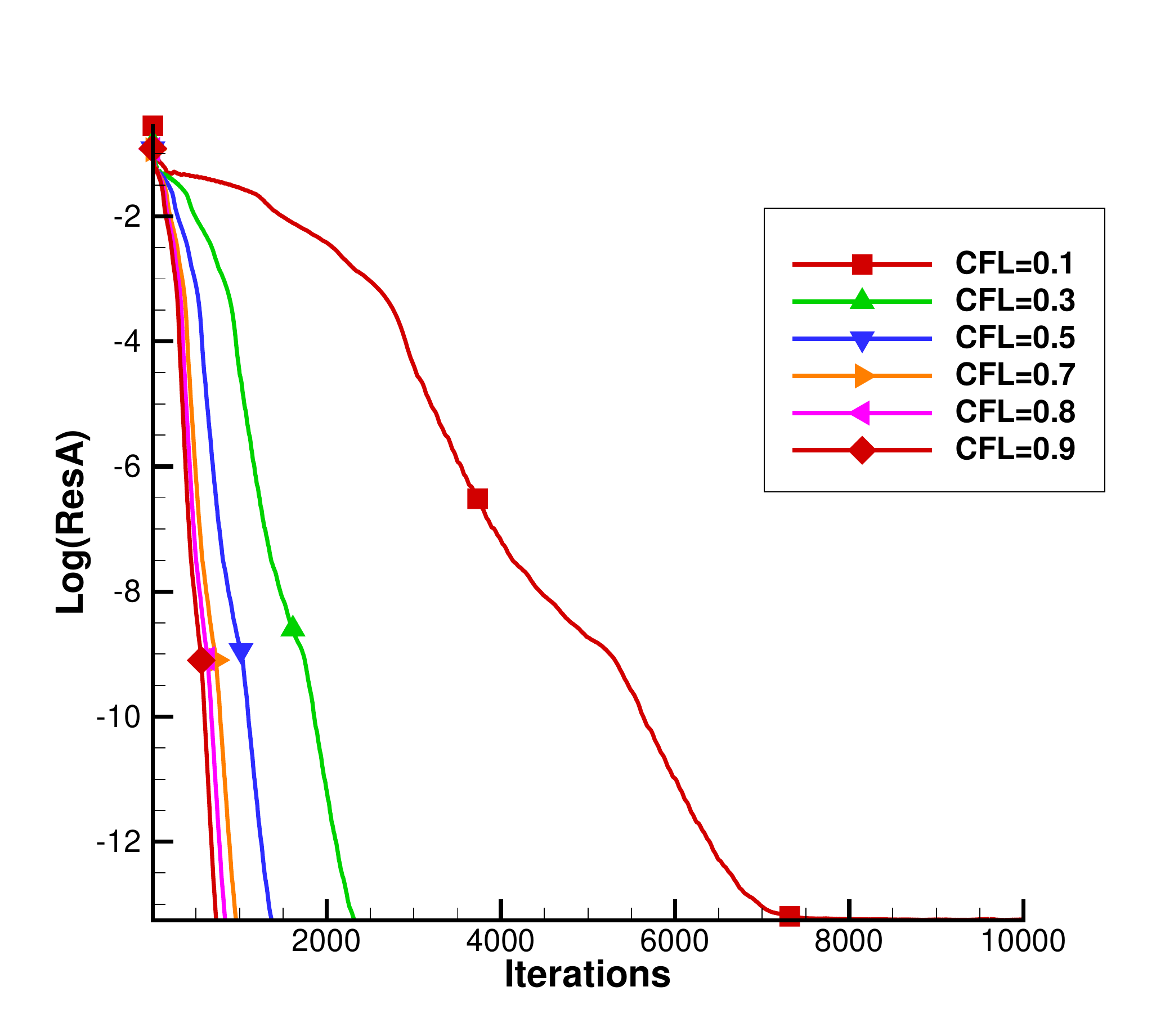}
\end{minipage}%
}%
\centering
\caption{Example 13, supersonic flow past three long plates with an attack angle. The convergence history of the residue
 as a function of number of iterations for two schemes with different CFL numbers.}
 \label{ex13fig1}
\end{figure}

\begin{figure}
\centering
\subfigure[RK Jacobi scheme]{
\begin{minipage}[t]{0.5\linewidth}
\centering
\includegraphics[width=2.1in]{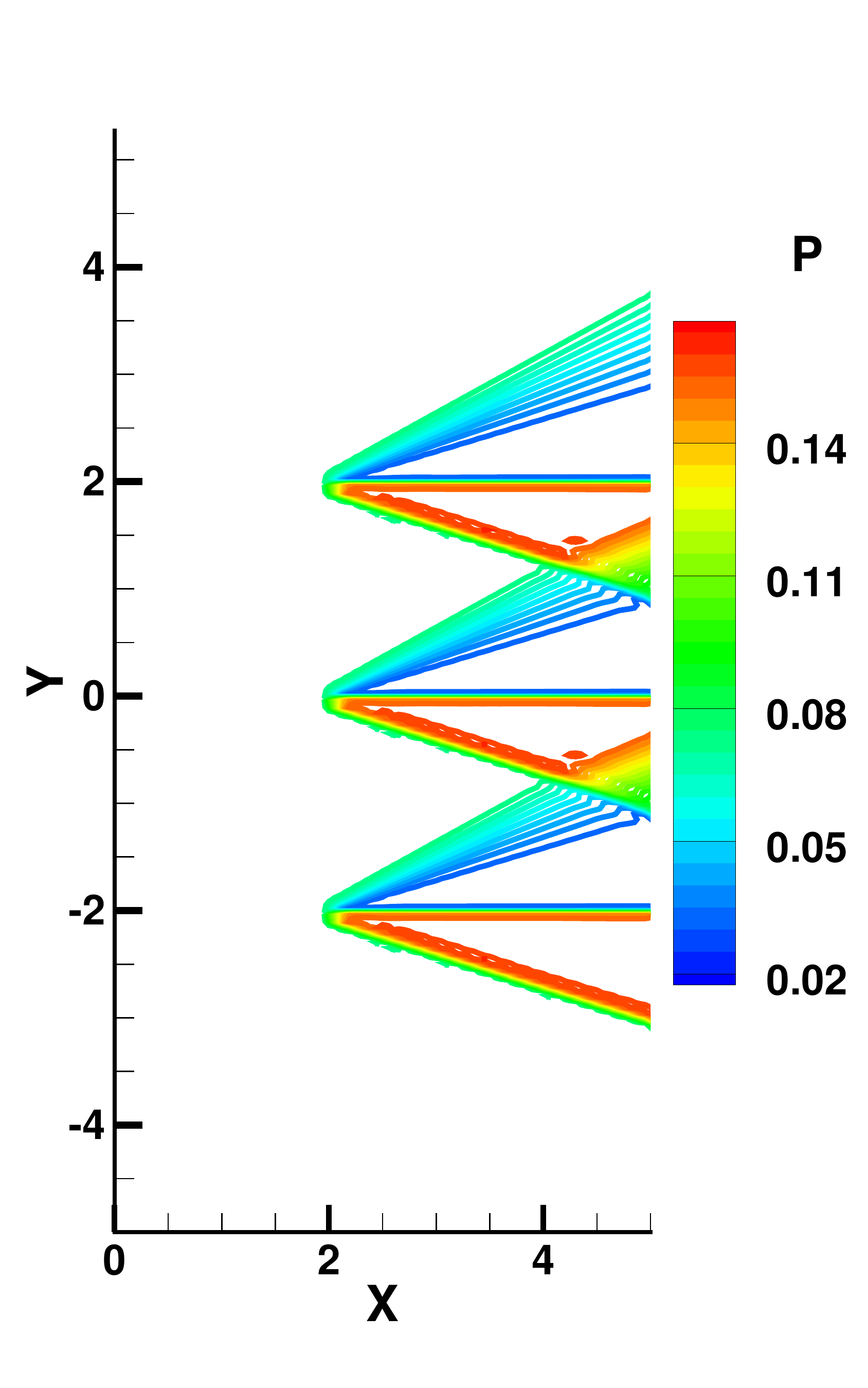}
\end{minipage}%
}%
\subfigure[FE fast sweeping scheme]{
\begin{minipage}[t]{0.5\linewidth}
\centering
\includegraphics[width=2.1in]{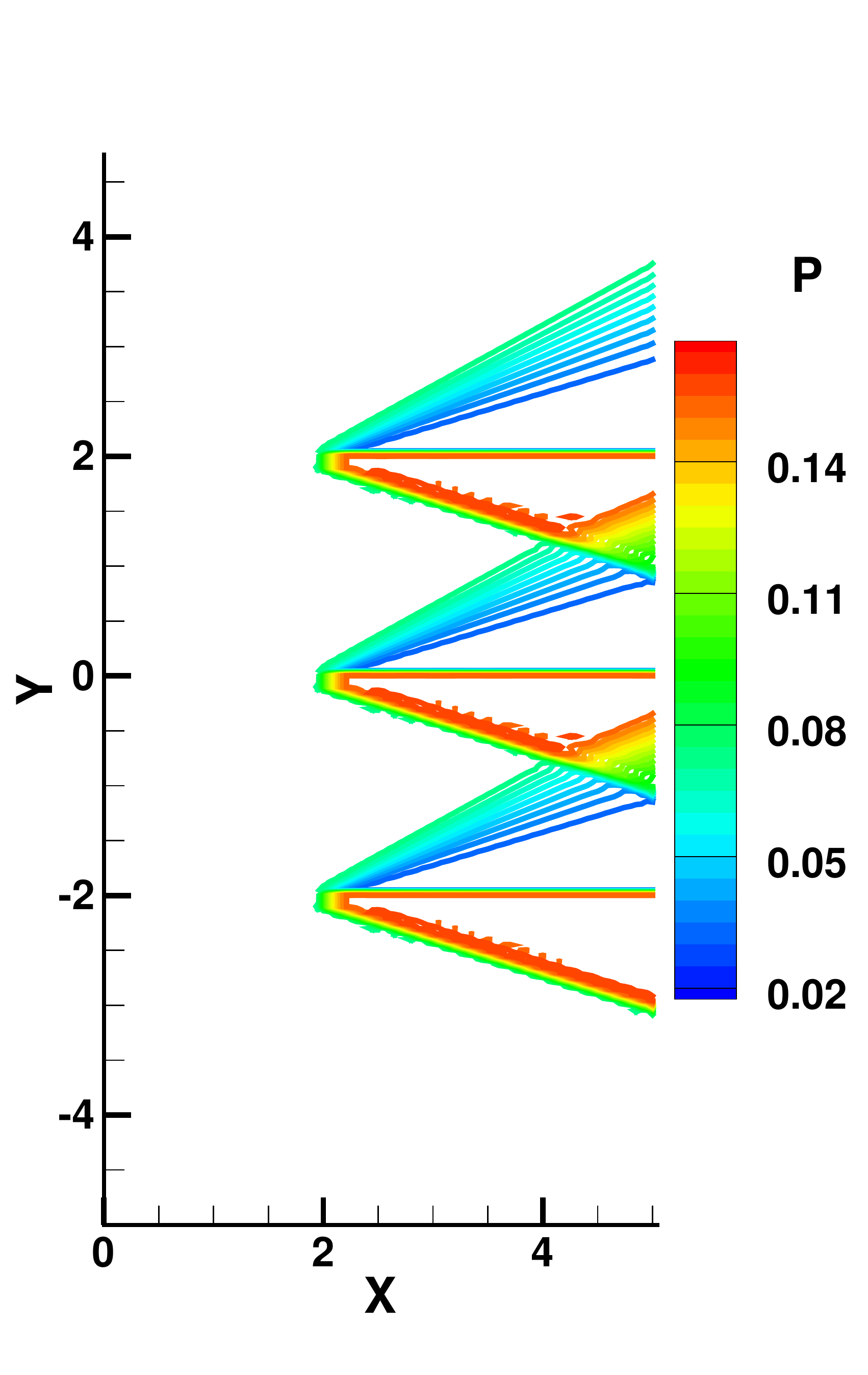}
\end{minipage}%
}%
\centering
\caption{Example 13, supersonic flow past three long plates with an attack angle. 30 equally pressure contour from 0.02 to 0.17 of the converged steady states of numerical solutions by two different iterative schemes.}
\label{ex13fig2}
\end{figure}

In Table \ref{tab:Margin_settings13}, number of iterations required
to reach the convergence threshold value $10^{-13}$, the final time and total CPU
time when convergence is obtained for these three iterative schemes with different
CFL numbers $\gamma$ are reported. Residue history in terms of iterations for the RK Jacobi and the FE fast sweeping schemes with different CFL numbers is
 presented in Figure \ref{ex13fig1}, and contour plots of the pressure variable of the converged steady state solutions of the RK Jacobi and the FE fast sweeping schemes are shown in Figure \ref{ex13fig2}. In this last example, we obtain the consistent conclusion with other examples, i.e.,
the proposed absolutely convergent fixed-point fast sweeping method (FE fast sweeping scheme) is the most efficient scheme among all three schemes studied here, in terms of both iteration numbers and CPU times to reach the convergence criterion. By using the largest CFL number permitted by each method to converge to steady states, the absolutely convergent fixed-point fast sweeping method saves more than $58\%$ CPU time of that by the TVD RK3 (RK Jacobi) scheme, while numerical steady states
obtained by these schemes are comparable.

\begin{table}
		\centering
		\begin{tabular}{|c|c|c|c|}\hline
			\multicolumn{4}{|c|}{FE Jacobi scheme }\\\hline
            $\gamma:$ CFL number & iteration number & final time & CPU time \\\hline
0.1	&10497	&14.26	&1481.52\\\hline
            0.2  & Not convergent & - & - \\\hline
		\end{tabular}
\begin{tabular}{|c|c|c|c|}\hline
			\multicolumn{4}{|c|}{RK Jacobi scheme }\\\hline
            $\gamma:$ CFL number & iteration number & final time & CPU time \\\hline
0.5&	6003	&13.59	&851.97\\\hline
0.7	&4257	&13.50	&603.83\\\hline
0.9	&3303	&13.46	&469.42\\\hline
1.0   & Not convergent &   -    &- \\
    \hline
\end{tabular}
\begin{tabular}{|c|c|c|c|}\hline
			\multicolumn{4}{|c|}{FE fast sweeping scheme }\\\hline
            $\gamma:$ CFL number & iteration number & final time & CPU time \\\hline
0.5	&1340	&12.33	&365.42\\\hline
0.8	&824	&12.12	&225.39\\\hline
0.9	&720	&11.91	&195.45\\\hline
1.0   & Not convergent &   -    & - \\    \hline
		\end{tabular}

		\caption{Example 13, supersonic flow past three long plates with an attack angle. Number of iterations, the final time and total CPU time when convergence is obtained. Convergence criterion threshold value  is $10^{-13}$. CPU time unit: second}
		\label{tab:Margin_settings13}
	\end{table}

\section{Concluding remarks}

In a recent work \cite{WuLiang}, the fast sweeping techniques were incorporated into a fifth order WENO method for
efficiently solving steady state problems of hyperbolic conservation laws. It was found that by using the
fast sweeping techniques, the forward Euler scheme with the fifth order WENO spatial discretization achieves a much larger CFL number than that in a regular time marching approach, hence it is much more efficient to converge to steady states of the high order WENO scheme. The forward Euler fast sweeping method is also more efficient than the popular TVD RK3 scheme to converge to steady states. However, an open problem in the fast sweeping WENO scheme in \cite{WuLiang} is that for some benchmark numerical examples, the iteration residue hangs at a truncation error level instead of converging to machine zero / round off errors. In this paper, we adopt the fifth order multi-resolution WENO scheme in \cite{JUNZ} and form a novel absolutely convergent forward Euler type fixed-point fast sweeping method for steady state of hyperbolic conservation laws. Extensive numerical experiments, including solving difficult problems which are challenging for high order schemes to converge to steady states, verify that the fast sweeping scheme can significantly enlarge the CFL number of the forward Euler scheme with a high order WENO spatial discretization (e.g., a fifth order WENO scheme) to the level of the TVD RK3 scheme. Furthermore, more than $40\%$ to $60\%$ CPU time is saved by using the forward Euler type fixed-point fast sweeping method rather than the TVD RK3 scheme to converge to steady states of the WENO scheme. The residue of the new absolutely convergent fast sweeping iterations converges to machine zero / round off errors for all benchmark problems tested in this paper. In the future, we will extend the proposed scheme to that on unstructured meshes for solving steady state problems on complex domains.


\begin{thebibliography}{99}

\bibitem{RB}
R. Borges, M. Carmona, B. Costa and W.S. Don, {\em An improved weighted essentially
non-oscillatory scheme for hyperbolic conservation laws}, J. Comput. Phys., 227 (2008),
3191-3211.

\bibitem{Capdeville2008}
G. Capdeville, {\em A central WENO scheme for solving hyperbolic conservation
laws on non-uniform meshes}, Journal of Computational Physics, 227 (2008), 2977-3014.

\bibitem{MC}
M. Castro, B. Costa and W.S. Don, {\em High order weighted essentially non-oscillatory
WENO-Z schemes for hyperbolic conservation laws}, J. Comput. Phys., 230 (2011),
1766-1792.

\bibitem{S.Chen}
S. Chen, {\em Fixed-point fast sweeping WENO methods for steady state solution of scalar hyperbolic
conservation laws}, International Journal of Numerical Analysis and Modeling, 11 (2014), 117-
130.

\bibitem{W.Chen}
W. Chen, C.-S. Chou and C.-Y. Kao, {\em Lax-Friedrichs fast sweeping methods for steady state problems for hyperbolic conservation laws}, Journal of Computational Physics, 234 (2012), 452-471.

\bibitem{M.Dumbser}
M. Dumbser and M. K\"aser, {\em Arbitrary high order non-oscillatory finite volume schemes
on unstructured meshes for linear hyperbolic systems}, J. Comput. Phys., 221 (2007),
693-723.

\bibitem{FLZ}
S. Fomel, S. Luo and H. Zhao,
{\em Fast sweeping method for the factored eikonal equation},
Journal of Computational Physics, 228 (2009), 6440--6455.

\bibitem{GST}
S. Gottlieb, C.-W. Shu and E. Tadmor, {\em Strong stability-preserving high-order time discretization methods}, SIAM Review, 43 (2001), 89--112.

\bibitem{HHSSXZ}
W. Hao, J.D. Hauenstein, C.-W. Shu, A.J. Sommese, Z. Xu and Y.-T. Zhang, {\em A homotopy method based on WENO schemes for solving steady state problems of hyperbolic conservation laws}, Journal of Computational Physics, 250, (2013), 332--346.

\bibitem{C.Hu}
C. Hu and C.-W. Shu, {\em Weighted essentially non-oscillatory schemes on triangular meshes}, Journal
of Computational Physics, 150 (1999), 97-127.

\bibitem{JiangShu}
G.-S. Jiang and C.-W. Shu, {\em Efficient implementation of weighted ENO schemes}, Journal of Computational Physics, 126 (1996), 202-228.

\bibitem{D.Levy}
D. Levy, S. Nayak, C.-W. Shu and Y.-T. Zhang, {\em Central WENO schemes for Hamilton-Jacobi
equations on triangular meshes}, SIAM Journal on Scientific Computing, 28 (2006), 2229-2247.

\bibitem{Levy1999}
D. Levy, G. Puppo and G. Russo, {\em Central WENO schemes for hyperbolic
systems of conservation laws}, Mathematical Modelling and Numerical Analysis, 33 (1999), 547-571.

\bibitem{D. L}
D. Levy, G. Puppo and G. Russo, {\em Compact central WENO schemes for multidimensional
conservation laws}, SIAM J. Sci. Comput., 22 (2000), 656-672.

\bibitem{LSZZ}
F. Li, C.-W. Shu, Y.-T. Zhang and H.-K. Zhao,
{\em A second order discontinuous Galerkin fast sweeping method for Eikonal equations},
Journal of Computational Physics, 227 (2008), 8191-8208.

\bibitem{LOC}
X.-D. Liu, S. Osher and T. Chan, {\em Weighted essentially non-oscillatory schemes}, Journal of Computational Physics, 115 (1994), 200-212.

\bibitem{Y.Liu}
Y. Liu and Y.-T. Zhang, {\em A robust reconstruction for unstructured WENO schemes}, Journal of
Scientific Computing, 54 (2013), 603-621.

\bibitem{QZZ1}
J. Qian, Y.-T. Zhang and H.-K. Zhao,
{\em Fast  sweeping  methods for Eikonal equations on triangular meshes},
SIAM Journal on Numerical Analysis, 45 (2007), 83--107.

\bibitem{QZZ2}
J. Qian, Y.-T. Zhang and H.-K. Zhao,
{\em A fast sweeping method for static convex Hamilton-Jacobi equations},
Journal of Scientific Computing, 31 (2007), 237--271.

\bibitem{M.A.Sa}
 M.A. Saad, Compressible Fluid Flow. Prentice Hall, New York, (1993).

\bibitem{CWS2}
  C.-W. Shu, Essentially non-oscillatory and weighted essentially non-oscillatory schemes
for hyperbolic conservation laws, in Advanced Numerical Approximation of Nonlinear
Hyperbolic Equations, B. Cockburn, C. Johnson, C.-W. Shu and E. Tadmor (Editor:
A. Quarteroni), Lecture Notes in Mathematics, volume 1697, Springer, Berlin, 1998,
325-432.

\bibitem{C.W.S}
 C.-W. Shu and S. Osher, {\em Efficient implementation of essentially non-oscillatory shock capturing
schemes}, Journal of Computational Physics, 77 (1988), 439-471.

\bibitem{WZ}
L. Wu and Y.-T. Zhang, {\em A third order fast sweeping method with linear computational complexity for Eikonal equations}, Journal of Scientific Computing, 62 (2015), 198-229.


\bibitem{WuLiang}
L. Wu, Y.-T. Zhang, S. Zhang and C.-W. Shu, {\em High order fixed-point sweeping WENO methods for steady state of hyperbolic conservation laws and its convergence study}, Communications in Computational Physics, 20 (2016), 835-869.

\bibitem{TMY}
T. Xiong, M. Zhang, Y.-T. Zhang and C.-W. Shu, {\em Fast sweeping fifth order WENO scheme for
static Hamilton-Jacobi equations with accurate boundary treatment}, Journal of Scientific Computing, 45 (2010), 514-536.

\bibitem{SSCW}
S. Zhang, S. Jiang and C.-W. Shu, {\em Improvement of convergence to steady state solutions of Euler
equations with the WENO schemes}, Journal of Scientific Computing, 47 (2011), 216-238.

\bibitem{SCW}
S. Zhang and C.-W. Shu, {\em A new smoothness indicator for the WENO schemes and its effect on the
convergence to steady state solutions}, Journal of Scientific Computing, 31 (2007), 273-305.

\bibitem{ZCLZS}
Y.-T. Zhang, S. Chen, F. Li, H. Zhao and C.-W. Shu,
{\em Uniformly accurate discontinuous Galerkin fast sweeping methods for Eikonal equations}, SIAM Journal on Scientific Computing, 33 (2011), 1873-1896.

\bibitem{Y.-T.Zhang}
Y.-T. Zhang and C.-W. Shu, {\em High order WENO schemes for Hamilton-Jacobi equations on triangular meshes}, SIAM Journal on Scientific Computing, 24 (2003), 1005-1030.

\bibitem{Y.-T.Zhangand C.-W.Shu}
Y.-T. Zhang and C.-W. Shu, {\em Third order WENO schemes on three dimensional tetrahedral meshes},
Communications in Computational Physics, 5 (2009), 836-848.

\bibitem{Y.-T.Zhang H.-K}
Y.-T. Zhang, H.-K. Zhao and S. Chen, {\em Fixed-point iterative sweeping methods for static Hamilton-Jacobi equations}, Methods and Applications of Analysis, 13 (2006), 299-320.

\bibitem{ZZQ}
Y.-T. Zhang, H.-K. Zhao and J. Qian, {\em High order fast sweeping methods for static Hamilton-Jacobi
equations}, Journal of Scientific Computing, 29 (2006), 25-56.

\bibitem{Z}
H.-K. Zhao,
{\em A fast sweeping method for Eikonal equations},
Math. Comp., 74 (2005), 603-627.

\bibitem{X. Zhong}
X. Zhong and C.-W. Shu, {\em A simple weighted essentially nonoscillatory limiter for Runge-
Kutta discontinuous Galerkin methods}, J. Comput. Phys., 232 (2013), 397-415.

\bibitem{JUN}
J. Zhu and J. Qiu, {\em A new type of finite volume WENO schemes for hyperbolic
conservation laws}, Journal of Scientific Computing, 73 (2017), 1338-1359.

\bibitem{JUNZ4}
J. Zhu and C.-W. Shu, {\em Numerical study on the convergence to steady state solutions of a new class of high order WENO schemes}, Journal of Computational Physics, 349 (2017), 80-96.

\bibitem{JUNZ}
J. Zhu and C.-W. Shu, {\em A new type of multi-resolution WENO schemes with increasingly higher order of accuracy}, Journal of Computational Physics, 375 (2018), 659-683.

\bibitem{JUNZ2}
J. Zhu and C.-W. Shu, {\em A new type of third-order finite volume multi-resolution WENO schemes on tetrahedral meshes}, Journal of Computational Physics, v406 (2020), 109212.

\bibitem{JUNZ3}
J. Zhu and C.-W. Shu, {\em Convergence to steady-state solutions of the new type of high-order multi-resolution WENO schemes: a numerical study}, Communications on Applied Mathematics and Computation, v2 (2020), 429-460.


\end{thebibliography}
\end{document}